\documentclass[10pt,letterpaper]{amsart} 
\usepackage{xcolor}
\usepackage[letterpaper, margin=1in]{geometry}
\usepackage{subcaption}
\usepackage{multicol}
\usepackage{booktabs}
\usepackage{listings}
\usepackage{amsmath, amsthm, amssymb, wasysym, verbatim, bbm, color, graphicx}
\usepackage{url}
\usepackage{mathtools}
\usepackage{thmtools}
\usepackage[colorlinks,linkcolor=blue]{hyperref}
\usepackage{enumerate}
\usepackage{xfrac}
\usepackage{float}
\usepackage{todonotes}
\usepackage{accents}
\usepackage{caption}
\usepackage{subcaption}

\newcommand{\doublehat}[1]{\bar{\bar{#1}}}

\newcommand{\R}{\mathbb{R}}

\newcommand{\N}{\mathbb{N}}
\newcommand{\Q}{\mathbb{Q}}

\newcommand{\tu}{\widetilde{u}}
\newcommand{\hu}{\widehat{u}}

\newcommand{\hM}{\overline{M}}

\newcommand{\hR}{\bar{R}}
\newcommand{\tR}{\widetilde{R}}

\DeclareMathOperator{\E}{\mathbb{E}}
\DeclareMathOperator{\tr}{\text{tr}}

\newcommand\norm[1]{\left\lVert#1\right\rVert}
\newcommand{\Grad}{\nabla}
\newcommand{\Div}{\operatorname{div}}
\newcommand{\dom}{D}
\newcommand{\Vh}{\mathbb{V}_h}
\newcommand{\Ph}{\mathbb{Q}_h}

\newcommand{\Ih}{\mathcal{I}_h}

\newcommand{\Uh}{\mathbb{U}_h}
\newcommand{\Sh}{\mathcal{S}_h}

\newcommand{\PUh}{\mathcal{P}_{\Uh}}
\newcommand{\PVh}{\mathcal{P}_{\Vh}}

\newcommand{\Ldiv}{L^2_{\Div}(\dom)}

\newcommand{\PP}{\mathbb{P}}

\newtheorem{lemma}{Lemma}[section]
\newtheorem{proposition}[lemma]{Proposition}

\newtheorem{theorem}[lemma]{Theorem}
\newtheorem{remark}[lemma]{Remark}

\newtheorem{assumption}[lemma]{Assumption}

\newtheorem*{maintheorem*}{Main Theorem}
\theoremstyle{definition}{\newtheorem{definition}[lemma]{Definition}}

\allowdisplaybreaks

\numberwithin{equation}{section}

\title[Numerics for the stochastic Navier--Stokes equations]{Convergence of a finite element Chorin projection method for the stochastic Navier--Stokes equations}

\author[F. Weber]{Franziska Weber}
\address[Franziska Weber]{\newline Department of Mathematics \newline University of California, Berkeley \newline Berkeley, CA 94720, USA.}
\email[]{fweber@berkeley.edu}
\date{\today}
\thanks{F.W. was partially supported by NSF DMS 2438083.}
\subjclass[2020]{65M12, 65C30, 65M60, 60H35, 35Q30}
\keywords{Stochastic Navier--Stokes equations, Chorin projection method, finite elements, convergence analysis}

\begin{document}
	\begin{abstract}
		We analyze a fully discrete finite element projection method for the stochastic incompressible Navier--Stokes equations with homogeneous Dirichlet boundary conditions and multiplicative noise. We show solvability and uniform moment estimates for the intermediate and the corrected velocity approximations under the assumption that the initial data is square-integrable and divergence free, and that the deterministic forcing term is bounded in $L^{q}(0,T;H^{-1}(\dom))$ for a sufficiently large $q$. Under standard stability and approximation assumptions on the finite element spaces and a coupling condition on the spatial and temporal discretization parameters, we prove convergence of a subsequence to a martingale solution of the stochastic Navier--Stokes equations. The analysis combines the complementary spatial and temporal estimates of the predicted and the corrected velocity approximations and uses the numerical dissipation to identify their common limit. We establish tightness of the piecewise constant temporal interpolations of the velocity field approximations using Aldous' criterion combined with the uniform discrete moment estimates. The convergence result holds in two and three space dimensions as it does not rely on any additional regularity of the limiting martingale solution beyond the one naturally available from the energy inequality. Finally, we demonstrate the performance of the scheme in numerical experiments in two space dimensions for smooth and vortex sheet initial data.
		\end{abstract}
\maketitle

\section{Introduction}

We consider the numerical approximation of the stochastic incompressible Navier--Stokes equations subjected to a random force 
\begin{subequations}\label{eq:SNSm}
	\begin{align}
		du + ((u\cdot\Grad) u +\Grad p-\nu\Delta u-f )dt &= G(u)dW, \\
		\Div u & = 0,
	\end{align}
\end{subequations}
in $(t,x,\omega)\in [0,T]\times\dom\times\Omega$. $T>0$ is a finite time horizon, $\nu>0$ is the viscosity and $W=\{W(t)\}_{t\in [0,T]}$ a $Q$-Wiener process. $u:[0,T]\times\dom\times\Omega\to \R^d$ and $p: [0,T]\times\dom\times\Omega\to \R$ are stochastic processes representing the velocity field and the pressure respectively. $f$ is a deterministic force for which we assume $f\in L^{2^{r_0}}(0,T;H^{-1}(\dom))$ for some $r_0\geq 3$. $G(u)$ is a Lipschitz continuous mapping and takes values in
$\mathcal L(\mathcal{K},L^2(\dom))$, where $\mathcal{K}$ is a suitable Hilbert space, the specific requirements will be given below. Here $\dom\subset \R^d$, $d=2,3$ is a bounded Lipschitz polygonal ($d=2$) or polyhedral ($d=3$) domain (cf. Assumption~\ref{ass:stokes} in Section~\ref{sec:spatial}).
We complement this system with initial and homogeneous Dirichlet boundary conditions
\begin{equation}
	\label{eq:initdata}
	u(0,\cdot) = u_0,\quad \text{in }\, \dom,\quad u=0,\quad \text{on }\, (0,T)\times\partial\dom\times\Omega,
\end{equation}
where $u_0\in \Ldiv$, the space of square-integrable functions that are divergence free almost everywhere.
Here we denoted $\Div u =\sum_{j=1}^d \partial_j u^{(j)}$ and $(u\cdot\Grad)v= \sum_{j=1}^d u^{(j)}\partial_j v$. We will also denote the gradient of a vector field $u$ by $\Grad u = (\partial_j u^{(i)})_{i,j=1}^d$.
The stochastic Navier--Stokes equations were first studied by Bensoussan and Temam~\cite{Bensoussan1973}, where existence of solutions for the stochastic Navier--Stokes equations with random initial data and random forcing, including white noise, was shown. Multiplicative noise was first studied by~\cite{Brzezniak1992} where an existence and uniqueness result for small initial data was shown. Flandoli and Gatarek~\cite{Flandoli1995} established the existence of global martingale solutions and stationary solutions under general assumptions on the diffusion coefficients using a priori estimates to deduce tightness of the Galerkin approximations. Well-posedness in 2D for additive noise was shown in~\cite{Menaldi2002}, and for multiplicative noise in~\cite{Liu2010}. Stronger well-posedness and regularity results are also available under more restrictive assumptions on the initial data and noise, e.g.,~\cite{GlattHoltz2009,Mikulevicius2009}. Global regularity and well-posedness of strong solutions in 3D remain open.

This is the starting point for the numerical approximation of~\eqref{eq:SNSm} we want to take in this article: Given initial data in $L^2_{\Div}(\dom)$, can the convergence to a martingale solution of~\eqref{eq:SNSm} be shown under no additional assumptions on the solution of~\eqref{eq:SNSm} that cannot be proved in the continuous case? We want to answer this question in the specific case of a fractional step or projection method as originally introduced for the deterministic system by Chorin and Temam~\cite{Chorin1967,Chorin1968,Temam1968}. From a numerical point of view, one of the key difficulties in approximating solutions to the incompressible Navier--Stokes equations is the combination of the nonlinear advection with the divergence constraint of the velocity field. Ensuring that the discrete velocity field satisfies the divergence-free condition, or an approximate version thereof, is crucial for maintaining stability of the numerical method. In Chorin's projection method, the evolution of the velocity due to the nonlinearity, dissipation and forcing terms is decoupled from the pressure computation. At each time step, one first solves a linearized convection-diffusion problem resulting in an intermediate velocity $\tu^{m+1}$ that does not necessarily satisfy the incompressibility constraint, followed by a projection step that enforces the (approximate) divergence-free constraint and yields the corrected velocity $u^{m+1}$. 
The scheme considered here has this prediction-projection structure. In semi-discrete form, the prediction step reads
\begin{equation}\label{eq:predictionstep}
	\tu^{m+1} - u^m + \Delta t\left[(\tu^m\cdot\Grad) \tu^{m+1}+\frac{1}{2}(\Div \tu^m) \tu^{m+1}\right] = \Delta t\nu \Delta \tu^{m+1} +\Delta t f^m+ G(u^m) \Delta W_m,
\end{equation}
where we denoted
\begin{equation*}
	\Delta W_m = W(t^{m+1}) - W(t^m),
\end{equation*}
followed by a projection step
\begin{equation}\label{eq:projectionstep}
	\frac{u^{m+1}-\tu^{m+1}}{\Delta t }= -\Grad p^{m+1},\quad 	\Div u^{m+1}  = 0.
\end{equation}	 
Projection methods have been popular both in practical computations and in theoretical studies  due to their simplicity and modularity. They have been studied extensively in the semi-discrete and the fully discrete setting, including convergence rates under the assumption that the solution is smooth, e.g.,~\cite{Temam1968,Temam1977,Temam1969,Guermond2006,Guermond1996,Guermond1998,Hou1992,Hou1993,Achdou2000,E2002,Shen1996}. More recently, they were also studied  under no additional uniqueness and regularity assumptions except those resulting from the energy inequality, and convergence up to a subsequence to a Leray--Hopf weak solution (whose existence is guaranteed by the work of Leray~\cite{Leray1934} and Hopf~\cite{Hopf1951}) was shown for finite difference and finite element methods~\cite{Kuroki2020,Maeda2022,Gallouet2023,Eymard2024,chorinprojection,Weber2026}. Recently, a computer-assisted proof of nonuniqueness of Leray--Hopf solutions in 3D has been announced in~\cite{Hou2025} and it is unclear whether global smooth solutions of the Navier--Stokes equations exist, convergence to a (distributional) Leray--Hopf solution is therefore a natural benchmark even in the deterministic setting.  The stochastic equivalent to Leray--Hopf solutions are martingale solutions (see Definition~\ref{def:martingalesolution} below). Numerical methods based on finite element spatial discretization whose approximations converge to a martingale solution of the stochastic Navier--Stokes equations were designed in~\cite{Brzezniak2013,Ondrejat2023}. Both approaches are based on monolithic discretizations in which velocity and pressure are solved for simultaneously (i.e., no decoupling as in the projection method~\eqref{eq:predictionstep} - \eqref{eq:projectionstep}). Both methods are fully discrete and the results do not require a priori regularity of the limiting solution. To the best of our knowledge, these are the only works for fully discrete schemes for the stochastic Navier--Stokes equations which establish convergence to martingale solutions without assuming additional regularity of the limiting solution. There are additional results for the 2D setting, and the 3D setting under the assumption that the limiting solution is sufficiently regular, where convergence rates and error estimates are proved~\cite{Carelli2012,Breit2021,Breit2024,Breit2024b,Bessaih2022,Bessaih2019}. A projection method time discretization similar to the one considered here was studied in~\cite{Hausenblas2019} for the 2D stochastic Navier--Stokes equations with additive solenoidal noise and periodic boundary conditions. 
They established convergence rates in probability of almost $1/4$ for both velocity and pressure, as well as strong convergence without an explicit rate. Their analysis is restricted to the time discretization and relies on the regularity of the unique strong solution in 2D and on the solenoidal character of the noise. 
In contrast, we here consider a fully discrete scheme for the 2 or 3D Navier--Stokes equations in a bounded domain and without any smoothness assumptions on the solution beyond those that hold for martingale solutions. More works are available for the stochastic Stokes equations and convergence rates have been established for different versions of the projection method~\cite{Feng2023,Carelli2012b}. Alternative discretizations based on splitting methods or artificial compressibility were studied in~\cite{Bessaih2014,Doghman2024}. For an overview of the analysis of numerical methods for the stochastic Stokes and Navier--Stokes equations, see also~\cite{Breit2025}.

The main contribution of this work is the convergence analysis of a fully discrete finite element scheme of Chorin's projection method \eqref{eq:predictionstep}--\eqref{eq:projectionstep} in 2 and 3D. Specifically, we prove
\begin{theorem}\label{thm:mainthm}
	Assume that $u_0\in \Ldiv$, $f\in L^{2^{r_0}}(0,T;H^{-1}(\dom))$ with $r_0\geq 3$ and where $\dom\subset \R^d$, $d=2,3$ is a bounded polygonal or polyhedral domain. Assume that $G$ in the noise term of~\eqref{eq:SNSm} is Lipschitz continuous. 
	
	Then up to a subsequence, the laws of the approximations of the intermediate and the corrected velocity field computed by an inf-sup stable finite element discretization of~\eqref{eq:predictionstep} - \eqref{eq:projectionstep} converge to the law of the same weak martingale solution of~\eqref{eq:SNSm} - \eqref{eq:initdata} as the time step $\Delta t\to 0$ and the mesh size $h\to 0$ under the condition that $h^{2(k+s)}=o(\Delta t)$, where $k$ is the approximation order of the finite element space of the velocity field and $s\in (0,1]$ is the elliptic regularity index of the Stokes problem. On a suitable probability space identified with Skorokhod's theorem, both sequences of approximations converge strongly in $L^2([0,T]\times\dom)$.
\end{theorem}
The precise assumptions are stated in Section~\ref{sec:spatial} and the proof follows from Proposition~\ref{prop:skorokhod}, Proposition~\ref{prop:tightness} and Theorem~\ref{lem:martingaleidentification} in Section~\ref{sec:conv}.
We note that the multiplicative noise is not required to be solenoidal. The coupling condition $h^{2(k+s)}=o(\Delta t)$ can be dropped for velocity-pressure pairs that admit a suitable Fortin operator, cf. Remark~\ref{rem:fortin}.

A central difficulty in the analysis compared with monolithic discretizations is that the projection method produces two discrete velocity fields with complementary properties. The intermediate velocity field $\tu^{m+1}_h$ satisfies $\tu_h\in L^{2^{r_0}}(\Omega;L^2(0,T;H^1_0(\dom))),$ uniformly in the discretization parameters, but is not discretely divergence free and does not admit a suitable direct estimate on its time increments. On the other hand, the corrected velocity $u_h$ is discretely divergence free and admits the temporal estimates needed to establish tightness, while the available energy bounds only yield $u_h \in L^{2^{r_0}}(\Omega;L^\infty(0,T;L^2(\dom)))$ uniformly in the discretization parameters.
Thus, the estimates required for compactness and tightness are distributed between the two velocity approximations. In the deterministic case, these can be combined pathwise to obtain precompactness, but in the stochastic case here the bounds only hold in expectation. We therefore first use Aldous' criterion to establish tightness of the laws of the piecewise constant time interpolants of the projected velocity approximations in the Skorokhod space $\mathbb{D}([0,T];H^{-1}(\dom))$ and use the numerical dissipation of the scheme to show that the difference of the two approximation sequences $\tu_h - u_h$ goes to zero in probability. We then recover the strong convergence of both approximation sequences in $L^2([0,T]\times\dom)$ after transferring to a new probability space using the Skorokhod representation theorem. In addition, the pressure generated in the projection step needs to be controlled. The energy estimate does not supply bounds on the pressure variable directly, therefore we use the numerical dissipation of the scheme to handle those terms. This is also where the condition $h^{2(k+s)}=o(\Delta t)$ enters.

Another difference compared to previous works concerns the time interpolation that is used to obtain tightness. In~\cite{Brzezniak2013,Ondrejat2023} tightness is proved via piecewise linear interpolations of the velocity field approximation by using bounds in fractional Sobolev spaces~\cite{Brzezniak2013} and by using H\"older bounds through the Kolmogorov continuity theorem~\cite{Ondrejat2023}, which are combined with the Aubin--Lions and Prokhorov theorems. Tightness is then transferred to the piecewise constant interpolations. The piecewise linear interpolants are not adapted to the filtration. In contrast, we verify Aldous' criterion directly for the piecewise constant interpolants of the projected velocity approximation, which are adapted to the filtration, so that the stochastic increments can be estimated at stopping times. Then we show separately that the limit is continuous in time almost surely using a result by Ethier and Kurtz~\cite[Theorem 3.10.2]{Ethier1986}. While the numerical dissipation in the scheme may allow comparing the piecewise constant and piecewise linear interpolants, our argument does not require it.  
Our proof also drew inspiration from~\cite{Motyl2013,Berselli2021}, and the identification of the quadratic variation follows~\cite{Brzezniak2013}.

The analysis is complemented with numerical experiments in 2D. We check the empirical spatial and temporal convergence rates for smooth initial data and then illustrate the influence of the multiplicative noise in an experiment with interacting vortices. Finally, we test the performance of the scheme on rough vortex sheet initial data which is allowed by the stability and convergence analysis. Together, these experiments illustrate the behavior of the scheme for initial data of different regularity and the variability induced by the stochastic forcing.

We hope that parts of our approach could extend beyond the present setting. Since the tightness argument does not require the approximations to be continuous in time, the approach could work for noise with jumps, e.g., L\'evy noise, where the limit itself is only a c\`adl\`ag process. Similarly, in schemes with less numerical dissipation, such as incremental projection methods, it may not be possible to show that the piecewise constant and piecewise linear temporal interpolations have the same limit, and therefore one might similarly have to work with piecewise constant interpolations. However, the occurrence of the pressure variable in the prediction step might lead to additional challenges, therefore we leave this for future work.
\subsection{Outline}
The rest of this article is organized as follows: In Section~\ref{sec:prelim}, we introduce the necessary notation, technical assumptions and the definition of martingale solutions. In Section~\ref{sec:num}, we introduce the numerical scheme and prove its solvability and stability estimate. Then in Section~\ref{sec:conv}, we prove the convergence to a martingale solution and finally, we conclude with numerical experiments in Section~\ref{sec:numexp}. In the appendix, we prove an auxiliary technical result. 
\section{Preliminaries}\label{sec:prelim}

We start by introducing general notation.
\subsection{Notation}
We denote the norm of a Banach space $X$ as $\|\cdot\|_X$ and its dual space by $X^*$.  By $C([0,T];X_w)$ we mean functions that are weakly continuous in time with values in $X$, i.e., $u\in C([0,T];X_w)$ if
the mapping
\begin{equation*}
	t\mapsto \langle u(t),v\rangle_{X,X^*}
\end{equation*}
is continuous for any $v\in X^*$.
The inner product on $L^2$ will be denoted by $(\cdot, \cdot)$. We will use $\langle\cdot,\cdot\rangle_{X,X^*} = \langle\cdot,\cdot\rangle$ to denote the duality bracket between $X$ and $X^*$ for a Banach space $X$. Sometimes we will omit the subscripts $X,X^*$.

We will use the subscript $\Div$  to indicate the divergence-free vector spaces, for example, 
\begin{equation*}
	\label{eq:sigma_space}
	\begin{aligned}
		C_{c,\Div}^\infty(\dom)=\{{\phi}\in C_c^\infty(\dom); \Div {\phi}=0\},&\quad L^2_{\Div}(\dom)=\{ {\phi}\in L^2(\dom):\Div  {\phi}=0,  {\phi}\cdot {n}|_{\partial\dom}=0\}=\overline{C_{c,\Div}^\infty(\dom)}^{L^2(\dom)},\\
		&\quad H^1_{ \Div}(\dom)=H^1_0(\dom)\cap L^2_{\Div}(\dom),
	\end{aligned}
\end{equation*}
where we used $n$ to denote the outward normal vector of the domain $\dom$ and we will use $H^{-1}(\dom):=(H^1_0(\dom))^*$, the dual space of $H^1_0(\dom)$, frequently. 
Additionally,
\begin{equation*}
	L^2_0(\dom)= \left\{\phi\in L^2(\dom); \, \int_{\dom} \phi dx =0\right\}
\end{equation*}
will be used to denote the space of $L^2$-functions that have zero average over the domain.
\subsection{Stochastic setting and assumptions on the data}
We let $\mathcal{K}$ be a separable real Hilbert space and
$Q:\mathcal K\to\mathcal K$ be a nonnegative, self-adjoint,
trace-class operator. We fix a stochastic basis
$(\Omega,\mathcal{G},(\mathcal{G}_t)_{t\in[0,T]},\mathbb P)$
satisfying the usual assumptions: The filtration is
right-continuous and $\mathcal G_0$ contains all the
$\mathbb P$-null sets of $\mathcal G$.
We denote the expectation with respect to the probability measure $\mathbb P$ by $\E$.
Next, we let $W$ be a $\mathcal{K}$-valued $Q$-Wiener process adapted to
$(\mathcal G_t)_{t\in[0,T]}$. Thus, $W(0)=0$ almost surely,
$W$ has continuous paths, and, for $0\le s<t\le T$,
the increment $W(t)-W(s)$ is Gaussian with mean zero and covariance $(t-s)Q$ and independent of $\mathcal G_s$.
For Hilbert spaces $E$ and $F$, we will use $\mathcal L(E,F)$ and $\mathcal L_2(E,F)$ to denote the spaces of bounded linear and Hilbert--Schmidt operators, respectively.
$A^*$ will be used to denote the adjoint of an operator $A$.
In particular, we have for $A\in\mathcal L(\mathcal K,L^2(\dom))$,
\begin{equation*}
	\norm{A Q^{1/2}}^2_{\mathcal{L}_2(\mathcal{K},L^2(\dom))}\leq \tr(Q)\norm{A}^2_{\mathcal{L}(\mathcal{K},L^2(\dom))}.
\end{equation*}
For an $E$-valued random variable $X$, we denote its law by $\mathrm{law}(X)=\mathbb P\circ X^{-1}$. By `convergence in law' we mean the weak convergence of the corresponding probability measures.

For the rest of this paper, we will assume that the initial velocity field $u_0\in \Ldiv$ and the body force $f\in L^{2^{r_0}}(0,T;H^{-1}(\dom))$ are deterministic.
We assume that the noise coefficient $G:L^2(\dom)\to \mathcal{L}(\mathcal{K},L^2(\dom))$ satisfies 
\begin{align}
	\norm{G(v)}_{\mathcal L(\mathcal K,L^2(\dom))}	&\le K_1\left(1+\norm{v}_{L^2(\dom)}\right),\label{eq:lingrowthG}\\
	\norm{G(v)-G(w)}_{\mathcal L(\mathcal K,L^2(\dom))}&\le K_2\norm{v-w}_{L^2(\dom)},	\label{eq:Glipschitz}
\end{align}
for all $v,w\in L^2(\dom)$, where $K_1,K_2>0$ are fixed constants.
 $G(v)$ is not required to take values in the solenoidal subspace. Next, we recall the It\^o isometry and the Burkholder--Davis--Gundy inequality~\cite[Theorem 1.1]{Marinelli2016}. We let $\Phi$ be a predictable process taking values in $\mathcal{L}(\mathcal{K},L^2(\dom))$ and satisfying
 \begin{equation*}
 	\E\int_0^T
 	\norm{\Phi(s)Q^{1/2}}_{\mathcal L_2(\mathcal K,L^2(\dom))}^2 ds 	<\infty.
 \end{equation*}
 Then $\int_0^t \Phi(s)dW(s)$ is a continuous square-integrable $L^2(\dom)$-valued martingale and satisfies the It\^o isometry~\cite{DaPrato2014}
 \begin{equation*}
 	\E\norm{\int_0^t \Phi(s) dW(s)}^2_{L^2(\dom)} = \E \int_0^t \norm{\Phi(s)Q^{1/2}}^2_{\mathcal{L}_2(\mathcal{K},L^2(\dom))} ds.
 \end{equation*}
 Furthermore, it satisfies the Burkholder--Davis--Gundy (BDG) inequality 
 \begin{equation}\label{eq:BDG}
 	\E\sup_{0\leq t\leq T}\norm{\int_0^t \Phi(s)dW(s)}^q_{L^2(\dom)}\leq C_q \E \left(\int_0^T \norm{\Phi(s)Q^{1/2}}^2_{\mathcal{L}_2(\mathcal{K},L^2(\dom))} ds\right)^{q/2},
 \end{equation}
 for every $q\geq 0$ for which the right-hand side is finite. 
 For a real-valued continuous local martingale $M$ with $M(0)=0$, this translates to the process
 \begin{equation*}
 	M_t^* = \max_{0\leq s \leq t} |M_s|.
 \end{equation*}
 and its quadratic variation $\langle M\rangle_t$ satisfying
 for any $q>0$,
 \begin{equation}\label{eq:BDGscalar}
 	k_q \E \langle M\rangle_\tau^{q/2} \leq \E (M_\tau^*)^q\leq K_q \E \langle M\rangle_\tau^{q/2},
 \end{equation}
where $k_q,K_q>0$ are constants, for every stopping time $\tau\leq T$, e.g.~\cite[Theorem 3.28]{Karatzas1991}.
 Next, we define weak martingale solutions of the Navier--Stokes equations. 
\begin{definition}
	\label{def:martingalesolution}
	
	Let $T>0$ be a finite time horizon, $u_0\in \Ldiv$, $f\in L^2(0,T;H^{-1}(\dom))$ and let $\mathcal{K}$, $Q$ and $G$ satisfy the assumptions stated above. A weak martingale solution of~\eqref{eq:SNSm} - \eqref{eq:initdata} is a stochastic basis $(\Omega,\mathcal{G},(\mathcal{G}_t)_{t\in [0,T]},\mathbb{P})$ satisfying the usual assumptions, a $\mathcal{K}$-valued Wiener process with covariance $Q$ that is adapted to the filtration $(\mathcal{G}_t)_{t\in [0,T]}$, and an $\Ldiv$-valued progressively measurable stochastic process $u$ satisfying
	\begin{equation*}
		u\in C([0,T];H^{-1}(\dom))\cap L^\infty(0,T;L^2_{\Div}(\dom))\cap L^2(0,T;H^1_0(\dom))
	\end{equation*}
	almost surely,
	\begin{equation*}
		\E \left(\sup_{0\leq t\leq T} \norm{u(t)}^2_{L^2(\dom)} + \int_0^T \norm{\Grad u(s)}^2_{L^2(\dom)}ds\right)<\infty,
	\end{equation*}
and $\mathbb{P}$-a.s. for all times $t\in [0,T]$,
	\begin{multline}
		\label{eq:weakformulation}
		(u(t),v)-(u_0,v)+\nu\int_0^t (\Grad u(s),\Grad v) ds -\int_0^t \int_{\dom}(u\cdot\Grad)v \cdot u dx ds \\
		= \int_0^t \langle f(s),v\rangle ds + \left(\int_0^t G(u(s))dW(s),v\right),\quad \forall \, v\in C^\infty_{c,\Div}(\dom),
	\end{multline}
	where $u(0)=u_0$ almost surely.
\end{definition}
\section{Numerical scheme}\label{sec:num}
We start by describing the time discretization.
We let $\Delta t=T/N>0$ be a time step size and let $t^m:=m\Delta t$, $m=0,1,\dots, N\in \N$ be the time steps. Given an initial condition $u_0\in \Ldiv$, we define the following semi-discrete algorithm:
\begin{enumerate}
	\item[{\bf Step 1}]
	For $m=0,1,2,\dots$, given $u^m$ and $\tu^m$, we compute the intermediate velocity field $\tu^{m+1}$ via
	\begin{equation}\label{eq:predictionsemi}
		\tu^{m+1} - u^m + \Delta t\left[(\tu^m\cdot\Grad) \tu^{m+1}+\frac{1}{2}(\Div \tu^m) \tu^{m+1}\right] = \Delta t\nu \Delta \tu^{m+1} +\Delta t f^m+ G(u^m) \Delta W_m,
	\end{equation}
	with boundary conditions 
	\begin{equation*}
		\left.\tu^{m+1}\right|_{\partial\dom}=0,
	\end{equation*}
	and where we have defined
	\begin{equation*}
		f^m=\frac{1}{\Delta t}\int_{t^m}^{t^{m+1}}f(s) ds,\quad \Delta W_m = W(t^{m+1}) - W(t^m).
	\end{equation*}
	\item[{\bf Step 2}] (Projection step) Next, we define $(u^{m+1},p^{m+1})$ via
		\begin{subequations}
			\label{eq:projectionsemi}
			\begin{align}
				\frac{u^{m+1}-\tu^{m+1}}{\Delta t }&= -\Grad p^{m+1},\\
				\Div u^{m+1} & = 0.
			\end{align}	
		\end{subequations}
		with boundary condition
		\begin{equation}
			\label{eq:bcstep2}
			\left.u^{m+1}\cdot n\right|_{\partial \dom}=0.
		\end{equation}

\end{enumerate}
	\subsection{Spatial discretization}\label{sec:spatial}
Next we describe the spatial discretization for the problem. 
We let $\mathcal{T}_h=\{K\}$ be families of conforming shape-regular triangulations of $\dom$ made of simplices with mesh size $h>0$. Here $h=\max_{K\in \mathcal{T}_h}\text{diam}(K)$.   We let $\Uh\subset H^1_0(\dom)^d$, and $\Ph\subset L^2_0(\dom)\cap H^1(\dom)$  be finite dimensional subspaces scaled by the meshsize $h>0$ that we will use for the spatial approximation of the velocity, and the pressure  respectively.  The pair $(\Uh,\Ph)$ has to satisfy certain assumptions that are needed in the following to prove stability and convergence of the scheme. First, we require the following:
\begin{assumption}
	\label{ass:LBB}
	The pair $(\Uh,\Ph)$ needs to satisfy the LBB (Ladyzhenskaya–Babu\v{s}ka–Brezzi) condition~\cite{Boffi2013}, meaning there exists $\beta^*>0$ independent of the mesh size $h>0$ such that  
	\begin{equation}\label{eq:LBB}
		\inf_{0\neq q\in \Ph}\sup_{0\neq v\in \Uh}\frac{(\Div v,q)}{\norm{\Grad v}_{L^2}\norm{q}_{L^2}}\geq \beta^*>0.
	\end{equation}
\end{assumption}
This condition is crucial for obtaining a stable approximation of the velocity-pressure pair of the Stokes equations (see e.g.~\cite{Boffi2013,Boyer2013,Guermond1998} and references therein). In particular, it is required to prove the solvability of the scheme. For examples of pairs of finite element spaces $(\Uh,\Ph)$ that satisfy the LBB condition~\eqref{eq:LBB}, see for instance~\cite{Girault1986,Ern2004}.  
One popular option is the Taylor--Hood finite element pair, but also the MINI element by Arnold, Brezzi and Fortin~\cite{Arnold1984} satisfies the requirement. Let us denote by $k$ the approximation order of the elements in the space $\Uh$ and by $\ell$ the approximation order of the basis functions in $\Ph$. We assume that $k,\ell\geq 1$. We will assume that the following approximation property holds:
\begin{assumption}
	\label{ass:approxorder}
	There exist approximation operators $\Ih$ onto $\Uh$ and $\Ih^p$ onto $\Ph$ such that for a constant $C>0$ independent of $h$,
	\begin{equation*}
		\norm{v-\Ih v}_{W^{m,a}(\dom)}\leq C h^{r-m}\norm{v}_{W^{r,a}(\dom)},
	\end{equation*}
for any $	v\in W^{r,a}(\dom)\cap H^1_0(\dom)$, $a\in [1,\infty]$
	where $r\in \{j,\dots, k+1\}$ and $m\in \{0,1\}$ and $j$ is the smallest positive integer such that $j>\frac{d}{a}$ if $a>1$ and  $j\geq d$ if $a=1$, and
	\begin{equation*}
		\norm{q-\Ih^p q}_{W^{m,a}(\dom)}\leq C h^{r-m}\norm{q}_{W^{r,a}(\dom)},
	\end{equation*}
	for any $	q\in W^{r,a}(\dom)\cap L^2_0(\dom)$,
	where $r\in \{j,\dots, \ell+1\}$ and $m\in \{0,1\}$ and $j$ is the smallest positive integer such that $j>\frac{d}{a}$ if $a>1$ and  $j\geq d$ if $a=1$.
\end{assumption}
For conforming Lagrange finite elements on shape-regular families of affine simplicial meshes, these estimates follow from the local interpolation estimates in~\cite[Theorem 11.13 and Section 11.5.1]{Ern2021} by assembling the interpolants and combining the elementwise
bounds for derivatives of order zero and one. The velocity operator preserves the homogeneous boundary condition. For the pressure operator, the zero-mean constraint can be enforced by subtracting the mean of the interpolant.
This assumption is not needed to prove the solvability and stability of the scheme but it is needed to show convergence. To introduce the next assumption, we define the space of discretely divergence-free functions:
\begin{equation}
	\label{eq:discdivfree}
	\Vh = \left\{v\in \Uh\, : \, (\Div v,q)=0\, \forall \, q\in \Ph \right\}.
\end{equation}
Next, we introduce the Stokes projection $\Sh: H^1_0(\dom)\times L^2_0(\dom)\to \Uh\times\Ph$: $(\varphi_h,\pi_h) = \Sh(\varphi,\pi)$ that satisfies
\begin{equation}\label{eq:stokesprojection}
\nu(\Grad (\varphi_h - \varphi),\Grad v) - (\pi_h-\pi,\Div v ) -(\Div(\varphi_h-\varphi),q) =0,\quad \forall\, (v,q)\in \Uh\times\Ph.
\end{equation}
For $\varphi\in H^1_{\Div}(\dom)$, we write $(\varphi_h,\pi_h) = \Sh(\varphi,0)$. Setting $v=0$ in~\eqref{eq:stokesprojection}, we obtain that $\varphi_h\in \Vh$.
Now let us make the following assumption:
\begin{assumption}
	\label{ass:stokes}
	There is $s\in(0,1]$ such that for every $g\in L^2(\dom)$ the solution $(z,\rho)$, where $\rho\in L^2_0(\dom)$, of the Stokes problem 
	\begin{equation*}
		-\nu\Delta z+\nabla\rho=g,\quad  \Div z=0 \quad \text{in }\, \dom,\quad  z=0 \quad \text{on }\, \partial \dom,
	\end{equation*}
	satisfies 
	\begin{equation*}
		\norm{z}_{H^{1+s}}+\norm{\rho}_{H^s}\le C\norm{g}_{L^2}.
	\end{equation*}
\end{assumption}
\begin{remark}
	\label{rem:stokesregularity}
	Assumption~\ref{ass:stokes} holds with $s=1$ if $\dom$ is a convex polygon~\cite{Kellogg1976}  or a convex polyhedron~\cite{Dauge1989}. On general polytopes, it holds for some $s>0$ depending on the angles of the corners~\cite{Dauge1989,Mazya2010}.
\end{remark}
The Stokes projection satisfies
\begin{equation}
	\label{eq:stokesprojectionestimate}
\norm{\Grad(\varphi-\varphi_h)}_{L^2}+\norm{\pi-\pi_h}_{L^2}\le C\Big(\inf_{v\in\Uh}\norm{\Grad(\varphi-v)}_{L^2}+\inf_{q\in\Ph}\norm{\pi-q}_{L^2}\Big),
\end{equation}
see~\cite[Chapter II, Theorem 1.1]{Girault1986}. 
Using Assumption~\ref{ass:stokes}, we can obtain the following additional estimates:
\begin{lemma}
	\label{lem:stokesestimate}
	If $\Div \varphi=0$ a.e. the Stokes projection $\varphi_h$ of $\varphi\in H^{k+1}(\dom)\cap H^1_0(\dom)$ satisfies
	\begin{equation}\label{eq:stokesimproved1}
	\norm{\Grad(\varphi-\varphi_h)}_{L^2}+ \norm{\pi_h}_{L^2}\le C\inf_{v\in \Uh}\norm{\Grad(\varphi-v)}_{L^2}\le Ch^k\norm{\varphi}_{H^{k+1}}.
	\end{equation}
	Additionally, under Assumption~\ref{ass:stokes}, we have
	\begin{equation}\label{eq:stokesimproved2}
	\norm{\varphi-\varphi_h}_{L^2}\le Ch^{k+s}\norm{\varphi}_{H^{k+1}}.
	\end{equation}
\end{lemma}
\begin{proof}
	For $\pi=0$, the second term on the right-hand side of~\eqref{eq:stokesprojectionestimate} vanishes and we obtain from standard interpolation estimates
	\begin{equation}
		\label{eq:stokesprojectionestimaterefined}
		\norm{\Grad(\varphi-\varphi_h)}_{L^2}+ \norm{\pi_h}_{L^2}\le C h^k \norm{\varphi }_{H^{k+1}(\dom)},
	\end{equation}
	which implies the first estimate. For the second estimate, we again set $\pi=0$ (since we assume that $\varphi$ is divergence free) and we denote $e = \varphi-\varphi_h$. Then~\eqref{eq:stokesprojection} reduces to
	\begin{equation}\label{eq:stokesprojectiondivfree}
		\nu(\Grad e,\Grad v) +(\pi_h,\Div v ) -(\Div e,q) =0,\quad \forall\, (v,q)\in \Uh\times\Ph.
	\end{equation}
	Now consider the dual problem: Find $(z,\rho)\in H^1_0(\dom)\times L^2_0(\dom)$ solving
	\begin{equation*}
		\nu (\Grad z,\Grad v) - (\rho,\Div v) = (e,v),\quad \forall\, v\in H^1_0(\dom),\quad (\Div z, q)=0,\quad \forall\, q\in L^2_0(\dom).
	\end{equation*}
	Using $v=e$ as a test function, we get
	\begin{equation*}
		\norm{e}_{L^2}^2 = \nu(\Grad z,\Grad e)-(\rho,\Div e).
	\end{equation*}
	Subtracting~\eqref{eq:stokesprojectiondivfree} with $(v,q)\in \Uh\times\Ph$ arbitrary, and using that $\Div z=0$, we get
	\begin{equation*}
		\norm{e}^2_{L^2}= \nu(\Grad e,\Grad(z-v)) - (\Div e,\rho-q) + (\pi_h,\Div (z-v)).
	\end{equation*}
	We use that $\norm{\Div w}_{L^2}\leq C \norm{\Grad w}_{L^2}$ to estimate the right-hand side:
	\begin{equation*}
		\norm{e}^2_{L^2} \leq C \left(\norm{\Grad e}_{L^2}+\norm{\pi_h}_{L^2}\right)\left(\norm{\Grad (z-v)}_{L^2}+ \norm{\rho- q}_{L^2}\right).
	\end{equation*}
	Since $(v,q)\in \Uh\times\Ph$ was arbitrary, we can take the infimum over $(v,q)\in \Uh\times\Ph$ on the right-hand side and use best approximation properties~\cite{Ern2017} to obtain
	\begin{equation*}
			\norm{e}^2_{L^2} \leq C \left(\norm{\Grad e}_{L^2}+\norm{\pi_h}_{L^2}\right)\left(h^s\norm{z}_{H^{1+s}}+ h^s\norm{\rho}_{H^s}\right).
	\end{equation*}
	By Assumption~\ref{ass:stokes}, we have
	\begin{equation}\label{eq:assinuse}
		\norm{z}_{H^{1+s}(\dom)}+\norm{\rho}_{H^s(\dom)}\leq C \norm{e}_{L^2(\dom)}. 
	\end{equation}
	Using this, the right-hand side can be bounded as follows:
	\begin{equation*}
		\norm{e}^2_{L^2} \leq C \left(\norm{\Grad e}_{L^2}+\norm{\pi_h}_{L^2}\right)h^s \norm{e}_{L^2}.
	\end{equation*}
	Now, we can divide by $\norm{e}_{L^2}$ and use the first estimate~\eqref{eq:stokesimproved1} to obtain~\eqref{eq:stokesimproved2}:
	\begin{equation*}
		\norm{e}_{L^2} \leq C \left(\norm{\Grad e}_{L^2}+\norm{\pi_h}_{L^2}\right)h^s \leq C h^{s+k}\norm{\varphi}_{H^{k+1}}.
	\end{equation*}
\end{proof}
Next, we define the skew-symmetric trilinear form $b$ for functions $u,v,w\in H^1(\dom)$ with $u\cdot n=0$ on $\partial\dom$ as
\begin{equation}
	\label{eq:defb}
	b(u,v,w) = \int_{\dom} (u\cdot\Grad) v \cdot w dx +\frac12 \int_{\dom} (\Div u) v\cdot w dx =\frac12\int_{\dom} (u\cdot\Grad) v \cdot w -(u\cdot \Grad)w \cdot v dx
\end{equation}
We observe that $b(u,w,v)= - b(u,v,w)$ and hence $b(u,v,v)=0$ for $u,v\in H^1(\dom)$ with $u\cdot n =0$ on $\partial\dom$ or $v=0$ on $\partial\dom$.  

	We start by approximating the initial data for $\tu$:
\begin{equation}
	\label{eq:initdataapprox}
	\tu^0_h = \PUh u_0, 
\end{equation} 
where $\PUh$ is the $L^2$-orthogonal projection, i.e.,
\begin{equation}\label{eq:L2proj}
	(\PUh u,v)=(u,v),\quad \forall v\in \Uh.
\end{equation}
Since $\tu_h^0$ may not be discretely divergence free, we project it to an approximately divergence free function using the following projection step: Find $(u_h^0,p_h^0)\in \Uh\times\Ph$ such that
\begin{subequations}\label{eq:1step}
	\begin{align}
	\left(\frac{u^{0}_h-\tu^{0}_h}{\Delta t }, v\right)&= (p^{0}_h,\Div v),\\
	\left(\Div u^{0}_h,q\right)&=0,	
	\end{align}
\end{subequations}
for all $(v,q)\in \Uh\times \Ph$.
Then we define the following iterative scheme: 
\begin{enumerate}
	\item[{\bf Step 1}] (Prediction step):
	For any $m\in \{0,1,\dots,N-1\}$, given $ u_h^{m},\tu_h^m \in \Uh $, find $\tu_h^{m+1} \in \Uh $, such that for all $v\in \Uh $,
	\begin{equation}
		\label{eq:predictionfully}
 \left (\tu^{m+1}_h-u^m_h ,v\right)+ \Delta tb(\tu^m_h,\tu^{m+1}_h,v)+\Delta t\nu(\Grad \tu_h^{m+1},\Grad v)=\Delta t \langle f^{m}, v\rangle  + (\PUh G(u^m_h)\Delta W_m,v),
	\end{equation}
	\item[{\bf Step 2}] (Projection step): Next we seek $(u^{m+1}_h,p^{m+1}_h)\in \Uh\times\Ph$ which satisfy for all $(v,q)\in \Uh\times\Ph$,
	\begin{subequations}
		\label{eq:projectionfully}
		\begin{align}
			\label{eq:projection1}
			\left(\frac{u^{m+1}_h-\tu^{m+1}_h}{\Delta t }, v\right)&= ( p^{m+1}_h,\Div v),\\
			\label{eq:projection2}
			\left(\Div u^{m+1}_h,q\right)&=0,	
		\end{align}	
	\end{subequations}
\end{enumerate}
Testing \eqref{eq:projection1} with $v\in \Vh$ shows that $u^{m+1}_h$ is the $L^2$-projection of $\tu^{m+1}_h$ onto $\Vh$, so $u^{m+1}_h-\widetilde u^{m+1}_h\perp\Vh$.
\subsection{Well-posedness of the scheme}
We start by showing that the algorithm~\eqref{eq:predictionfully} - \eqref{eq:projectionfully} is well-posed, i.e., that given the approximations at time $m$, the approximations at time $m+1$ can be computed.
\begin{lemma}
	\label{lem:solvability}
	Assume that $u_0\in L^2_{\Div}(\dom)$ and that $f\in L^2(0,T;H^{-1}(\dom))$. Then for any $0<\Delta t,h<1$, there exists a sequence of random variables $(\tu_h^m,u_h^m,p_h^m)_{m=0}^N\subset \Uh\times\Uh\times \Ph$ which for almost all $\omega\in \Omega$ solve~\eqref{eq:predictionfully} - \eqref{eq:projectionfully} and for each $m$, the maps $\omega\mapsto (\tu_h^m(\omega),u_h^m(\omega),p_h^m(\omega))$ are $\mathcal{G}_{t^m}$-measurable. 
\end{lemma}
\begin{proof}
	We first show that given a fixed $\omega\in \Omega$, and $(\tu^m_h(\omega),u^m_h(\omega))$, for $m\geq 0$, there exists a unique $(\tu^{m+1}_h,u^{m+1}_h,p^{m+1}_h)$ solving \eqref{eq:predictionfully} - \eqref{eq:projectionfully}. The zeroth step is identical to the projection step~\eqref{eq:projectionfully}, so the solvability of this step will follow once we have proved the general case. So we fix $\omega\in \Omega$ and first show that given $u_h^m(\omega)$ bounded, we can find a unique $\tu^{m+1}_h(\omega)$ solving~\eqref{eq:predictionfully}. Firstly, for a given $\omega$, we have $\norm{\Delta W_m(\omega)}_{\mathcal{K}}\leq C$ (with the constant depending on $\omega$). Now we write the scheme~\eqref{eq:predictionfully} as
	\begin{equation*}
		(A^m\tu_h^{m+1},v)= a^m(\tu_h^{m+1},v) = F^m(v),
	\end{equation*}
	where $a^m:\Uh\times\Uh \to \R$ and $F^m:\Uh\to \R$ are defined by
	\begin{equation*}
		a^m(\tu,v)= (\tu,v)+\Delta t b(\tu_h^m(\omega),\tu,v) +\Delta t \nu (\Grad \tu,\Grad v),
	\end{equation*}
	and
	\begin{equation*}
		F^m(v) = (u^m_h(\omega),v)+\Delta t \langle f^m,v\rangle+(G(u^m_h(\omega))\Delta W_m(\omega),v).
	\end{equation*}
	We show that $a^m$ is coercive and bounded and that $F^m$ is bounded on $\Uh$. Then it follows from the Lax--Milgram theorem that $\tu^{m+1}_h(\omega)$ exists and is unique. We will omit the dependence on $\omega$ for the following calculations.
	To show coercivity, we plug in $v=\tu$:
	\begin{equation*}
		a^m(\tu,\tu) =  	\norm{\tu}_{L^2}^2+ \Delta t b(\tu^m_h,\tu ,\tu)+\Delta t\nu(\Grad \tu ,\Grad \tu)
		= \norm{\tu}_{L^2}^2+ \Delta t\nu\norm{\Grad \tu}_{L^2}^2
		\geq c\norm{\tu}_{H^1_0}^2,
	\end{equation*}
	using the skew-symmetry of $b$. 
	This proves the coercivity on $\Uh$. For the boundedness, we have
	\begin{align*}
		\left| 	a^m(\tu,v)\right| \leq  &  	\norm{\tu}_{L^2}\norm{v}_{L^2}+ \Delta t \left( \norm{\tu^m_h}_{L^4}\norm{\tu}_{H^1_0}\norm{v}_{L^4}+\norm{\tu^m_h}_{L^4}\norm{v}_{H^1_0}\norm{\tu}_{L^4}\right)+\Delta t\nu\norm{\Grad \tu}_{L^2}\norm{\Grad v}_{L^2}\\
		\leq & \norm{\tu}_{L^2}\norm{v}_{L^2}+ \Delta t( \norm{\tu^m_h}_{L^4}+\nu)\norm{\tu}_{H^1_0}\norm{v}_{H^1_0}\leq C(\omega) \norm{\tu}_{H^1_0}\norm{v}_{H^1_0}, 
	\end{align*}
	where we used the Sobolev embedding theorem and the Poincar\'e inequality. Here we also used that $\tu^m_h\in\Uh\subset H^1_0(\dom)$.
	Hence $a^m$ is bounded on $\Uh$. 
	To show boundedness of $F^m$ (for every $m$), we estimate:
	\begin{align*}
		|F^m(v)|&\leq \norm{u^m_h}_{L^2}\norm{v}_{L^2}+\Delta t \norm{  f^{m}}_{H^{-1}}\norm{v}_{H^1_0}+ \norm{G(u^m_h)}_{\mathcal{L}(\mathcal{K},L^2)}\norm{\Delta W_m}_{\mathcal{K}}\norm{v}_{L^2}\\
		&\leq \norm{u^m_h}_{L^2}\norm{v}_{L^2}+\Delta t \norm{  f^{m}}_{H^{-1}}\norm{v}_{H^1_0}+ C K_1\left(\norm{u^m_h}_{L^2}+1\right)\norm{\Delta W_m}_{\mathcal{K}}\norm{v}_{L^2}.
	\end{align*}
	Thus, with the Poincar\'e inequality, we obtain boundedness of $F^m$ on $\Uh$.
	This implies the solvability of Step 1 of the algorithm. 
		Step 2 follows from the choice of the LBB-stable pair $(\Uh,\Ph)$, cf.~\cite{Guermond1996,Guermond1998} and also~\cite{chorinprojection}, since this step is done exactly as in the deterministic case. We provide it here for completeness.
		We first define the operator $B_h : \Uh\to \Ph$ by 
		\begin{equation*}
			(B_h v ,q ):=-(\Div v ,q ) = ( v,B_h^\top q ) ,\quad \forall\, q \in \Ph,\, v\in \Uh,
		\end{equation*}
		We note that $\Vh$ as defined in~\eqref{eq:discdivfree} is exactly the kernel of $B_h$, i.e., $\Vh = \ker B_h$ and also $(\operatorname{range}B_h^\top)^\perp = \ker B_h$. Since $\Uh$ is finite dimensional, we obtain
		\begin{equation*}
			\Uh = \Vh\oplus \text{range}B_h^\top.
		\end{equation*}
		Thanks to the LBB-condition~\eqref{eq:LBB}, $B_h^\top$ is injective. Therefore, we can decompose $\tu_h^{m+1}\in \Uh$ as
		\begin{equation*}
			\tu_h^{m+1} = \PVh \tu_h^{m+1} + B_h^\top q_h^{m+1},
		\end{equation*}
		for a unique $ \PVh \tu_h^{m+1}$ (as it is the orthogonal projection) and a unique $q_h^{m+1}\in \Ph$ by the injectivity of $B_h^\top$. Comparing to~\eqref{eq:projectionfully}, we see that $u_h^{m+1}=  \PVh \tu_h^{m+1}$ and $q_h^{m+1} = {\Delta t}p_h^{m+1}$.
			Next, we show the measurability. Clearly, since $u_0$ is deterministic,  it is measurable with respect to $\mathcal{G}_{t^0}$. Now, assume that $\tu^m_h$ and $u^m_h$ are $\mathcal{G}_{t^m}$-measurable. We fix a basis $\{\phi_i\}_{i=1}^{n_h}$ of $\Uh$ and write~\eqref{eq:predictionfully} as $A^m(\omega)\xi =F^m(\omega)$, where
	\begin{align*}
		A_{ij}^m(\omega) &= (\phi_j,\phi_i) + \Delta t b(\tu_h^m(\omega),\phi_j,\phi_i)+\Delta t \nu (\Grad \phi_j,\Grad\phi_i),\\
		F^m_i(\omega) &= (u_h^m(\omega),\phi_i) + \Delta t \langle f^m,\phi_i\rangle + (G(u^m_h(\omega))\Delta W_m(\omega),\phi_i).
	\end{align*}
	We note that $A^m$ is an affine function of $\tu_h^m(\omega)$ and therefore $\mathcal{G}_{t^m}$-measurable, and it is continuously invertible on $GL(n_h)$ for each $\omega$ due to the coercivity proved above, and so $(A^m(\omega))^{-1}$ is $\mathcal{G}_{t^m}$-measurable. $F^m$ is $\mathcal{G}_{t^{m+1}}$-measurable because $u^m_h$ is $\mathcal{G}_{t^m}$-measurable and because $(u,r)\mapsto G(u)r$ is continuous on $L^2(\dom)\times\mathcal{K}\to L^2(\dom)$ by~\eqref{eq:Glipschitz}. Therefore $(A^m)^{-1}F^m$ is $\mathcal{G}_{t^{m+1}}$-measurable and so $\tu_h^{m+1}$ is $\mathcal{G}_{t^{m+1}}$-measurable. 
	 The second, projection step of the algorithm also corresponds to solving a linear system $A_2((u^{m+1}_h,p^{m+1}_h)) = \tu_h^{m+1}$, therefore $u^{m+1}_h$ and $p^{m+1}_h$ are $\mathcal{G}_{t^{m+1}}$-measurable too.
\end{proof}
\subsection{Energy estimate}
Next, we prove a discrete energy estimate which yields a priori estimates for the numerical method that will be important in passing to the limit $h,\Delta t\to 0$ later.
\begin{lemma}
	\label{lem:energyestimate}
	
	Assume that $u_0\in \Ldiv$ and that $f\in L^{2^{r_0}}(0,T;H^{-1}(\dom))$ for some $r_0\in \N$, $r_0\geq 3$.
The approximations computed by the scheme~\eqref{eq:predictionfully} - \eqref{eq:projectionfully} satisfy for constants $C$ and $C(T,u_0,f,\tr Q)$ not depending on $h$ or $\Delta t$ (but depending on $T$, $u_0$, $f$, $\tr Q$, $\nu$, $K_1$ and $r$),
	\begin{enumerate}[(i)]
		\item 
		\begin{multline*}
				\max_{1\leq m\leq N} \E\norm{u^m_h}^2_{L^2} + \Delta t\nu \sum_{m=1}^N\E\norm{\Grad \tu_h^m}^2_{L^2} +\sum_{m=0}^{N-1}\E\norm{\tu^{m+1}_h - u^m_h}_{L^2}^2+\sum_{m=0}^{N-1}\E\norm{\tu^{m+1}_h - u^{m+1}_h}_{L^2}^2\\
			\leq  C e^{CT} \left(\norm{u_0}_{L^2}^2 + \norm{f}^2_{L^2(0,T;H^{-1}(\dom))}+1\right),
		\end{multline*}
		\item \begin{equation*}
			\E\max_{1\leq i \leq N}\norm{u^i_h}_{L^2}^2 + \E \max_{1\leq i\leq N}\norm{\tu^{i}_h - u^{i-1}_h}_{L^2}^2+\E \max_{1\leq i\leq N}\norm{\tu^{i}_h - u^{i}_h}_{L^2}^2\leq C(T, u_0,f,\tr Q),
		\end{equation*}		
		\item We have for any $r\in \N$, $1\leq r\leq r_0$,
		\begin{multline*}
			\max_{1\leq i \leq N}\E\norm{u^i_h}_{L^2}^{2^r} + \E\sum_{m=0}^{N-1}\left|\norm{u^{m+1}_h}_{L^2}^{2^{r-1}}-\norm{u^{m}_h}_{L^2}^{2^{r-1}}\right|^2+\E \sum_{m=0}^{N-1}\norm{u_h^{m+1}}_{L^2}^{2^r-2}\norm{\tu_h^{m+1}-u^m_h}_{L^2}^2\\
			+
			\E \sum_{m=0}^{N-1}\norm{u_h^{m+1}}_{L^2}^{2^r-2}\norm{\tu_h^{m+1}-u^{m+1}_h}_{L^2}^2 +
			\nu\Delta t\E \sum_{m=0}^{N-1}\norm{u_h^{m+1}}_{L^2}^{2^r-2}\norm{\Grad \tu_h^{m+1}}_{L^2}^2
			\leq C(T,f,u_0,\tr Q)
		\end{multline*}
		and
			 \begin{equation*}
			\max_{1\leq i \leq N}	\E\norm{\tu_h^{i}}_{L^2}^{2^r}\leq C(T,u_0,f,\tr Q).
		\end{equation*}
		\item We have for any $r\in \N$, $1\leq r\leq r_0$,
		\begin{equation*}
			\E\max_{1\leq m\leq N}\norm{u^m_h}_{L^2}^{2^r} \leq C(T, u_0,f,\tr Q),\quad \E\max_{1\leq m\leq N}\norm{\tu^m_h}_{L^2}^{2^r} \leq C(T, u_0,f,\tr Q).
		\end{equation*}
		\item For any $r\in \N$, $1\leq r\leq r_0$, we have 
		\begin{align*}
		\E\left(\nu \Delta t \sum_{m=0}^{N-1}\norm{\Grad \tu_h^{m+1}}_{L^2}^2\right)^{2^{r-1}} &\leq C(T,u_0,f,\tr Q),\\
		\E\left(\sum_{m=0}^{N-1}\norm{\tu^{m+1}_h - u^m_h}_{L^2}^2\right)^{2^{r-1}}&\leq C(T,u_0,f,\tr Q),\\
		\E \left(\sum_{m=0}^{N-1}\norm{\tu^{m+1}_h - u^{m+1}_h}_{L^2}^2\right)^{2^{r-1}}&\leq C(T,u_0,f,\tr Q).
		\end{align*}
	\end{enumerate}
\end{lemma}
\begin{proof}
	To prove (i), we take $ v= \tu^{m+1}_h$ as a test function in~\eqref{eq:predictionfully} and use the skew-symmetry of $b$ and the identity $a(a-b) = \frac12 (a^2- b^2+(a-b)^2)$:
	\begin{multline}\label{eq:energy1step}
		\frac{1}{2}\norm{\tu^{m+1}_h}_{L^2}^2 -\frac{1}{2} \norm{u^m_h}_{L^2}^2+\frac12\norm{\tu^{m+1}_h - u^m_h}_{L^2}^2 +\Delta t \nu \norm{\Grad \tu^{m+1}_h}_{L^2}^2 \\
		= \Delta t \langle f^m,\tu^{m+1}_h\rangle + (\PUh G(u^m_h)\Delta W_{m},\tu_h^{m+1}).
	\end{multline}
	Then we take $u^{m+1}_h$ as a test function in~\eqref{eq:projection1} and use~\eqref{eq:projection2} and again $a(a-b) = \frac12 (a^2- b^2+(a-b)^2)$ to obtain
	\begin{equation}\label{eq:projectionestimate}
		\frac12\norm{u^{m+1}_h}^2_{L^2}-\frac12 \norm{\tu^{m+1}_h}_{L^2}^2 +\frac12\norm{u^{m+1}_h - \tu^{m+1}_h}_{L^2}^2 = 0.
	\end{equation}
	Using this in~\eqref{eq:energy1step}, we obtain
	\begin{multline}\label{eq:energy1astep}
			\frac{1}{2}\norm{u^{m+1}_h}_{L^2}^2 -\frac{1}{2} \norm{u^m_h}_{L^2}^2+\frac12\norm{\tu^{m+1}_h - u^m_h}_{L^2}^2+\frac12\norm{\tu^{m+1}_h - u^{m+1}_h}_{L^2}^2 +\Delta t \nu \norm{\Grad \tu^{m+1}_h}_{L^2}^2 \\
			= \Delta t \langle f^m,\tu^{m+1}_h\rangle + (\PUh G(u^m_h)\Delta W_{m},\tu_h^{m+1}-u_h^m)+(\PUh G(u^m_h)\Delta W_{m},u_h^m).
	\end{multline}
	From~\eqref{eq:energy1step}, by induction and independence of $\Delta W_m$ from $\mathcal{G}_{t^m}$, and the boundedness of moments of Gaussian random fields~\cite[Proposition 2.19]{DaPrato2014}, we obtain that all moments of $\tu_h^{m+1}$ are bounded for each fixed $\Delta t,h$, and therefore the moments of $u^{m+1}_h$ are finite too. So we can take expectations in~\eqref{eq:energy1astep}, use the independence of the increments of the Wiener process, and sum from $m=0$ to $m=N-1$:
	\begin{multline}\label{eq:energy2step}
		\frac{1}{2}\E\norm{u^{N}_h}_{L^2}^2 -\frac{1}{2} \norm{u^0_h}_{L^2}^2+\frac12\sum_{m=0}^{N-1}\E\norm{\tu^{m+1}_h - u^m_h}_{L^2}^2+\frac12\sum_{m=0}^{N-1}\E\norm{\tu^{m+1}_h - u^{m+1}_h}_{L^2}^2 +\Delta t \nu \sum_{m=0}^{N-1}\E\norm{\Grad \tu^{m+1}_h}_{L^2}^2 \\
		= \Delta t \sum_{m=0}^{N-1}\E \langle f^m,\tu^{m+1}_h\rangle + \sum_{m=0}^{N-1}\E( G(u^m_h)\Delta W_{m},\tu_h^{m+1}-u_h^m).
	\end{multline}
	Now we bound the terms on the right-hand side as follows:
	\begin{equation*}
 \left|\Delta t \sum_{m=0}^{N-1}\E\langle f^m,\tu^{m+1}_h\rangle \right| \leq \frac12\Delta t \left( \frac{1}{\nu}\sum_{m=0}^{N-1}\norm{ f^m}_{H^{-1}(\dom)}^2 + \nu\sum_{m=0}^{N-1}\E\norm{\Grad \tu^{m+1}_h}_{L^2}^2    \right)
	\end{equation*}
	and using the law of iterated expectations, we can bound the other term,
	\begin{equation}\label{eq:Gestimate}
		\begin{split}
	&	\left|   \sum_{m=0}^{N-1}\E( G(u^m_h)\Delta W_{m},\tu_h^{m+1}-u_h^m)\right| \\
	&\qquad  \leq  \sum_{m=0}^{N-1}\E\norm{ G(u^m_h)\Delta W_{m}}_{L^2}^2 + \frac14\sum_{m=0}^{N-1}\E\norm{\tu_h^{m+1}-u_h^m}_{L^2}^2\\
		& \qquad\leq  \sum_{m=0}^{N-1}\E\left(\norm{ G(u^m_h)}_{\mathcal{L}(\mathcal{K},L^2)}^2\norm{\Delta W_{m}}_{\mathcal{K}}^2 \right)+ \frac14\sum_{m=0}^{N-1}\E\norm{\tu_h^{m+1}-u_h^m}_{L^2}^2\\
		&\qquad =  \sum_{m=0}^{N-1}\E \left[\E\left[\norm{ G(u^m_h)}_{\mathcal{L}(\mathcal{K},L^2)}^2\norm{\Delta W_{m}}_{\mathcal{K}}^2 |\mathcal{G}_{t^m}\right]\right]+ \frac14\sum_{m=0}^{N-1}\E\norm{\tu_h^{m+1}-u_h^m}_{L^2}^2\\
		&\qquad = \sum_{m=0}^{N-1}\E \left[\norm{ G(u^m_h)}_{\mathcal{L}(\mathcal{K},L^2)}^2\E\left[\norm{\Delta W_{m}}_{\mathcal{K}}^2 |\mathcal{G}_{t^m}\right]\right]+ \frac14\sum_{m=0}^{N-1}\E\norm{\tu_h^{m+1}-u_h^m}_{L^2}^2\\
			&\qquad = \sum_{m=0}^{N-1}\Delta t (\tr Q)\E \left[\norm{ G(u^m_h)}_{\mathcal{L}(\mathcal{K},L^2)}^2\right]+ \frac14\sum_{m=0}^{N-1}\E\norm{\tu_h^{m+1}-u_h^m}_{L^2}^2\\
		& \qquad \leq 2\Delta t(\tr Q)K_1^2\sum_{m=0}^{N-1}\left( \E\norm{u^m_h}_{L^2}^2 +1\right)+ \frac14\sum_{m=0}^{N-1}\E\norm{\tu_h^{m+1}-u_h^m}_{L^2}^2.
		\end{split}
	\end{equation}
	The first step~\eqref{eq:1step} and~\eqref{eq:initdataapprox} yield
	\begin{equation*}
		\norm{u_h^0}_{L^2}^2 + \norm{u_h^0 -\tu_h^0}_{L^2}^2 = \norm{\tu_h^0}_{L^2}^2 = \norm{\PUh u_0}_{L^2}^2\leq \norm{u_0}_{L^2}^2.
	\end{equation*}
	Using these estimates in~\eqref{eq:energy2step}, we obtain
		\begin{multline}\label{eq:energy3step}
		\frac{1}{2}\E\norm{u^{N}_h}_{L^2}^2+\frac14\sum_{m=0}^{N-1}\E\norm{\tu^{m+1}_h - u^m_h}_{L^2}^2+\frac12 \sum_{m=0}^{N-1}\E\norm{\tu^{m+1}_h - u^{m+1}_h}_{L^2}^2\\
		 +\frac{\Delta t\nu}{2}  \sum_{m=0}^{N-1}\E\norm{\Grad \tu^{m+1}_h}_{L^2}^2 \\
		\leq  \frac{1}{2} \norm{u_0}_{L^2}^2 + \frac{\Delta t }{2\nu}\sum_{m=0}^{N-1}\norm{f^m}_{H^{-1}(\dom)}^2 +  2\Delta t(\tr Q)K_1^2\sum_{m=0}^{N-1}\left( \E\norm{u^m_h}_{L^2}^2+1 \right) .
	\end{multline}
Thus, a discrete version of Gr\"onwall's inequality,~\cite[Proposition 4.1]{Emmrich1999}, yields
\begin{equation}\label{eq:uestimate1}
	\max_{1\leq m\leq N} \E\norm{u^m_h}^2_{L^2} \leq C e^{CT} \left(\norm{u_0}_{L^2}^2 + \norm{f}^2_{L^2(0,T;H^{-1}(\dom))}+1 \right).
\end{equation}
Furthermore, this estimate, combined with~\eqref{eq:energy3step} implies
\begin{multline*}
	\Delta t\nu \sum_{m=1}^N\E\norm{\Grad \tu_h^m}_{L^2}^2 +\sum_{m=0}^{N-1}\E\norm{\tu^{m+1}_h - u^m_h}_{L^2}^2+\sum_{m=0}^{N-1}\E\norm{\tu^{m+1}_h - u^{m+1}_h}_{L^2}^2\\
	 \leq  C e^{CT} \left(\norm{u_0}_{L^2}^2 + \norm{f}^2_{L^2(0,T;H^{-1}(\dom))} +1\right),
\end{multline*}
which implies with~\eqref{eq:projectionestimate} 
\begin{equation}\label{eq:utildeestimate}
		\max_{1\leq m\leq N} \E\norm{\tu^m_h}^2_{L^2} =
		\max_{1\leq m\leq N} \E\left(\norm{u^m_h}^2_{L^2} +\norm{u^m_h-\tu_h^m}_{L^2}^2 \right)  \leq C e^{CT} \left(\norm{u_0}_{L^2}^2 + \norm{f}^2_{L^2(0,T;H^{-1}(\dom))}+1 \right).
\end{equation}
This proves the first part of the lemma. 
To prove estimate (ii), we start from~\eqref{eq:energy1astep} and apply Young's inequality to the terms on the right-hand side
	\begin{multline*}
	\frac{1}{2}\norm{u^{m+1}_h}_{L^2}^2 -\frac{1}{2} \norm{u^m_h}_{L^2}^2+\frac14\norm{\tu^{m+1}_h - u^m_h}_{L^2}^2+\frac12\norm{\tu^{m+1}_h - u^{m+1}_h}_{L^2}^2   \\
	\leq  \frac{\Delta t}{2\nu}\norm{f^m}_{H^{-1}}^2 +\norm{ G(u^m_h)\Delta W_{m}}_{L^2}^2+( G(u^m_h)\Delta W_{m},u_h^m).
\end{multline*}
Summing from $m=0$ to $m=i-1$,  we have
	\begin{multline}\label{eq:dontknow}
	\frac{1}{2}\norm{u^{i}_h}_{L^2}^2 +\frac14\sum_{m=0}^{i-1}\norm{\tu^{m+1}_h - u^m_h}_{L^2}^2+\frac12\sum_{m=0}^{i-1}\norm{\tu^{m+1}_h - u^{m+1}_h}_{L^2}^2   \\
	\leq  \frac{1}{2}\norm{u^0_h}_{L^2}^2+ \frac{\Delta t}{2\nu}\sum_{m=0}^{i-1}\norm{f^m}_{H^{-1}}^2 +\sum_{m=0}^{i-1}\norm{ G(u^m_h)\Delta W_{m}}_{L^2}^2+\sum_{m=0}^{i-1}( G(u^m_h)\Delta W_{m},u_h^m).
\end{multline}
Taking the maximum over $i=1,\dots, N$ and using that the maximum of sums of positive terms is reached for $i=N$, and then taking expectations, we obtain
\begin{multline*}
\E\max_{1\leq i \leq N}\norm{u^i_h}_{L^2}^2
	\leq   \norm{u_h^0}_{L^2}^2 + \frac{\Delta t}{\nu} \sum_{m=0}^{N-1}\norm{ f^m}_{H^{-1}(\dom)}^2  \\
	+2\sum_{m=0}^{N-1} \E \norm{ G(u^m_h)\Delta W_{m}}_{L^2}^2+2\E\max_{1\leq i\leq N}\sum_{m=0}^{i-1}( G(u^m_h)\Delta W_{m},u_h^m).
\end{multline*} 
The first two terms on the right-hand side are bounded by the assumptions on $u_0$ and $f$. To bound the third term on the right-hand side, we use the estimates in~\eqref{eq:Gestimate} and (i) of the lemma that we have already shown. Finally, for the fourth term, we use the Burkholder--Davis--Gundy inequality~\eqref{eq:BDGscalar}. We define
\begin{equation*}
	M(t) = \int_0^{t}(G(u_h(s))dW(s),u_h(s)),
\end{equation*}
where 
\begin{equation*}
	u_h(t) = u_h^m,\quad t\in [t^m,t^{m+1}),\quad u_h(T)=u_h^N.
\end{equation*}
Then
\begin{equation*}
	M(t^i) = \int_0^{t^i}(G(u_h(s))dW(s),u_h(s)) = \sum_{m=0}^{i-1}( G(u_h^m)\Delta W_{m},u_h^m),
\end{equation*}
$M(t)$ is a continuous martingale since $u_h$ is piecewise constant and adapted and the quadratic variation satisfies
\begin{equation*}
	\langle M\rangle_{t^N}=\int_0^{t^N}\norm{Q^{1/2}G(u_h)^*u_h}^2_{\mathcal{K}} ds\leq \tr Q\Delta t \sum_{m=0}^{N-1}\norm{G(u_h^m)}_{\mathcal{L}(\mathcal{K},L^2)}^2\norm{u_h^m}_{L^2}^2,
\end{equation*}
and so we have
\begin{align*}
2\E\max_{1\leq i\leq N}\sum_{m=0}^{i-1}( G(u_h^m)\Delta W_{m},u_h^m)& \leq C \E \left( \Delta t  \sum_{m=0}^{N-1}\norm{G(u^m_h)}_{\mathcal{L}(\mathcal{K},L^2)}^2 \norm{u_h^m}_{L^2}^2\right)^{1/2}\\
& \leq C  \E \max_{0\leq i\leq  N-1}\norm{u_h^i}_{L^2} \left(\Delta t \sum_{m=0}^{N-1}\norm{G(u^m_h)}_{\mathcal{L}(\mathcal{K},L^2)}^2 \right)^{1/2}\\
&  \leq \frac14\E  \max_{1\leq i\leq  N}\norm{u_h^i}^2_{L^2} + \frac14 \norm{u_h^0}_{L^2}^2  + C \Delta t\E  \sum_{m=0}^{N-1}\norm{G(u^m_h)}_{\mathcal{L}(\mathcal{K},L^2)}^2 \\
&  \leq \frac14\E  \max_{1\leq i\leq  N}\norm{u_h^i}^2_{L^2} + \frac14 \norm{u_h^0}_{L^2}^2  + CK_1^2 \Delta t\E  \sum_{m=0}^{N-1}\left(\norm{u^m_h}_{L^2}^2+1\right).
\end{align*}
Now we can absorb the first term on the right-hand side by the left-hand side and estimate the third term on the right-hand side by the previously proved estimate~\eqref{eq:uestimate1}. This implies (ii). From (i), we also get
\begin{equation*}
	\E \max_{1\leq i \leq N} \norm{\tu^i_h-u^{i}_h}^2_{L^2}\leq \sum_{m=1}^N \E \norm{\tu_h^m-u^m_h}_{L^2}^2 \leq C(T,u_0,f,\tr Q),
\end{equation*}
which implies
\begin{equation}
	\label{eq:secondutildeestimate}
	\E\max_{1\leq i \leq N}\norm{\tu^i_h}_{L^2}^2 \leq C(T,u_0,f,\tr Q).
\end{equation}
We first prove assertion (iii) for $r=2$ and for $u_h^m$. To do so, we multiply~\eqref{eq:energy1astep} by $\norm{u^{m+1}_h}_{L^2}^2$ and use $a(a-b) = \frac12(a^2-b^2 +(a-b)^2)$ once more,
\begin{multline}\label{eq:provingiii}
	\frac14\left(\norm{u^{m+1}_h}_{L^2}^4-\norm{u^m_h}_{L^2}^4 + \left| \norm{u^{m+1}_h}_{L^2}^2-\norm{u^m_h}_{L^2}^2\right|^2\right) +\frac12\norm{u^{m+1}_h-\tu^{m+1}_h}_{L^2}^2\norm{u^{m+1}_h}_{L^2}^2\\ +\frac12\norm{\tu_h^{m+1} - u^m_h}_{L^2}^2\norm{u^{m+1}_h}_{L^2}^2 + \nu\Delta t \norm{\Grad \tu_h^{m+1}}_{L^2}^2\norm{u^{m+1}_h}_{L^2}^2\\
	=\Delta t \underbrace{ \langle f^m,\tu^{m+1}_h\rangle \norm{u^{m+1}_h}_{L^2}^2 }_{\text{I}}+ \underbrace{(G(u_h^m)\Delta W_m,\tu^{m+1}_h-u^m_h)\norm{u^{m+1}_h}_{L^2}^2}_{\text{II}}  + \underbrace{(G(u_h^m)\Delta W_m,u^{m}_h)\norm{u^{m+1}_h}_{L^2}^2 }_{\text{III}}.
\end{multline}
We estimate terms I, II and III. For I, we have using Young's and Poincar\'e's inequality and~\eqref{eq:projectionestimate}:
\begin{align*}
	\text{I} & \leq \norm{f^m}_{H^{-1}}\norm{\Grad \tu_h^{m+1}}_{L^2}\norm{u^{m+1}_h}_{L^2}^2\\
	& \leq \frac{1}{\nu}\norm{f^m}_{H^{-1}}^2 \norm{u^{m+1}_h}_{L^2}^2 + \frac{\nu}{4}\norm{\Grad \tu_h^{m+1}}_{L^2}^2\norm{u^{m+1}_h}_{L^2}^2\\
	& \leq \frac{C}{\nu^3}\norm{f^m}_{H^{-1}}^4 + \frac{\nu}{4 C(\dom)^2}\norm{u^{m+1}_h}_{L^2}^4 + \frac{\nu}{4}\norm{\Grad \tu_h^{m+1}}_{L^2}^2\norm{u^{m+1}_h}_{L^2}^2\\
	& \leq \frac{C}{\nu^3}\norm{f^m}_{H^{-1}}^4 + \frac{\nu}{4 C(\dom)^2}\norm{u^{m+1}_h}_{L^2}^2\norm{\tu_h^{m+1}}_{L^2}^2 + \frac{\nu}{4}\norm{\Grad \tu_h^{m+1}}_{L^2}^2\norm{u^{m+1}_h}_{L^2}^2\\
	& \leq \frac{C}{\nu^3}\norm{f^m}_{H^{-1}}^4 + \frac{\nu}{4}\norm{\Grad \tu^{m+1}_h}_{L^2}^2\norm{u_h^{m+1}}_{L^2}^2 + \frac{\nu}{4}\norm{\Grad \tu_h^{m+1}}_{L^2}^2\norm{u^{m+1}_h}_{L^2}^2.
\end{align*}
The second term we estimate as follows:
\begin{align*}
	\text{II} & \leq \norm{G(u_h^m)}^2_{\mathcal{L}(\mathcal{K},L^2)}\norm{\Delta W_m}^2_{\mathcal{K}}\norm{u^{m+1}_h}_{L^2}^2 + \frac14\norm{\tu^{m+1}_h-u_h^m}_{L^2}^2\norm{u_h^{m+1}}^2_{L^2}\\
	& =\norm{G(u_h^m)}^2_{\mathcal{L}(\mathcal{K},L^2)}\norm{\Delta W_m}^2_{\mathcal{K}}\left( \norm{u^m_h}_{L^2}^2+\left[\norm{u^{m+1}_h}_{L^2}^2-\norm{u^m_h}_{L^2}^2\right]\right) + \frac14\norm{\tu^{m+1}_h-u_h^m}_{L^2}^2\norm{u_h^{m+1}}^2_{L^2}\\
	& \leq 	\norm{G(u_h^m)}^2_{\mathcal{L}(\mathcal{K},L^2)}\norm{\Delta W_m}^2_{\mathcal{K}}\norm{u^m_h}_{L^2}^2 + 4\norm{G(u_h^m)}^4_{\mathcal{L}(\mathcal{K},L^2)}\norm{\Delta W_m}^4_{\mathcal{K}}+\frac{1}{16}\left| \norm{u^{m+1}_h}_{L^2}^2-\norm{u^m_h}_{L^2}^2\right|^2\\
	&\qquad  + \frac14\norm{\tu^{m+1}_h-u_h^m}_{L^2}^2\norm{u_h^{m+1}}^2_{L^2}.
\end{align*} 
For the third term, we have 
\begin{align*}
	\text{III}&=(G(u_h^m)\Delta W_m,u^{m}_h)\left( \norm{u^m_h}_{L^2}^2+\left[\norm{u^{m+1}_h}_{L^2}^2  -\norm{u^m_h}_{L^2}^2\right]  \right)\\
	&  \leq (G(u_h^m)\Delta W_m,u^{m}_h) \norm{u^m_h}_{L^2}^2  + 4 \norm{G(u_h^m)}_{\mathcal{L}(\mathcal{K},L^2)}^2 \norm{\Delta W_m}_{\mathcal{K}}^2 \norm{u^m_h}_{L^2}^2 + \frac{1}{16}\left|\norm{u^{m+1}_h}_{L^2}^2 - \norm{u^m_h}^2_{L^2}\right|^2
\end{align*}
Using these three estimates in~\eqref{eq:provingiii}, we obtain
\begin{multline}\label{eq:provingiii2}
	\frac14\left(\norm{u^{m+1}_h}_{L^2}^4-\norm{u^m_h}_{L^2}^4\right) + \frac18\left| \norm{u^{m+1}_h}_{L^2}^2-\norm{u^m_h}_{L^2}^2\right|^2 +\frac12\norm{u^{m+1}_h-\tu^{m+1}_h}_{L^2}^2\norm{u^{m+1}_h}_{L^2}^2\\ +\frac14\norm{\tu_h^{m+1} - u^m_h}_{L^2}^2\norm{u^{m+1}_h}_{L^2}^2 + \frac{\nu\Delta t}{2} \norm{\Grad \tu_h^{m+1}}_{L^2}^2\norm{u^{m+1}_h}_{L^2}^2\\
	\leq \frac{C}{\nu^3}\Delta t \norm{f^m}_{H^{-1}}^4 +5\norm{G(u_h^m)}^2_{\mathcal{L}(\mathcal{K},L^2)}\norm{\Delta W_m}^2_{\mathcal{K}}\norm{u^m_h}_{L^2}^2 + 4\norm{G(u_h^m)}^4_{\mathcal{L}(\mathcal{K},L^2)}\norm{\Delta W_m}^4_{\mathcal{K}}\\
	+ (G(u_h^m)\Delta W_m,u^{m}_h) \norm{u^m_h}_{L^2}^2. 
\end{multline}
After taking expectations, the last term on the right-hand side vanishes. The other two stochastic terms on the right-hand side, we estimate using the law of iterated conditional expectations as in~\eqref{eq:Gestimate}:
\begin{align}\label{eq:no1}
	\begin{split}
	\E \norm{G(u_h^m)}^2_{\mathcal{L}(\mathcal{K},L^2)}\norm{\Delta W_m}^2_{\mathcal{K}}\norm{u^m_h}_{L^2}^2 & = \E\left[\E\left[ \norm{G(u_h^m)}^2_{\mathcal{L}(\mathcal{K},L^2)}\norm{\Delta W_m}^2_{\mathcal{K}}\norm{u^m_h}_{L^2}^2|\mathcal{G}_{t^m}\right]   \right]\\
	& = \E\left[\norm{G(u_h^m)}^2_{\mathcal{L}(\mathcal{K},L^2)}\norm{u^m_h}_{L^2}^2\E\left[ \norm{\Delta W_m}^2_{\mathcal{K}}|\mathcal{G}_{t^m}\right]   \right]\\
	& = \Delta t \tr Q \E\left[\norm{G(u_h^m)}^2_{\mathcal{L}(\mathcal{K},L^2)}\norm{u^m_h}_{L^2}^2\right]\\
		& \leq  C\Delta t \tr Q K_1^2 \left(\E\norm{u^m_h}_{L^2}^4+1\right),
	\end{split}
\end{align} 
and 
\begin{align}\label{eq:no2}
	\begin{split}
	\E\norm{G(u_h^m)}^4_{\mathcal{L}(\mathcal{K},L^2)}\norm{\Delta W_m}^4_{\mathcal{K}}& = \E\left[\norm{G(u_h^m)}^4_{\mathcal{L}(\mathcal{K},L^2)} \E\left[\norm{\Delta W_m}^4_{\mathcal{K}}|\mathcal{G}_{t^m}\right]\right]\\
	& \leq C (\tr Q)^2 \Delta t^2 \E\norm{G(u_h^m)}^4_{\mathcal{L}(\mathcal{K},L^2)}\\
	& \leq C K_1^4  (\tr Q)^2 \Delta t^2\left( \E\norm{ u_h^m}^4_{L^2}+1\right),
	\end{split}
\end{align}
using~\cite[Proposition 2.19]{DaPrato2014}.
Thus, we obtain,
\begin{multline}\label{eq:provingiii3}
	\frac14\left(\E\norm{u^{m+1}_h}_{L^2}^4-\E\norm{u^m_h}_{L^2}^4\right) + \frac18\E\left| \norm{u^{m+1}_h}_{L^2}^2-\norm{u^m_h}_{L^2}^2\right|^2 +\frac12\E\norm{u^{m+1}_h-\tu^{m+1}_h}_{L^2}^2\norm{u^{m+1}_h}_{L^2}^2\\ +\frac14\E\norm{\tu_h^{m+1} - u^m_h}_{L^2}^2\norm{u^{m+1}_h}_{L^2}^2 + \frac{\nu\Delta t}{2} \E\norm{\Grad \tu_h^{m+1}}_{L^2}^2\norm{u^{m+1}_h}_{L^2}^2\\
	\leq \frac{C}{\nu^3}\Delta t \norm{f^m}_{H^{-1}}^4 +5 C\Delta t \tr Q K_1^2 \left(\E\norm{u^m_h}_{L^2}^4+1\right)+ 4C K_1^4  (\tr Q)^2 \Delta t^2 \left(\E\norm{ u_h^m}^4_{L^2}+1\right).
\end{multline}
We sum over $m$ and use the discrete Gr\"onwall inequality,~\cite[Proposition 4.1]{Emmrich1999},  to obtain
\begin{multline}
	\label{eq:firstmomentestimate}
	\max_{1\leq m\leq N} \E\norm{u_h^m}_{L^2}^4 +\sum_{m=0}^{N-1} \E\left| \norm{u^{m+1}_h}_{L^2}^2-\norm{u^m_h}_{L^2}^2\right|^2 +\sum_{m=0}^{N-1} \E\norm{u^{m+1}_h-\tu^{m+1}_h}_{L^2}^2\norm{u^{m+1}_h}_{L^2}^2\\
	+\sum_{m=0}^{N-1} \E\norm{\tu_h^{m+1} - u^m_h}_{L^2}^2\norm{u^{m+1}_h}_{L^2}^2 
	+\nu \Delta t\sum_{m=0}^{N-1} \E\norm{\Grad \tu_h^{m+1}}_{L^2}^2\norm{u^{m+1}_h}_{L^2}^2\leq C(T,u_0,f,\tr Q).
\end{multline}
To prove estimate (iv) for $r=2$, we use the same strategy as the one we used to prove (ii):
We start by summing~\eqref{eq:provingiii2} over $m$ and then take the maximum and expectations
\begin{multline}\label{eq:provingiv1}
	\frac14\E\max_{1\leq i\leq N}\norm{u^{i}_h}_{L^2}^4
	\leq \frac14\norm{u^0_h}_{L^2}^4+\frac{C}{\nu^3}\Delta t \sum_{m=0}^{N-1}\norm{f^m}_{H^{-1}}^4 +5\sum_{m=0}^{N-1}\E\norm{G(u_h^m)}^2_{\mathcal{L}(\mathcal{K},L^2)}\norm{\Delta W_m}^2_{\mathcal{K}}\norm{u^m_h}_{L^2}^2 \\
	+ 4\sum_{m=0}^{N-1}\E\norm{G(u_h^m)}^4_{\mathcal{L}(\mathcal{K},L^2)}\norm{\Delta W_m}^4_{\mathcal{K}}
	+ \E\max_{1\leq i\leq N}\sum_{m=0}^{i-1}(G(u_h^m)\Delta W_m,u^{m}_h) \norm{u^m_h}_{L^2}^2.  
\end{multline}
We have already bounded the third and fourth term on the right-hand side in~\eqref{eq:no1} and~\eqref{eq:no2}. The first and second term are also bounded by the assumptions. So it remains to estimate the fifth term, for which we apply the Burkholder--Davis--Gundy inequality~\eqref{eq:BDGscalar} to the martingale
\begin{equation*}
	M(t) = \int_0^t \norm{u_h(s)}_{L^2}^2 (G(u_h(s))dW(s),u_h(s)),
\end{equation*}
and then Young's and Jensen's inequality,
\begin{align}\label{eq:martingaleterm1}
	 \E\max_{1\leq i\leq N}\sum_{m=0}^{i-1}(G(u_h^m)\Delta W_m,u^{m}_h) \norm{u^m_h}_{L^2}^2 & \leq C \E\left(\Delta t \sum_{m=0}^{N-1} \norm{G(u^m_h)}_{\mathcal{L}(\mathcal{K},L^2)}^2 \norm{u^m_h}_{L^2}^6\right)^{1/2}  \\
	 & \leq C \E\max_{0\leq i\leq N-1}\norm{u^i_h}_{L^2}^3\left(\Delta t \sum_{m=0}^{N-1}\norm{G(u_h^m)}_{\mathcal{L}(\mathcal{K},L^2)}^2\right)^{1/2}\notag\\
	 & \leq \frac18 \E \max_{0\leq i\leq N-1}\norm{u^i_h}_{L^2}^4 + C \Delta t^2 N \E\sum_{m=0}^{N-1}\norm{G(u_h^m)}_{\mathcal{L}(\mathcal{K},L^2)}^4\notag\\
	 &  \leq \frac18 \E \max_{1\leq i\leq N}\norm{u^i_h}_{L^2}^4 + \frac18 \norm{u_h^0}_{L^2}^4 + C K_1^4\Delta t T \E\sum_{m=0}^{N-1}\left(\norm{ u_h^m}_{L^2}^4+1\right).\notag
\end{align}
The first term on the right-hand side can be compensated by the left-hand side of~\eqref{eq:provingiv1} and the other two terms are bounded thanks to the previous estimate (iii). 
Thus (iv) follows for $r=2$ and $u^m_h$. 
To obtain the estimates for $r\geq 3$, we multiply~\eqref{eq:provingiii} successively by $\norm{u^{m+1}_h}_{L^2}^{2^{r-1}}$ which yields, in the same way as for $r=2$, 
\begin{multline}
	\E\max_{1\leq i \leq N}\norm{u^i_h}_{L^2}^{2^r} + \E\sum_{m=0}^{N-1}\left|\norm{u^{m+1}_h}_{L^2}^{2^{r-1}}-\norm{u^{m}_h}_{L^2}^{2^{r-1}}\right|^2+\E\sum_{m=0}^{N-1}\norm{u_h^{m+1}}_{L^2}^{2^r-2}\norm{\tu_h^{m+1}-u^m_h}_{L^2}^2\\
	+\E\sum_{m=0}^{N-1}\norm{u_h^{m+1}}_{L^2}^{2^r-2}\norm{\tu_h^{m+1}-u^{m+1}_h}_{L^2}^2 +
	\nu\Delta t\E \sum_{m=0}^{N-1}\norm{u_h^{m+1}}_{L^2}^{2^r-2}\norm{\Grad \tu_h^{m+1}}_{L^2}^2
	\leq C(T,f,u_0,\tr Q)
\end{multline}
thus proving (iii) and (iv) for general $r$ and $u^m_h$. 
We start by proving estimate (v) for $r=2$. To do so, we start with identity~\eqref{eq:energy1astep}, sum it from $m=0$ to $m=N-1$ and then square the resulting identity. We then use that
\begin{equation}
	\label{eq:basicstuff}
	\sum_{i=1}^I a_i^2 \leq \left(\sum_{i=1}^I a_i\right)^2 \leq I \sum_{i=1}^I a_i^2
\end{equation}
 holds for $a_i\geq 0$, and obtain
\begin{align}\label{eq:gradientbound1}
	\begin{split}
	&	\frac{1}{4}\norm{u^{N}_h}_{L^2}^4 +\frac14\left(\sum_{m=0}^{N-1}\norm{\tu^{m+1}_h - u^m_h}_{L^2}^2\right)^2+\frac14\left(\sum_{m=0}^{N-1}\norm{\tu^{m+1}_h - u^{m+1}_h}_{L^2}^2\right)^2 +\left(\Delta t \nu\sum_{m=0}^{N-1} \norm{\Grad \tu^{m+1}_h}_{L^2}^2\right)^2 \\
&	\qquad \leq 2 \norm{u^0_h}_{L^2}^4+4\left( \Delta t \sum_{m=0}^{N-1}\langle f^m,\tu^{m+1}_h\rangle\right)^2 \\
&\qquad 	+ 4\left(\sum_{m=0}^{N-1}(\PUh G(u^m_h)\Delta W_{m},\tu_h^{m+1}-u_h^m)\right)^2+4\left(\sum_{m=0}^{N-1}(\PUh G(u^m_h)\Delta W_{m},u_h^m)\right)^2\\
	&\qquad \leq 2 \norm{u^0_h}_{L^2}^4+8\left( \frac{\Delta t}{\nu} \sum_{m=0}^{N-1} \norm{ f^m}_{H^{-1}}^2\right)^2 +\frac{1}{2}
\left(\Delta t\nu \sum_{m=0}^{N-1}\norm{\Grad\tu^{m+1}_h}_{L^2}^2\right)^2 \\
&\qquad 	+ 32\left(\sum_{m=0}^{N-1}\norm{G(u^m_h)}_{\mathcal{L}(\mathcal{K},L^2)}^2\norm{\Delta W_{m}}_{\mathcal{K}}^2\right)^2 +\frac18 \left(\sum_{m=0}^{N-1}\norm{\tu_h^{m+1}-u_h^m}_{L^2}^2\right)^2+4\left(\sum_{m=0}^{N-1}( G(u^m_h)\Delta W_{m},u_h^m)\right)^2.
\end{split}
\end{align}
Rearranging terms and then taking expectations, we obtain
\begin{multline}\label{eq:gradientbound2}
		\frac{1}{4}\E\norm{u^{N}_h}_{L^2}^4 +\frac18\E\left(\sum_{m=0}^{N-1}\norm{\tu^{m+1}_h - u^m_h}_{L^2}^2\right)^2+\frac14\E\left(\sum_{m=0}^{N-1}\norm{\tu^{m+1}_h - u^{m+1}_h}_{L^2}^2\right)^2 +\frac12\E\left(\Delta t \nu\sum_{m=0}^{N-1} \norm{\Grad \tu^{m+1}_h}_{L^2}^2\right)^2 \\
		\leq 2  \norm{u^0_h}_{L^2}^4+8\left( \frac{\Delta t}{\nu} \sum_{m=0}^{N-1} \norm{ f^m}_{H^{-1}}^2\right)^2  \\
		\qquad 	+ 32\E\left(\sum_{m=0}^{N-1}\norm{G(u^m_h)}_{\mathcal{L}(\mathcal{K},L^2)}^2\norm{\Delta W_{m}}_{\mathcal{K}}^2\right)^2 + 4\E\left(\sum_{m=0}^{N-1}( G(u^m_h)\Delta W_{m},u_h^m)\right)^2.
\end{multline}
The first two terms on the right-hand side are bounded, so it remains to estimate the other two terms. For the third term, we use the law of iterated expectations and~\eqref{eq:basicstuff} and~\cite[Prop. 2.19]{DaPrato2014} again,
\begin{align*}
	\E\left(\sum_{m=0}^{N-1}\norm{G(u^m_h)}_{\mathcal{L}(\mathcal{K},L^2)}^2\norm{\Delta W_{m}}_{\mathcal{K}}^2\right)^2 & \leq N \sum_{m=0}^{N-1}\E\norm{G(u^m_h)}_{\mathcal{L}(\mathcal{K},L^2)}^4\norm{\Delta W_{m}}_{\mathcal{K}}^4\\
	& = N \sum_{m=0}^{N-1}\E\left[\norm{G(u^m_h)}_{\mathcal{L}(\mathcal{K},L^2)}^4\E\left[\norm{\Delta W_{m}}_{\mathcal{K}}^4|\mathcal{G}_{t^m}\right]\right]\\
	& \leq C\Delta t^2  N (\tr Q)^2\sum_{m=0}^{N-1}\E\left[\norm{G(u^m_h)}_{\mathcal{L}(\mathcal{K},L^2)}^4\right]\\
	& \leq C\Delta t   T K_1^4 (\tr Q)^2\sum_{m=0}^{N-1}\E\left[\norm{u^m_h}_{L^2}^4+1\right]\leq C(T,u_0,f,\tr Q),
\end{align*}
where we used (iii) and~\cite[Proposition 2.19]{DaPrato2014} to bound the right-hand side. For the fourth term in~\eqref{eq:gradientbound2}, we use the Burkholder--Davis--Gundy inequality~\eqref{eq:BDGscalar} once more:
\begin{align*}
	\E\left(\sum_{m=0}^{N-1}( G(u^m_h)\Delta W_{m},u_h^m)\right)^2 & \leq C \Delta t \E\sum_{m=0}^{N-1}\norm{G(u_h^m)}_{\mathcal{L}(\mathcal{K},L^2)}^2 \norm{u_h^m}_{L^2}^2\\
	& \leq C K_1^2\Delta t \sum_{m=0}^{N-1}\left(\E\norm{u_h^m}_{L^2}^4+1\right)\leq C(T,u_0,f,\tr Q),
\end{align*}
using (iii) once more for $r=2$. Thus the right-hand side of~\eqref{eq:gradientbound2} is bounded and (v) follows for $r=2$. For $r=3$, we square~\eqref{eq:gradientbound1} (after rearranging terms), and use again~\eqref{eq:basicstuff}:
\begin{multline}\label{eq:gradientbound3}
	\frac{1}{16} \norm{u^{N}_h}_{L^2}^8 +\frac{1}{64} \left(\sum_{m=0}^{N-1}\norm{\tu^{m+1}_h - u^m_h}_{L^2}^2\right)^4+\frac{1}{16} \left(\sum_{m=0}^{N-1}\norm{\tu^{m+1}_h - u^{m+1}_h}_{L^2}^2\right)^4 +\frac14 \left(\Delta t \nu\sum_{m=0}^{N-1} \norm{\Grad \tu^{m+1}_h}_{L^2}^2\right)^4 \\
	\leq 16 \norm{u^0_h}_{L^2}^8+256\left( \frac{\Delta t}{\nu} \sum_{m=0}^{N-1} \norm{ f^m}_{H^{-1}}^2\right)^4  \\
	\qquad 	+ 4096\left(\sum_{m=0}^{N-1}\norm{G(u^m_h)}_{\mathcal{L}(\mathcal{K},L^2)}^2\norm{\Delta W_{m}}_{\mathcal{K}}^2\right)^4 + 64\left(\sum_{m=0}^{N-1}(\PUh G(u^m_h)\Delta W_{m},u_h^m)\right)^4,
\end{multline}
then we take expectations and bound the terms on the right-hand side in the same way as for $r=2$ using the already proved claim (iii). For $r>3$, we successively square~\eqref{eq:gradientbound3}, take expectations and bound the right-hand side using (iii).
Now the estimates for $\tu_h^m$ follow: 
From~\eqref{eq:projectionestimate}, we have
\begin{equation}\label{eq:tuh}
	\norm{\tu_h^{m}}_{L^2}^q \leq C_q \norm{u_h^m}_{L^2}^q + C_q \norm{u^m_h-\tu_h^m}_{L^2}^q\leq C_q \norm{u_h^m}_{L^2}^q + C_q\left(\sum_{j=1}^N  \norm{u^j_h-\tu_h^j}_{L^2}^2\right)^{q/2}
\end{equation}
Taking expectations and using (iii) for $u_h^m$ and (v), the right-hand side is bounded. Therefore the second part of (iii) follows for any $1\leq q\leq 2^{r_0}$.
 For (iv), we first take the maximum over $m$ in~\eqref{eq:tuh} and then take expectations and instead use (iv) for $u_h^m$ and (v).
\end{proof}
\section{Convergence of the scheme}\label{sec:conv}
We start by defining piecewise constant functions based on the approximations computed by the scheme~\eqref{eq:predictionfully} - \eqref{eq:projectionfully}: We denote
	\begin{alignat}{2}
	u_h(t,x,\omega)& = u_h^m(x,\omega),\qquad &t\in [t^m,t^{m+1}),\quad x\in \dom,\quad \omega\in\Omega,\label{eq:defuh}\\
	\tu_h(t,x,\omega) &= \tu_h^{m+1}(x,\omega),\qquad &t\in [t^m,t^{m+1}),\quad x\in \dom,\quad \omega\in\Omega,\\
	p_h(t,x,\omega)& = p_h^{m}(x,\omega),\qquad & t\in [t^m,t^{m+1}),\quad x\in \dom,\quad \omega\in\Omega,\\
	f_h(t,x)& = \frac{1}{\Delta t}\int_{t^m}^{t^{m+1}} f(s,x)ds,\qquad &t\in [t^m,t^{m+1}),\quad x\in \dom,\label{eq:defph}
\end{alignat}
for $m=0,1,2,\dots,N-1$. We also denote for $x\in \dom$, $\omega\in \Omega$, $r\geq T$,
\begin{equation*}
	u_h(r,x,\omega) = u^N_h(x,\omega),\quad \tu_h(r,x,\omega) = \tu^N_h(x,\omega),\quad p_h(r,x,\omega) = p^N_h(x,\omega),
\end{equation*}
and for $x\in \dom$, $\omega\in \Omega$ and $r< 0$,
\begin{equation*}
	u_h(r,x,\omega) = u^0_h(x),\quad \tu_h(r,x,\omega) = \tu^0_h(x),\quad p_h(r,x,\omega) = p^0_h(x).
\end{equation*}
$u_h$ and $p_h$ are adapted to the filtration $(\mathcal{G}_t)_{t\in [0,T]}$ thanks to Lemma~\ref{lem:solvability}. 
Throughout this section, we will consider a sequence of discretization parameters $(h_n,\Delta t_n)_{n\in \N}\to 0$ and subsequences thereof which satisfy certain coupling conditions stated later on. For the ease of notation, we will omit writing the subscript $n$ and `for all $(h,\Delta t)\leq (h_0,\Delta t_0)$' means `for all $n\geq n_0$ with $n_0$ such that $h_n\leq h_0$ and $\Delta t_n \leq \Delta t_0$ for $n\geq n_0$'.
To show convergence of the scheme to a martingale solution of the stochastic Navier--Stokes equations, we first verify Aldous' condition and combine it with the discrete energy estimates to establish tightness of the laws of $u_h$ on the Skorokhod space $\mathbb D([0,T];H^{-1}(\dom))$. The difference $\tu_h-u_h$ vanishes as $h,\Delta t\to 0$ by the energy estimate, Lemma~\ref{lem:energyestimate}. We can then apply
Skorokhod's representation theorem to deduce the existence of copies of the approximations on a new probability space. Using the a priori energy estimate, Lemma~\ref{lem:energyestimate}, we obtain convergence to a common limit $\bar u$.
We establish the divergence constraint and time continuity of this limit and pass to the limit in the non-stochastic integrals of the scheme.
Finally, the martingale representation theorem will provide a
$Q$-Wiener process on an extension of the probability space,
with respect to which $\bar u$ is a martingale solution
of~\eqref{eq:SNSm}.

\subsection{Aldous' condition}\label{sec:aldous}
To show the tightness, we will use the Aldous condition combined with the energy estimates from Lemma~\ref{lem:energyestimate}. Before we can state the tightness criteria, we need to introduce some notation: For $(\mathbb{S},d)$ a separable and complete metric space, we let $\mathbb{D}([0,T];\mathbb{S})$ denote  the space of all $\mathbb{S}$-valued c\`adl\`ag functions defined on $[0,T]$, i.e., the functions which are right continuous and have left limits at every $t\in [0,T]$. If $u\in \mathbb{D}([0,T];\mathbb{S})$, we denote by $u(t)$ its value at time $t$ and by $u(t-)$ its left-hand limit at time $t$.
We endow $\mathbb{D}([0,T];\mathbb{S})$ with the Skorokhod topology~\cite{Joffe1986}. This means that a sequence $\{u_n\}\subset \mathbb{D}([0,T];\mathbb{S})$ converges to $u\in \mathbb{D}([0,T];\mathbb{S})$ iff there exists a sequence $\{\lambda_n\}$ of homeomorphisms of $[0,T]$ such that $\lambda_n(s)\to s$ uniformly on $[0,T]$ and $u_n\circ \lambda_n\to u$ uniformly on $[0,T]$. This topology is metrizable with the metric $d_T$ defined by~\cite{Motyl2013,Joffe1986}
\begin{equation*}
	d_T(u,v):= \inf_{\lambda \in \Lambda_T}\left(\sup_{t\in [0,T]} d(u(t),v\circ \lambda(t)) + \sup_{t\in [0,T]} |t-\lambda(t)| + \sup_{s\neq t}\left|\log \left(\frac{\lambda(t)-\lambda(s)}{t-s}\right)\right|\right),
\end{equation*}
where $\Lambda_T$ is the set of increasing homeomorphisms of $[0,T]$. $(\mathbb{D}([0,T];\mathbb{S}),d_T)$ is a Polish space~\cite{Joffe1986}.

\begin{definition}[Aldous condition]
	\label{def:aldous}
	Let $(X_n)_{n\in \N}$ be a sequence of $\mathbb{S}$-valued c\`adl\`ag processes and for each $n$ let $(\mathcal{G}^n_t)_{t\in [0,T]}$ be a filtration to which $X_n$ is adapted. Then $(X_n)$  satisfies the Aldous condition iff
	\begin{equation}
		\label{eq:aldous}
		\begin{split}
		&	\forall \, \epsilon>0, \quad \forall \,\eta>0,\quad \exists\, \delta>0,\quad \exists\, n_0\in\N,\quad 
		\text{such that for every sequence }\, (\tau_n)_{n\in \N}\, \\
		&\text{ of $\mathcal{G}^n$-stopping times with } \, \tau_n\leq T\, \text{one has }\\
		&\sup_{n\geq n_0}\sup_{0\leq \theta\leq \delta} \mathbb{P}\left[d(X_n(\tau_n+\theta),X_n(\tau_n))\geq \eta\right]\leq \epsilon,
		\end{split}
	\end{equation}
	where we let $X_n(t)=X_n(T)$ for $t>T$.
\end{definition}
We will now show that the approximations $u_h$ computed by the scheme~\eqref{eq:predictionfully} - \eqref{eq:projectionfully} satisfy the Aldous condition with $\mathbb{S} = H^{-1}(\dom)$. First we need the following technical lemma:
\begin{restatable}{lemma}{rellichsmooth}
	\label{lem:Rellich3}
	For any $\zeta>0$ there exists $L\in \N$ and $\phi_1,\dots, \phi_L\in C^\infty_{c,\Div}(\dom)\oplus \Grad C^\infty(\overline{\dom})$ such that for any $v\in L^2(\dom)$, we have
	\begin{equation}
		\label{eq:rellich2}
		\norm{v}_{H^{-1}(\dom)}\leq \sqrt{\sum_{k=1}^L (v,\phi_k)^2} + \zeta \norm{v}_{L^2(\dom)}.
	\end{equation}	
\end{restatable}
Using this, we can prove:
\begin{lemma}
	\label{lem:aldous}
	Under Assumption~\ref{ass:stokes} and the inverse CFL-type condition $h^{2(k+s)}\leq C \Delta t$, where $k$ is the approximation order of the space $\Uh$,
 for any $\epsilon>0$, there exist $\delta>0$ and $(h_0,\Delta t_0)>0$ such that $\forall \tau_h\leq T$ $(\mathcal{G}_t)$-stopping times, we have 
\begin{equation}\label{eq:aldousuh}
	\sup_{0 < (h,\Delta t) \leq (h_0,\Delta t_0)}\sup_{0\leq \theta\leq \delta} \E \norm{u_h(\tau_h+\theta)-u_h(\tau_h)}^\alpha_{H^{-1}(\dom)}\leq  \epsilon,
\end{equation}
where $1\leq \alpha\leq 2^{r_0-1}$. 
\end{lemma}
This lemma implies that the Aldous condition from Definition~\ref{def:aldous} is satisfied with $\mathcal{G}^n_t = \mathcal{G}_t$ for all $n$ after using Markov's inequality.
\begin{proof}
	First, we note that by Lemma~\ref{lem:Rellich3}, for any given $\zeta>0$, we can find $L\in \N$ and $\phi_1,\dots, \phi_L\in C^\infty_{c,\Div}(\dom)\oplus \Grad C^\infty(\overline{\dom})$ such that
	\begin{equation*}
			\E\norm{u_h(\tau_h+\theta)-u_h(\tau_h)}_{H^{-1}(\dom)}^\alpha\leq C_\alpha\E \left(\sum_{j=1}^L (u_h(\tau_h+\theta)-u_h(\tau_h),\phi_j)^2\right)^{\alpha/2} + \zeta^\alpha C_\alpha\E \norm{u_h(\tau_h+\theta)-u_h(\tau_h)}_{L^2(\dom)}^\alpha.
	\end{equation*}
	By Lemma~\ref{lem:energyestimate}, (iv), we can bound the last right-hand side term,
		\begin{equation}\label{eq:halftheestimate}
		\E\norm{u_h(\tau_h+\theta)-u_h(\tau_h)}_{H^{-1}(\dom)}^\alpha\leq C_\alpha\E \left(\sum_{j=1}^L (u_h(\tau_h+\theta)-u_h(\tau_h),\phi_j)^2\right)^{\alpha/2} + \zeta^\alpha C_\alpha C(T,u_0,f,\tr Q).
	\end{equation}
	Thus we can choose $\zeta^\alpha< \epsilon (2C_\alpha C(T,u_0,f,\tr Q))^{-1}$ and it remains to estimate $\E \left(\sum_{j=1}^L (u_h(\tau_h+\theta)-u_h(\tau_h),\phi_j)^2\right)^{\alpha/2}$.
	So we fix $\omega\in \Omega$. For a given $\theta>0$ and a stopping time $\tau_h$, we can write
	\begin{equation}\label{eq:telescopic}
		u_h(\tau_h+\theta)-u_h(\tau_h)=u_h^{m_2+1}-u_h^{m_1} = \sum_{m=m_1}^{m_2} (u_h^{m+1}-u_h^m),
	\end{equation}
where 
\begin{equation*}
	m_1:=\Big\lfloor\frac{\tau_h}{\Delta t}\Big\rfloor,\qquad m_2:=\Big\lfloor\frac{(\tau_h+\theta)\wedge T}{\Delta t}\Big\rfloor-1.
\end{equation*}
Without loss of generality, we assume that $m_2\geq m_1$. 
	Now given a function $v\in C^\infty_{c,\Div}(\dom)\oplus \Grad C^\infty(\overline{\dom})$, we write  $v= w +\Grad g$ where $\Div w= 0$ and $g\in C^\infty(\overline{\dom})$. So we have
	\begin{align}
		\label{eq:decomp1}
(u_h(\tau_h+\theta)-u_h(\tau_h),v)&= (u_h(\tau_h+\theta)-u_h(\tau_h),\Ih w) + (u_h(\tau_h+\theta)-u_h(\tau_h), w - \Ih w)\\
&\qquad + (u_h(\tau_h+\theta)-u_h(\tau_h),\Grad g)\notag \\
	& = -\underbrace{\Delta t\sum_{m=m_1}^{m_2} b(\tu_h^m,\tu_h^{m+1},\Ih  w)}_{\text{(A.I)}} - \underbrace{\nu\Delta t \sum_{m=m_1}^{m_2} (\Grad \tu_h^{m+1},\Grad \Ih w)}_{\text{(A.II)}}\notag\\
	& +\underbrace{\Delta t \sum_{m=m_1}^{m_2} ( p_h^{m+1},\Div \Ih w)}_{\text{(A.III)}} + \underbrace{\Delta t \sum_{m=m_1}^{m_2}\langle f^m,\Ih w\rangle }_{\text{(A.IV)}}\notag\\
	&	+ \underbrace{\sum_{m=m_1}^{m_2}(G(u_h^m)\Delta W_m,\Ih w)}_{\text{(A.V)}}+ \underbrace{(u_h(\tau_h+\theta)-u_h(\tau_h), \Grad g )}_{\text{(B)}} \notag\\
	&+ \underbrace{(u_h(\tau_h+\theta)-u_h(\tau_h),w-\Ih w)}_{\text{(C)}},\notag
	\end{align}
	where $\Ih: H^2(\dom)\to \Uh$ is the operator from Assumption~\ref{ass:approxorder}.
	We estimate each of the terms:
	\begin{align*}
		|\text{(A.I)}| & \leq \Delta t\sum_{m=m_1}^{m_2} \norm{\tu_h^{m+1}}_{H^1} \norm{\tu_h^m}_{L^2} \norm{\Ih w}_{W^{1,\infty}}\\
		& \leq C \Delta t\sum_{m=m_1}^{m_2} \norm{\Grad \tu_h^{m+1}}_{L^2} \norm{\tu_h^m}_{L^2} \norm{ w}_{W^{1,\infty}}\\
			& \leq C \Delta t\sum_{m=m_1}^{m_2} \norm{\Grad \tu_h^{m+1}}_{L^2} \norm{\tu_h^m}_{L^2} \norm{w}_{H^3}\\
			& \leq C \sqrt{(m_2-m_1+1)\Delta t}\left(\Delta t\sum_{m=1}^N\norm{\Grad \tu^m_h}_{L^2(\dom)}^2\right)^{1/2}\norm{\tu_h}_{L^\infty(0,T;L^2(\dom))}\norm{w}_{H^3},
	\end{align*}
	where we used the Poincar\'e inequality and a standard interpolation estimate. We continue to the next term,
	\begin{align*}
		|\text{(A.II)}| & \leq \nu \Delta t \sum_{m=m_1}^{m_2} \norm{\Grad \tu_h^{m+1}}_{L^2}\norm{\Grad \Ih w}_{L^2}\\
		& \leq C \sqrt{(m_2-m_1+1)\Delta t}\left(\Delta t\sum_{m=1}^N\norm{\Grad \tu^m_h}_{L^2(\dom)}^2\right)^{1/2}\norm{w}_{H^2}.
	\end{align*}
	Next, we let $(w_h,\pi_h)=\Sh(w,0)$ be the Stokes projection of $w$. Then $w_h$ is discretely divergence free and therefore $(p_h,\Div w_h)=0$. By~\eqref{eq:projection1}, this implies
	\begin{equation}\label{eq:rewrite}
		(p_h^{m+1},\Div w_h)= \left(\frac{u^{m+1}_h - \tu_h^{m+1}}{\Delta t},w_h\right) = 0.
	\end{equation}
	 We use this in term (A.III) together with Lemma~\ref{lem:stokesestimate}:
		\begin{align}\label{eq:LBBpressureestimate}
		|\text{(A.III)}| & \leq \Delta t \left|\sum_{m=m_1}^{m_2} ( p^{m+1}_h,\Div (\Ih w- w_h))\right|\\
		& = \Delta t \left|\sum_{m=m_1}^{m_2}  \left(\frac{u^{m+1}_h - \tu_h^{m+1}}{\Delta t},\Ih w- w_h\right) \right|\notag\\
		& \leq  \sum_{m=m_1}^{m_2}  \norm{u^{m+1}_h - \tu_h^{m+1}}_{L^2}\norm{\Ih w-w_h}_{L^2}\notag\\
		& \leq \sqrt{m_2-m_1 +1} \left( \sum_{m=m_1}^{m_2}  \norm{u^{m+1}_h - \tu_h^{m+1}}_{L^2}^2\right)^{1/2}\left(\norm{\Ih w-w}_{L^2}+\norm{w-w_h}_{L^2}\right)\notag\\
		& \leq \frac{\sqrt{(m_2-m_1 +1)\Delta t}}{\sqrt{\Delta t}} \left( \sum_{m=m_1}^{m_2}  \norm{u^{m+1}_h - \tu_h^{m+1}}_{L^2}^2\right)^{1/2}\left(h^{k+1}\norm{w}_{H^{k+1}}+h^{s+k}\norm{w}_{H^{k+1}}\right)\notag
	\end{align}
	Thus, we can estimate (A.III) as follows using the inverse CFL-type condition $h^{2(k+s)}\leq C\Delta t$:
		\begin{align*}
		|\text{(A.III)}| & \leq  \frac{\sqrt{(m_2-m_1 +1)\Delta t}}{\sqrt{\Delta t}} \left( \sum_{m=m_1}^{m_2}  \norm{u^{m+1}_h - \tu_h^{m+1}}_{L^2}^2\right)^{1/2}\left(h^{k+1}\norm{w}_{H^{k+1}}+h^{s+k}\norm{w}_{H^{k+1}}\right)\\
		& \leq \sqrt{(m_2-m_1 +1)\Delta t}\left( \sum_{m=m_1}^{m_2}  \norm{u^{m+1}_h - \tu_h^{m+1}}_{L^2}^2\right)^{1/2}\norm{w}_{H^{k+1}}.		
	\end{align*}
	We continue to estimate term (A.IV):
	\begin{align*}
		|\text{(A.IV)}|& \leq \Delta t \sum_{m=m_1}^{m_2} \norm{f^m}_{H^{-1}}\norm{\Ih w}_{H^1}\\
		& \leq C \Delta t \sum_{m=m_1}^{m_2} \norm{f^m}_{H^{-1}}\norm{ w}_{H^2}\\
		& \leq C \sqrt{(m_2-m_1+1)\Delta t} \norm{f}_{L^2(0,T;H^{-1})}\norm{w}_{H^2},
	\end{align*}
	where we used interpolation estimates again. We estimate the term involving the stochastic forcing:
	\begin{align*}
		|\text{(A.V)}|& \leq \norm{\sum_{m=m_1}^{m_2}G(u_h^m)\Delta W_m}_{L^2}\norm{\Ih w}_{L^2}\\
		& \leq C\norm{\sum_{m=m_1}^{m_2}G(u_h^m)\Delta W_m}_{L^2}\norm{w}_{H^2}.
	\end{align*}
	For the term (C), we use the properties of the interpolation operator,~\cite[Theorem 11.13]{Ern2021}:
	\begin{align}\label{eq:termC}
		|\text{(C)}|& \leq \left(\norm{u_h(\tau_h+\theta)}_{L^2}+\norm{u_h(\tau_h)}_{L^2}\right) \norm{w-\Ih w}_{L^2}\\
		& \leq 
		C \left(\norm{u_h(\tau_h+\theta)}_{L^2}+\norm{u_h(\tau_h)}_{L^2}\right)h^{k+1} \norm{w}_{H^{k+1}}.\notag
	\end{align}
	In order to estimate term (B), we split it up into two terms  (that can be bounded in the same way) using triangle inequality and then use the scheme again:
	\begin{equation*}
		|\text{(B)}| \leq \left|(u_h(\tau_h+\theta), \Grad g )\right| + \left|(u_h(\tau_h), \Grad g )\right|.
	\end{equation*}
	We consider the first term, the second is estimated in the same way. We have $u_h(\tau_h+\theta) = u_h^{m}$ for some $m\in \N$. Then we use that $u_h^m$ is weakly divergence free,
		\begin{align}\label{eq:termB}
	\left|(u_h(\tau_h+\theta), \Grad g )\right| &\leq |(u^m_h,\Grad g-\Grad \Ih^p g)|+|(u^m_h,\Grad \Ih^p g)|\\
		& \leq  C\norm{u_h}_{L^\infty(0,T;L^2(\dom))}h^{\ell} \norm{g}_{H^{\ell+1}},\notag
	\end{align} 
	where  we  denoted by $\Ih^p$ the  interpolation operator on the space $\Ph$ as in Assumption~\ref{ass:approxorder} (cf.~\cite[Section 9.1]{Ern2021}).  
	Next, we denote $\sigma:= \max\{3,k+1\}$. Combining all these estimates and using that $(m_2-m_1+1)\Delta t\leq \theta +\Delta t$
		\begin{align*}
		\left|\left(u_h(\tau_h+\theta)-u_h(\tau_h),v\right)\right|&   \leq C \norm{w}_{H^\sigma}\Bigg((\sqrt{\theta} +\sqrt{\Delta t})\norm{ \tu_h}_{L^2([0,T];H^1(\dom))}\left(\norm{\tu_h}_{L^\infty(0,T;L^2(\dom))}+1\right)\\
		& + (\sqrt{\theta} +\sqrt{\Delta t})\left(\left(\sum_{m=m_1}^{m_2}  \norm{u^{m+1}_h-\tu^{m+1}_h}_{L^2}^2\right)^{1/2} + \norm{f}_{L^2(0,T;H^{-1})} \right) \\
		&+\norm{\sum_{m=m_1}^{m_2}G(u_h^m)\Delta W_m}_{L^2}\Bigg)+C \norm{u_h}_{L^\infty(0,T;L^2(\dom))}\left(h^\ell \norm{g}_{H^{\ell+1}}+h^{k+1}\norm{w}_{H^\sigma}\right).
	\end{align*}
	We choose $v=\phi_j=w_j + \Grad g_j$, $j=1,\dots, L$, (where $\Div w_j=0$), raise everything to the power 2, then sum over $j$ and raise to the power $\alpha/2$ and take expectations to obtain
		\begin{align*}
		&	\E \left(\sum_{j=1}^L (u_h(\tau_h+\theta)-u_h(\tau_h),\phi_j)^2\right)^{\alpha/2} \\
		& \leq C\left(\sum_{j=1}^L \norm{w_j}_{H^\sigma}^2\right)^{\alpha/2}\Bigg( ( {\theta}^{\alpha/2} +{\Delta t}^{\alpha/2})\underbrace{\E\left(\Delta t\sum_{m=1}^N\norm{\Grad \tu^m_h}_{L^2(\dom)}^2\right)^{\alpha/2}\left(\norm{\tu_h}^\alpha_{L^\infty(0,T;L^2(\dom))}+1\right)}_{(a)}\\
		& + ({\theta}^{\alpha/2} +{\Delta t}^{\alpha/2})\left(\underbrace{\E\left(\sum_{m=m_1}^{m_2}  \norm{u^{m+1}_h-\tu^{m+1}_h}_{L^2}^2\right)^{\alpha/2}}_{(b)} + \norm{f}^\alpha_{L^2(0,T;H^{-1})} \right) \\
		&+\underbrace{\E\norm{\sum_{m=m_1}^{m_2}G(u_h^m)\Delta W_m}^\alpha_{L^2}}_{(c)}\Bigg)\\
		&+C \underbrace{\E\norm{u_h}^\alpha_{L^\infty(0,T;L^2)}}_{(d)}\left(h^{\alpha \ell}\left(\sum_{j=1}^L \norm{g_j}_{H^{\ell+1}}^2\right)^{\alpha/2} +h^{(k+1)\alpha}\left(\sum_{j=1}^L \norm{w_j}_{H^{\sigma}}^2\right)^{\alpha/2} \right).
	\end{align*}
	We estimate these terms using the energy estimate, Lemma~\ref{lem:energyestimate}: Using Young's inequality, we have
	\begin{align*}
		(a)& \leq C\left(\E\norm{\tu_h}_{L^\infty(0,T;L^2)}^{2\alpha}  +\E \left(\Delta t\sum_{m=1}^N\norm{\Grad \tu^m_h}_{L^2(\dom)}^2\right)^{\alpha}+1\right)\leq C,
	\end{align*}
	using Lemma~\ref{lem:energyestimate}, (iv) and (v). Term (b) is bounded by Lemma~\ref{lem:energyestimate}, (v), and term  (d) is bounded by Lemma~\ref{lem:energyestimate} (iv). It remains to estimate term (c). We first define the martingale 
	\begin{equation*}
		M_h(t):= \int_0^t G(u_h(s))dW(s).
	\end{equation*}
	Then 
	\begin{align*}
		\sum_{m=m_1}^{m_2} G(u_h^m)\Delta W_m &= M_h(t^{m_2+1})-M_h(t^{m_1})\\
		& = M_h\left(\Delta t \left\lfloor\frac{(\tau_h+\theta)\wedge T}{\Delta t}\right\rfloor \right) - M_h \left(\Delta t\left\lfloor\frac{\tau_h}{\Delta t}\right\rfloor \right)\\
		& = M_h\left(\Delta t \left\lfloor\frac{(\tau_h+\theta)\wedge T}{\Delta t}\right\rfloor \right)-
		M_h\left(\Delta t \left\lceil\frac{(\tau_h+\theta)\wedge T}{\Delta t}\right\rceil \right)\\
		&\quad +M_h\left(\Delta t \left\lceil\frac{(\tau_h+\theta)\wedge T}{\Delta t}\right\rceil \right)-M_h \left(\Delta t\left\lceil\frac{\tau_h}{\Delta t}\right\rceil \right)\\
		&\quad +M_h \left(\Delta t\left\lceil\frac{\tau_h}{\Delta t}\right\rceil \right)
		 -  M_h \left(\Delta t\left\lfloor\frac{\tau_h}{\Delta t}\right\rfloor \right).
	\end{align*}
	We have for $\max\{2,\alpha\}<q\leq 2^{r_0}$,
	\begin{align*}
		&\E\norm{M_h\left(\Delta t \left\lfloor\frac{(\tau_h+\theta)\wedge T}{\Delta t}\right\rfloor \right)-
			M_h\left(\Delta t \left\lceil\frac{(\tau_h+\theta)\wedge T}{\Delta t}\right\rceil \right)}_{L^2}^\alpha \\
			& \leq \E \max_{0\leq m\leq N-1} \norm{G(u_h^m)\Delta W_m}_{L^2}^\alpha\\
			& \leq \left(\E\sum_{m=0}^{N-1} \norm{G(u_h^m)\Delta W_m}_{L^2(\dom)}^{q}\right)^{\frac{\alpha}{q}}\\
			& \leq C\left(\E\sum_{m=0}^{N-1} \norm{G(u_h^m)}^{q}_{\mathcal{L}(\mathcal{K},L^2)}\norm{\Delta W_m}_{\mathcal{K}}^{q}\right)^{\frac{\alpha}{q}}\\
			& \leq CK_1^\alpha \left(\E\sum_{m=0}^{N-1} \left(\norm{ u_h^m}^{q}_{L^2(\dom)}+1\right)\norm{\Delta W_m}_{\mathcal{K}}^{q}\right)^{\frac{\alpha}{q}}\\
			& \leq K_1^\alpha C_q \Delta t^{\frac{\alpha}{2}-\frac{\alpha}{q}}\left(\left(\E  \norm{ u_h }^{q}_{L^{q}(0,T;L^2(\dom))}\right)^{\frac{\alpha}{q}}+1\right)\\
			& \leq C \Delta t^{\frac{\alpha}{2}-\frac{\alpha}{q}}.
	\end{align*}
	Using the same estimate, we get
	\begin{equation*}
		\E\norm{M_h \left(\Delta t\left\lceil\frac{\tau_h}{\Delta t}\right\rceil \right)
		-  M_h \left(\Delta t\left\lfloor\frac{\tau_h}{\Delta t}\right\rfloor \right)}_{L^2}^\alpha\leq C \Delta t^{\frac{\alpha}{2}-\frac{\alpha}{q}}.
	\end{equation*}
	Now $\Delta t \left\lceil\frac{(\tau_h+\theta)\wedge T}{\Delta t}\right\rceil$ and $\Delta t\left\lceil\frac{\tau_h}{\Delta t}\right\rceil$ are both stopping times and satisfy
	\begin{equation*}
		0\leq \Delta t \left\lceil\frac{(\tau_h+\theta)\wedge T}{\Delta t}\right\rceil - \Delta t\left\lceil\frac{\tau_h}{\Delta t}\right\rceil \leq \theta +\Delta t,
	\end{equation*} 
	and hence we can use the Burkholder--Davis--Gundy inequality,~\eqref{eq:BDG}, to estimate
	\begin{align*}
		&\E\norm{M_h\left(\Delta t \left\lceil\frac{(\tau_h+\theta)\wedge T}{\Delta t}\right\rceil \right) - M_h \left(\Delta t\left\lceil\frac{\tau_h}{\Delta t}\right\rceil \right)}_{L^2}^\alpha\\
		& \leq C\E  \left(\int_{\Delta t\left\lceil\frac{\tau_h}{\Delta t}\right\rceil}^{\Delta t \left\lceil\frac{(\tau_h+\theta)\wedge T}{\Delta t}\right\rceil } \norm{G(u_h(s))}_{\mathcal{L}(\mathcal{K},L^2(\dom))}^2 ds \right)^{\alpha/2}\\
		& \leq C K_1^\alpha \E\left( \int_{\Delta t\left\lceil\frac{\tau_h}{\Delta t}\right\rceil}^{\Delta t \left\lceil\frac{(\tau_h+\theta)\wedge T}{\Delta t}\right\rceil }\left( \norm{u_h}_{L^2}^2  +1\right)ds  \right)^{\alpha/2}\\
		& \leq C K_1^\alpha (\theta^{\alpha/2}+\Delta t^{\alpha/2})\E\left(\norm{u_h}_{L^\infty(0,T;L^2(\dom))}^2+1\right)^{\alpha/2}.
	\end{align*}
	Thus with the energy estimate, Lemma~\ref{lem:energyestimate},
	\begin{equation*}
		(c)\leq C (\theta^{\alpha/2}+\Delta t^{\alpha/2}+\Delta t^{\frac{\alpha}{2}-\frac{\alpha}{q}} ).
	\end{equation*}
Now we can choose $\delta$, $\Delta t_0$ and $h_0$ small enough such that for all $\theta\leq \delta$, $\Delta t \leq \Delta t_0$ and $h\leq h_0$ (satisfying $h^{2(k+s)}\leq C\Delta t$)
	\begin{equation*}
			C_\alpha\E\left(\sum_{j=1}^L(u_h(\tau_h+\theta)-u_h(\tau_h),\phi_j)^2\right)^{\alpha/2}\leq \frac{\epsilon}{2}.
	\end{equation*}
	Now together with~\eqref{eq:halftheestimate} and the specific choice of $\zeta$, the estimate follows.
\end{proof}
\begin{remark}\label{rem:fortin}
	The condition $h^{2(k+s)}\leq C \Delta t$ and Assumption~\ref{ass:stokes} are needed to bound the pressure term (A.III) because $\Ih w$ is not necessarily discretely divergence free even if $w$ is divergence free. For a similar reason, the condition $h^{2(k+s)} = o(\Delta t)$ appears later in Lemma~\ref{lem:limits}. If $(\Uh,\Ph)$ admits a Fortin operator $\Pi_h: H^1_0(\dom)\to \Uh$ that is a projection onto $\Uh$, local and locally $W^{1,\infty}$-stable, these assumptions are not needed. Then, replacing $\Ih$ by $\Pi_h$ in~\eqref{eq:decomp1} and following, and using that $(\Div \Pi_h w,q) = (\Div w,q)=0$ for all $q\in \Ph$, term (A.III) vanishes. The remaining terms use $\norm{\Pi_h w}_{W^{1,\infty}}\leq C \norm{w}_{W^{1,\infty}}$ and $\norm{w-\Pi_h w}_{W^{1,\infty}} + \norm{w-\Pi_h w}_{H^1(\dom)} \to 0$ which follow from the locality and a Bramble--Hilbert argument on shape-regular meshes. Such an operator is known to exist for the MINI element~\cite{Boffi2013} and the lowest-order Taylor--Hood elements in any dimension~\cite{Diening2022}. A related construction is used in the analysis of a finite element incremental projection method in~\cite{Eymard2024}. For general inf-sup stable pairs, such a Fortin operator is not available in general, therefore we chose to use Assumption~\ref{ass:stokes} and the inverse CFL-type condition $h^{2(k+s)}\leq C \Delta t$ (and $h^{2(k+s)}= o(\Delta t)$) instead.
	\end{remark}
Next, we will adapt the tightness condition from~\cite[Corollary 1]{Motyl2013} to show that the laws of $u_h$ are tight on $\mathbb{D}([0,T];H^{-1}(\dom))$.
\begin{proposition}[Tightness criterion]
	\label{prop:tightness}
	Let $U$ and $H$ be Hilbert spaces satisfying $U\stackrel{c}{\hookrightarrow}  H \simeq H^* {\hookrightarrow}  U^*$ and assume $U$ is dense in $H$,
	and let $(X_n)_{n\in \N}$ be a sequence of c\`adl\`ag and $\mathcal{G}_t$-adapted $U^*$-valued processes such that
	\begin{enumerate}[(a)]
		\item $\exists C_1$ such that
		\begin{equation}
			\label{eq:supcond}
			\sup_{n\in \N} \E\sup_{s\in [0,T]}\norm{X_n(s)}_H \leq C_1,
		\end{equation}
		\item $(X_n)_{n\in \N}$ satisfies the Aldous condition,~\eqref{eq:aldous} in $U^*$.
	\end{enumerate}
	Let $\mu_n$ be the law of $X_n$ on $\mathbb{D}([0,T];U^*)$. 
	Then for all $\epsilon>0$, there exists a compact set $K_\epsilon$ of $\mathbb{D}([0,T];U^*)$  such that
	\begin{equation*}
		\mu_n(K_\epsilon)\geq 1-\epsilon,
	\end{equation*}
	for all $n\in \N$. 
	
\end{proposition}
We recall the relevant results from~\cite{Motyl2013,Joffe1986} and show how they combine to yield the proof of this result. We first recall the definition of modulus of continuity and modulus for c\`adl\`ag functions: 
\begin{definition}[Modulus of continuity]
	\label{def:moc}
	Given a function $u:[0,T]\to E$, we denote for an interval $I\subset [0,T]$,
	\begin{equation}\label{eq:mocI}
		w_E(u;I):= \sup_{s,t\in I}\norm{u(s)-u(t)}_E,
	\end{equation}
	and then define the \emph{modulus of continuity} for $\theta>0$,
	\begin{equation}
		\label{eq:moc}
		w_{[0,T],E}(u,\theta)= \sup_{0\leq t\leq t+\theta\leq T}w_E(u;[t,t+\theta]).
	\end{equation}
\end{definition}
A function $u$ is continuous on $[0,T]$ if and only if 
\begin{equation*}
	\lim_{\theta\to0} w_{[0,T],E}(u,\theta)=0.
\end{equation*}
However, functions in $\mathbb{D}([0,T];E)$ are not a priori continuous, therefore we need a different measure of their regularity, which is the following:
\begin{definition}[Modulus for c\`adl\`ag functions]
	\label{def:moccadlag}
	Given a function $u:[0,T]\to E$, 
	we denote for $\theta>0$,
	\begin{equation}
		\label{eq:moccadlag}
		w_{[0,T],E}'(u,\theta)= \inf\left\{\max_{i\leq r} w_E(u;[t_{i-1},t_i)): 0 =t_0<t_1<\dots< t_r = T,\, \inf(t_i-t_{i-1})\geq \theta\right\},
	\end{equation}
	where $w_E(u;I)$ was defined in~\eqref{eq:mocI}. When it is clear from the context, we omit the subscript $E$ and just write $w_{[0,T]}'(u,\theta)$.
\end{definition}
With this definition at hand, we recall an alternative to the Aldous condition:
\begin{definition}[Definition 4 in~\cite{Motyl2013}]
	\label{def:Ttilde}
	
	We say a sequence $\{X_n\}$ of $E$-valued c\`adl\`ag processes satisfies condition $\mathbf{[\tilde{T}]}$, if
	\begin{equation}
		\label{eq:conditionT}
		\begin{split}
				\forall \epsilon>0\quad \forall \eta>0\quad \exists\, \delta>0:\qquad \sup_{n\in \N} \mathbb{P}[w'_{[0,T]}(X_n,\delta)>\eta]\leq \epsilon.
		\end{split}
	\end{equation}
\end{definition}
Theorem 2.2.2 in~\cite{Joffe1986} yields that~\eqref{eq:aldous} implies~\eqref{eq:conditionT} with $\sup_{n\in\mathbb N}$ replaced by $\sup_{n\ge n_0}$. Since $w'_{[0,T]}(X_n,\delta)\to 0$ as $\delta\to 0$ for every c\`adl\`ag path by~\cite[Lemma 3.6.2]{Ethier1986}, the dominated convergence theorem implies that $\delta$ can be decreased so that $\mathbb{P}[w'_{[0,T]}(X_n,\delta)>\eta]\leq\epsilon$ also holds for the finitely many $n<n_0$. Thus~\eqref{eq:aldous} implies~\eqref{eq:conditionT}. Now, we can introduce the following tightness criterion:
\begin{lemma}[Theorem 3.7.2 in~\cite{Ethier1986}]
	\label{lem:compactinD}
	
	Let $\{X_n\}$ be a sequence of $E$-valued random variables with sample paths in $\mathbb{D}([0,T];E)$. Then $\{X_n\}$ is tight iff the following two conditions are satisfied:
	\begin{enumerate}[(a)]
		\item For all $\epsilon>0$, and all $t\in (\Q\cap[0,T])\cup \{0,T\}$, there exists a compact set $K_{\epsilon,t}\subset E$ such that 
		\begin{equation*}
			\inf_{n\in \N} \mathbb{P}[X_n(t)\in K_{\epsilon,t}]\geq 1-\epsilon,
		\end{equation*}
		\item Condition \eqref{eq:conditionT} is satisfied.
	\end{enumerate}
	
\end{lemma}
Now we are ready to prove Proposition~\ref{prop:tightness}. 
\begin{proof}[Proof of Proposition~\ref{prop:tightness}]
Let $\epsilon>0$ be arbitrary. Markov's inequality combined with~\eqref{eq:supcond} yields for any $r>0$,
	\begin{equation*}
		\mathbb{P}[\sup_{t\in [0,T]}\norm{X_n(t)}_H>r]\leq \frac{\E\sup_{t\in [0,T]}\norm{X_n(t)}_H}{r}\leq \frac{C_1}{r}.
	\end{equation*}
	We choose $R_1>0$ so large that $C_1/R_1 <\epsilon$. 
		By Theorem 2.2.2 in~\cite{Joffe1986}, we have that condition~\eqref{eq:conditionT} is satisfied. Since $U$ compactly embeds into $H$, we have that $H\equiv H^*$ compactly embeds into $U^*$ by Schauder's theorem. 
		Compact operators on a Hilbert space map the closed ball onto a compact set, therefore, the closed ball $\overline{B_{R_1}^{H}(0)}$ is a compact set in $U^*$ and satisfies $\mathbb{P}[X_n(t)\in B_{R_1}^{H}(0)]\geq 1-\epsilon$ by construction. This implies (a) in Lemma~\ref{lem:compactinD} is satisfied. Therefore the assumptions of  Lemma~\ref{lem:compactinD} are satisfied and we conclude that $X_n$ is tight in $\mathbb{D}([0,T];U^*)$.
\end{proof}
Now we can use the Skorokhod theorem to extract an almost surely convergent subsequence.
We define $r_h:= \tu_h - u_h$. Then we have the following:

\begin{proposition}[Skorokhod representation]
	\label{prop:skorokhod}
	Assume $r_0\geq 2$ and the assumptions of Lemma~\ref{lem:aldous}.
	There exist a subsequence of $(h,\Delta t)$, not relabeled, a probability space $(\bar\Omega,\bar{\mathcal G},\bar{\mathbb P})$, and random variables $\bar u_h$, $\bar r_h$, $\bar u$ on this probability space such that:
	\begin{enumerate}[(i)]
		\item $\mathrm{law}(\bar u_h,\bar r_h)=\mathrm{law}(u_h,r_h)=\PP\circ (u_h,r_h)^{-1}$ on $\mathbb{D}([0,T];H^{-1}(\dom))\times L^2([0,T]\times \dom)$ for every $h$,
		\item
$\bar u_h$ and $\bar{r}_h$ are constant on each $[t^m,t^{m+1})$ with values in $\Uh$ almost surely, $\bar{p}_h$ defined from $(\bar{u}_h,\bar{r}_h)$ via~\eqref{eq:projectionfully} is constant on $[t^m,t^{m+1})$ with values in $\Ph$ almost surely. Defining $\hu_h:=\bar{u}_h + \bar{r}_h$, the  triple $(\bar u_h,\hu_h,\bar p_h)$ has the law of $(u_h,\tu_h,p_h)$.
		\item The estimates of Lemma~\ref{lem:energyestimate} hold for $(\bar u_h,\hu_h)$ with the same constants and $\E$ replaced by $\bar{\E}$ (the expected value on $(\bar\Omega,\bar{\mathcal G},\bar{\mathbb P})$);
		
		\item $\bar{u}_h\to \bar{u}$ in $\mathbb{D}([0,T];H^{-1}(\dom))$ a.s. and $\bar{r}_h\to 0$ in $L^2([0,T]\times \dom)$ a.s.
		\item $\bar u\in L^q(\bar\Omega;L^2(0,T;H^1_0(\dom)))\cap L^q(\bar\Omega;L^\infty(0,T;L^2(\dom)))$ for $1\leq q \leq 2^{r_0}$.
		\item $\bar{\E}\|\bar u_h-\bar u\|_{L^2([0,T]\times \dom)}^q\to0$ and $\bar{\E}\|\hu_h-\bar u\|_{L^2([0,T]\times \dom)}^q\to0$ for $1\leq q < 2^{r_0}$, and along a further subsequence both converge strongly in $L^2([0,T]\times\dom)$ a.s.
	\end{enumerate}
\end{proposition}
\begin{proof}
	We apply Proposition~\ref{prop:tightness} to $\{X_n\}_{n} =\{u_h\}_{h>0}$ with $H=L^2(\dom)$ and $U = H^1_0(\dom)$. Condition (a) in Proposition~\ref{prop:tightness} follows from Lemma~\ref{lem:energyestimate}, (ii) and Condition (b) follows from Lemma~\ref{lem:aldous}.
	Thus by Prokhorov's theorem, $\mu_h:= \PP\circ (u_h)^{-1}\to \mu$ along a subsequence on $ \mathbb{D}([0,T];H^{-1}(\dom))$, where $\mu$ is a probability measure on $\mathbb{D}([0,T];H^{-1}(\dom))$. 
	Recall that $r_h:= \tu_h - u_h$. By Lemma~\ref{lem:energyestimate} (i), we have $\E\norm{r_h}^2_{L^2([0,T]\times\dom)}\leq C \Delta t$, therefore $r_h\to 0$ as $h,\Delta t\to 0$ in probability on $L^2([0,T]\times\dom)$.  Thus $\PP\circ(u_h,r_h)^{-1}\to \mu\otimes\delta_0$  on $ \mathbb{D}([0,T];H^{-1}(\dom))\times L^2([0,T]\times\dom)$ by~\cite[Theorem 2.7]{vanderVaart1998}  along a subsequence. The product $\mathbb{D}([0,T];H^{-1}(\dom))\times L^2([0,T]\times\dom)$ is a Polish space, therefore, we can use the Skorokhod representation theorem~\cite[Chapter I, Theorem 2.7]{Ikeda2014} or~\cite[Theorem I.6.7]{Billingsley1999} to conclude the existence of a probability space $(\bar{\Omega},\bar{\mathcal{G}},\bar{\mathbb{P}})$
	and random variables $\{(\bar{u}_h,\bar{r}_h)\}_{h>0}$ such that $(\bar{u}_h,\bar{r}_h)$ have the same law $\PP\circ(u_h,r_h)^{-1}$ as $(u_h,r_h)$ for all $h>0$, and $(\bar{u}_h(\omega),\bar{r}_h(\omega))\to (\bar{u}(\omega),\bar{r}(\omega))$ on $\mathbb{D}([0,T];H^{-1}(\dom))\times L^2([0,T]\times\dom)$ for almost every $\omega\in \bar{\Omega}$ and the law of $(\bar{u},\bar{r})$ is $\mu\otimes\delta_0$.  Thus in particular $\bar{r}=0$ a.s. This proves points (i) and (iv). 
	 Now define $\hu_h =\bar{u}_h +\bar{r}_h$.  
	Define $\bar{p}_h\in \Ph$ from $\bar{r}_h$ through~\eqref{eq:projection1}, which is unique by Lemma~\ref{lem:solvability}. 
	The set 
	\begin{equation*}
		A_{h,T}:= \{v\in \mathbb{D}([0,T];H^{-1}(\dom))\,:\, v \text{ is constant on each }\, [t^m,t^{m+1}) \, \text{ with values in }\, \Uh\}
	\end{equation*}
	is closed in $\mathbb{D}([0,T];H^{-1}(\dom))$ and therefore Borel-measurable, 
	therefore, since $\mathbb P(u_h\in A_{h,T})=1$ and (i), $\bar{u}_h$ is in $A_{h,T}$ almost surely. Similarly, $\bar{r}_h$ is in the space
		\begin{equation*}
		R_{h,T}:= \{v\in L^2([0,T]\times\dom)\,:\, v \text{ is constant on each }\, [t^m,t^{m+1}) \, \text{ with values in }\, \Uh\},
	\end{equation*}
	a.s.	because it is closed in $L^2([0,T]\times\dom)$.
	The mapping $\Phi:\Uh\times \Uh\to \Uh\times\Uh\times\Ph$ defined via~\eqref{eq:projection1}, $\Phi(u_h,r_h) = (u_h,u_h+r_h,p_h)$ on that finite dimensional space is continuous, as it is linear and bounded (by the LBB condition~\eqref{eq:LBB}). Therefore, since $(\bar{u}_h,\hu_h,\bar{p}_h )=\Phi(\bar{u}_h,\bar{r}_h)$, $(\bar{u}_h,\hu_h,\bar{p}_h )$ has the same law as $({u}_h,\tu_h,{p}_h )$ and $\hu_h$ and $\bar{p}_h $ are almost surely piecewise constant in time and taking values in $\Uh$ and $\Ph$ respectively. This proves (ii). 
	Then all the quantities in Lemma~\ref{lem:energyestimate} are continuous functions of $\bar{u}_h$, $\bar{r}_h$ and $\bar{p}_h$ and so since they have the same laws as $(u_h,r_h,p_h)$, the same estimates hold for $(\bar{u}_h,\bar{r}_h,\bar{p}_h)$ with $\E$ replaced by $\bar{\E}$. This implies (iii).
	Since $\bar{u}_h\to \bar{u}$ in $\mathbb{D}([0,T];H^{-1}(\dom))$ a.s., we have $\bar{u}_h(t)\to \bar{u}(t)$ in $H^{-1}(\dom)$ at every continuity point of $\bar{u}$ a.s. by~\cite[Prop. VI.2.4]{Jacod2013}. By the Skorokhod convergence, 
	\begin{equation*}
		\sup_{t\in [0,T]}\norm{\bar u_h(t)}_{H^{-1}(\dom)}=\sup_{t\in [0,T]}\norm{\bar u_h(\lambda_h(t))}_{H^{-1}(\dom)}\leq\sup_{t\in [0,T]}\norm{\bar u(t)}_{H^{-1}(\dom)}+\sup_{t\in [0,T]}\norm{\bar u_h\circ\lambda_h-\bar u}_{H^{-1}(\dom)},
	\end{equation*}
	where $\lambda_h$ is the bijection of $[0,T]$ from the Skorokhod metric,
	 therefore we also have $\sup_{h>0}\sup_{t\in [0,T]} \norm{\bar{u}_h}_{H^{-1}(\dom)}<\infty$ pathwise, and therefore by the dominated convergence theorem, we obtain $\bar{u}_h\to \bar{u}$ in $L^2(0,T;H^{-1}(\dom))$ a.s. Then also $\hu_h\to \bar{u}$ in $L^2(0,T;H^{-1}(\dom))$ a.s. since $\bar{r}_h\to 0$ in $L^2([0,T]\times\dom)\subset L^2(0,T;H^{-1}(\dom))$, so $\hu_h$ and $\bar{u}_h$ have the same limit a.s.
	 By weak lower semi-continuity of the norm and the uniqueness of the weak limit, we have
	 \begin{equation*}
	 	\bar{\E}\norm{\bar{u}}^q_{L^2(0,T;H^1_0(\dom))}\leq \liminf \bar{\E}\norm{\hu_h}^q_{L^2(0,T;H^1_0(\dom))}\leq C,
	 \end{equation*}
	 for $1\leq q \leq 2^{r_0}$ using the a priori estimates from Lemma~\ref{lem:energyestimate} (v).
	 To prove the second part of (v), we first fix $\omega$ with $\bar{u}_h\to \bar{u}$ in $L^2(0,T;H^{-1}(\dom))$.
	 By weak lower semi-continuity and the uniqueness of the weak* limit, we have
	 \begin{equation*}
	 	\norm{\bar{u}}_{L^\infty(0,T;L^2(\dom))}\leq \liminf_{h>0}\norm{\bar{u}_h}_{L^\infty(0,T;L^2(\dom))}.
	 \end{equation*}
	 We raise this to the power $q\in [1,2^{r_0}]$ and then use Fatou's lemma to obtain that
	 \begin{equation*}
	 	\bar{\E}\norm{\bar{u}}_{L^\infty(0,T;L^2(\dom))}^q\leq \liminf_{h>0}\bar{\E}\norm{\bar{u}_h}_{L^\infty(0,T;L^2(\dom))}^q<\infty,
	 \end{equation*}
	 where we used Lemma~\ref{lem:energyestimate} (iv) to bound the right-hand side.
	 This implies (v).
	  It remains to prove (vi). 
	Since 
	\begin{equation*}
		\norm{\hu_h-\bar{u}}_{L^2(0,T;H^{-1}(\dom))}^q\leq C\norm{\hu_h-\bar{u}}_{L^2([0,T]\times\dom)}^q\leq C_q(\norm{\hu_h}_{L^2([0,T]\times\dom)}^q+\norm{\bar{u}}_{L^2([0,T]\times\dom)}^q),
	\end{equation*} 
	for some $q\in (2,2^{r_0}]$, since $r_0\geq 2$,	it follows that
	\begin{equation*}
		\bar{\E}\norm{\hu_h-\bar{u}}_{L^2(0,T;H^{-1}(\dom))}^q \leq C,
	\end{equation*}
	and therefore by Vitali's convergence theorem
	\begin{equation*}
			\bar{\E}\norm{\hu_h-\bar{u}}_{L^2(0,T;H^{-1}(\dom))}^2\stackrel{h,\Delta t\to 0}{\longrightarrow} 0.
	\end{equation*}
	Next, we use Lions' lemma~\cite[Chapter III, Lemma 2.1]{Temam1977} which yields that for any given $\epsilon>0$ there exists $C_\epsilon>0$ such that for all $v\in L^2(0,T;H^1_0(\dom))$, one has
	\begin{equation*}
		\norm{v}_{L^2([0,T]\times\dom)}\leq \epsilon \norm{v}_{L^2(0,T;H^1_0(\dom))}+ C_\epsilon \norm{v}_{L^2(0,T;H^{-1}(\dom))}.
	\end{equation*}
	Applying this in expectation to $\hu_h-\bar{u}$, we have
	for any given $\epsilon>0$,
	\begin{align*}
		\bar{\E}\norm{\hu_h - \bar{u}}^2_{L^2([0,T]\times\dom)}&\leq 2\epsilon^2 \bar{\E}\norm{\hu_h - \bar{u}}_{L^2(0,T;H^1_0(\dom))}^2+ 2C_\epsilon^2 \bar{\E}\norm{\hu_h - \bar{u}}_{L^2(0,T;H^{-1}(\dom))}^2\\
		& \leq 2C \epsilon^2 +2C_\epsilon^2 \bar{\E}\norm{\hu_h - \bar{u}}_{L^2(0,T;H^{-1}(\dom))}^2\stackrel{h,\Delta t \to 0}{\longrightarrow} 2C\epsilon^2.
	\end{align*}
	Since $\epsilon>0$ was arbitrary, this implies the second part of (vi). The first part follows since
	\begin{equation*}
		\bar{\E}\norm{\bar{u}_h-\bar{u}}^2_{L^2([0,T]\times\dom)} \leq 2 \bar{\E}\norm{\bar{r}_h}^2_{L^2([0,T]\times\dom)}+ 2\bar{\E}\norm{\hu_h-\bar{u}}^2_{L^2([0,T]\times\dom)} \stackrel{h,\Delta t \to 0}{\longrightarrow} 0.
	\end{equation*}
	This implies convergence in probability. Using the energy estimate, Lemma~\ref{lem:energyestimate} (iv) and (v), and (v) of the current proposition, and the Vitali convergence theorem, we obtain also that
	\begin{equation*}
		\bar{\E}\norm{\bar{u}_h-\bar{u}}^q_{L^2([0,T]\times\dom)} , \bar{\E}\norm{\hu_h-\bar{u}}^q_{L^2([0,T]\times\dom)} \stackrel{h,\Delta t \to 0}{\longrightarrow} 0,
	\end{equation*}
	for $1\leq q <2^{r_0}$.
	By~\cite[Theorem 20.5]{Billingsley1995}, we obtain that a subsequence of $\bar{u}_h$ and $\hu_h$ converges a.s. in $L^2([0,T]\times\dom)$.
\end{proof}
Next, we will show that the limit $\bar{u}$ is almost surely divergence free.
\begin{lemma}
	\label{lem:divfree}
	
	The limit $\bar{u}$ satisfies $\bar{u}(t)\in L^2_{\Div}(\dom)$ for a.e. $t\in [0,T]$ almost surely. 
\end{lemma}
\begin{proof}
	Let $q\in H^2(\dom)$, $g\in L^2([0,T])$. Then $\Grad q g\in L^2([0,T]\times\dom)$. By Proposition~\ref{prop:skorokhod}, (vi), $\bar{u}_h\to \bar{u}$ in $L^2([0,T]\times\dom)$ on a set of full measure. So we fix $\omega\in \bar{\Omega}$ for which $\bar{u}_h(\omega)\to \bar{u}(\omega)$ and for which~\eqref{eq:projection2} holds.  Then
	\begin{equation*}
		\int_0^T\int_{\dom} \bar{u} \cdot \Grad q g dx dt  = \lim_{h\to 0}	\int_0^T\int_{\dom} \bar{u}_h \cdot \Grad q g dx dt, 
	\end{equation*}
	by the strong and therefore also weak convergence of $\bar{u}_h$ in $L^2([0,T]\times\dom)$. Using that  by~\eqref{eq:projection2}, $(\Div \bar{u}_h,q_h)=0$ for all $q_h\in \Ph$, we obtain that $(\Div \bar{u}_h,\Ih^p q)=0$. Thus
	\begin{align*}
	\left|	\int_0^T \int_{\dom} \bar{u}_h \cdot \Grad q g dx dt \right| & = \left|	\int_0^T \int_{\dom} \bar{u}_h \cdot (\Grad q -\Grad \Ih^p q) g dx dt \right| \\
	& \leq C \norm{\bar{u}_h}_{L^2([0,T]\times\dom)}\norm{g}_{L^2([0,T])}\norm{\Grad q - \Grad \Ih^p q}_{L^2(\dom)}\\
	& \leq C h \norm{\bar{u}_h}_{L^2([0,T]\times\dom)}\norm{g}_{L^2([0,T])}\norm{q}_{H^2(\dom)},
	\end{align*}
	where we have used Assumption~\ref{ass:approxorder}.
	Since $\bar{u}_h$ is bounded in $L^2([0,T]\times\dom)$ along the convergent subsequence, the right-hand side goes to zero as $h\to 0$. Thus
	\begin{equation*}
			\int_0^T\int_{\dom} \bar{u} \cdot \Grad q g dx dt = 0.
	\end{equation*}
	The choice of $\omega\in \bar{\Omega}$ was independent of $g$ and $q$, therefore $(\bar{u}(t),\Grad q)=0$ for almost any $t\in [0,T]$ and any $q\in H^2(\dom)$. 
	Now take a countable subset $\{q_j\}\subset H^2(\dom)$ that is dense in $H^1(\dom)$. Since $(\bar{u}(t),\Grad q_j)=0$ for all $j$, for $t$ outside the union of the countably many exceptional null sets, and $\{q_j\}$ is dense in $H^1(\dom)$, we have that $\bar{u}(t)$ is orthogonal to $\{\Grad q\, : \, q\in H^1(\dom)\}$, therefore $\bar{u}(t)\in L^2_{\Div}(\dom)$ for a.e. $t\in [0,T]$ by the Helmholtz decomposition and because $\bar{u}\in L^2(0,T;H^1_0(\dom))$ a.s. This can be repeated for every $\omega\in \bar{\Omega}$ outside a $\bar{\PP}$-null set where $\bar{u}_h$ is not convergent or not satisfying~\eqref{eq:projection2}.
\end{proof}

Next, we show that the limit $\bar{u}$ is almost surely continuous in time with values in $U^*=H^{-1}(\dom)$. To do so, we use Theorem 3.10.2(a) from~\cite{Ethier1986} which we recall here for convenience:

\begin{theorem}[cf. Theorem 10.2(a) in Chapter 3 of~\cite{Ethier1986}]
	\label{thm:ethiercontinuity}
	
	Let $(\mathbb{S},d)$ be a metric space and let
	$X_n$, $n=1,2,\dots$ and $X$ be processes with sample paths in $\mathbb{D}([0,T];\mathbb{S})$ and suppose that $X_n$ converges to $X$ in law. For $u\in\mathbb D([0,T];\mathbb{S})$, let
	\begin{equation*}
	J_T(u):=\int_0^T e^{-s}\min\{J(u,s),1\}\,ds+	e^{-T}	\min\{J(u,T),1\},
\end{equation*}
	where
\begin{equation*}
		J(u,s)	:=	\sup_{0<t\le s}d(u(t),u(t-)).
\end{equation*}
	Then $X$ has almost surely continuous sample paths on
	$[0,T]$ if and only if $J_T(X_n)$ converges to zero in law.
\end{theorem}
This follows from \cite[Theorem 3.10.2(a)]{Ethier1986}, by extending each path constantly from $T$ to $[T,\infty)$.
We want to apply this result to the sequence $\{\bar{u}_h\}_{h>0}$ with $\mathbb{S} = H^{-1}(\dom)$, $d(u,v)=\norm{u-v}_{H^{-1}}$, thus we need to show that $J_T(\bar{u}_h)$ converges in law to zero. This follows in a similar way as the verification of Aldous condition in Lemma~\ref{lem:aldous}:
\begin{lemma}
	\label{lem:uhcontinuity}
	$\{\bar{u}_h\}_{h>0}$ satisfies the conditions of Theorem~\ref{thm:ethiercontinuity} for $\mathbb{S}=H^{-1}(\dom)$, $d(u,v)=\norm{u-v}_{H^{-1}}$, thus the limit $\bar{u}$ is almost surely in $C([0,T];H^{-1}(\dom))$. In fact, it is also almost surely in $C([0,T];L^2(\dom)_w)$.
\end{lemma}
\begin{proof}
	We want to use Theorem~\ref{thm:ethiercontinuity}, thus we need to show that those conditions are satisfied.
	We already know that $\{\bar{u}_h\}_{h>0}$ is tight on $\mathbb{D}([0,T];H^{-1}(\dom))$ and therefore converges a.s. to $\bar{u}$ on $\mathbb{D}([0,T];H^{-1}(\dom))$ by  Proposition~\ref{prop:skorokhod} (iv). Thus if we can show that $\bar{\E}J_T(\bar{u}_h)\to 0$ as $h,\Delta t \to 0$, the result follows from Theorem~\ref{thm:ethiercontinuity}.
So we bound $\bar{\mathbb{E}}J_T(\bar{u}_h)$: We have using that $\bar{u}_h\sim u_h$, 
	\begin{align}
	\bar{\mathbb{E}}|J_T(\bar{u}_h)|
		& \leq  C_T\bar{\mathbb{E}}\sup_{0\leq t\leq T}\norm{\bar{u}_h(t)-\bar{u}_h(t-)}_{H^{-1}}  \notag\\
		&\leq C_T{\mathbb{E}}\sup_{0\leq m\leq N-1}\norm{u_h^{m+1}-u_h^m}_{H^{-1}} .\label{eq:jumpcond2}
	\end{align}
	So if we can show that
	\begin{equation}\label{eq:jumpcond}
		\E \left[\sup_{0\leq m\leq N-1}\norm{u^{m+1}_h-u_h^m}_{H^{-1}}\right]\stackrel{h,\Delta t\to 0}{\longrightarrow} 0,
	\end{equation}
the result of the lemma follows.
	To prove~\eqref{eq:jumpcond}, we note that using Lemma~\ref{lem:Rellich3}, it reduces to showing
	\begin{equation}\label{eq:easier}
		\E \max_{0\leq m\leq N-1}\left|(u^{m+1}_h-u_h^m,\phi)\right| \stackrel{h,\Delta t\to 0}{\longrightarrow} 0,
	\end{equation}
for $\phi\in C^\infty_{c,\Div}(\dom)\oplus \Grad C^\infty(\overline{\dom})$. As in Lemma~\ref{lem:aldous}, we decompose $\phi=w+\Grad g$ where $w\in C^\infty_{c,\Div}(\dom)$ and $g\in C^\infty(\overline{\dom})$. We plug in the scheme~\eqref{eq:predictionfully} -- \eqref{eq:projectionfully}:
\begin{align}\label{eq:Jestimate}
(u^{m+1}_h-u_h^m,\phi) & = (u^{m+1}_h-u_h^m,\Ih w) + (u^{m+1}_h-u_h^m, w - \Ih w)\\
	&\qquad + (u^{m+1}_h-u_h^m,\Grad g)\notag\\
	& = -\underbrace{\Delta t b(\tu_h^m,\tu_h^{m+1},\Ih  w)}_{\text{(A.I)}} - \underbrace{\nu\Delta t (\Grad \tu_h^{m+1},\Grad \Ih w)}_{\text{(A.II)}}\notag\\
	& +\underbrace{\Delta t ( p_h^{m+1},\Div \Ih w)}_{\text{(A.III)}} + \underbrace{\Delta t \langle f^m,\Ih w\rangle }_{\text{(A.IV)}}
	+ \underbrace{(G(u_h^m)\Delta W_m,\Ih w)}_{\text{(A.V)}} \notag\\
	&+ \underbrace{(u^{m+1}_h-u_h^m, \Grad g )}_{\text{(B)}}+ \underbrace{(u^{m+1}_h-u_h^m,w-\Ih w)}_{\text{(C)}}.\notag
\end{align}	
We estimate the terms on the right-hand side independently of $m$. All terms except for (A.V) are similar to Lemma~\ref{lem:aldous}. Using Assumption~\ref{ass:approxorder},
\begin{equation*}
	|\text{(A.I)}| \leq \sqrt{\Delta t} C\norm{\tu_h}_{L^\infty(0,T;L^2(\dom))}\norm{\tu_h}_{L^2(0,T;H^1_0(\dom))}\norm{w}_{W^{1,\infty}(\dom)},
\end{equation*}
\begin{equation*}
	|\text{(A.II)}| \leq \sqrt{\Delta t}\nu C \norm{\Grad \tu_h}_{L^2([0,T]\times\dom)}\norm{w}_{H^2(\dom)}.
\end{equation*}
Using a similar estimate as in~\eqref{eq:LBBpressureestimate}, 
we have for (A.III)
\begin{align*}
	|\text{(A.III)}| & \leq C   \norm{u^{m+1}_h - \tu_h^{m+1}}_{L^2}\norm{\Ih w-w_h}_{L^2}\\
	&\leq C \norm{u^{m+1}_h - \tu_h^{m+1}}_{L^2}\left(h^{k+1}\norm{w}_{H^{k+1}}+h^{s+k}\norm{w}_{H^{k+1}}\right)\\
	& \leq C  {h^{k+s}}  \left(\norm{u_h}_{L^\infty(0,T;L^2(\dom))}+ \norm{\tu_h}_{L^\infty(0,T;L^2(\dom))}\right)\norm{w}_{H^{k+1}},
\end{align*}
where $w_h$ denotes the Stokes projection of $w$.
For (A.IV), we have
\begin{equation*}
	|\text{(A.IV)}| \leq \sqrt{\Delta t}C\norm{f}_{L^2(0,T;H^{-1}(\dom))}\norm{w}_{H^2},
\end{equation*}
and for (A.V), we have
\begin{equation*}
	|\text{(A.V)}| \leq C\max_{0\leq m\leq N-1}\norm{G(u_h^m)\Delta W_m}_{L^2(\dom)}\norm{w}_{H^2}.
\end{equation*}
Terms (B) and (C) are the same as in Lemma~\ref{lem:aldous} above and estimated as in~\eqref{eq:termB} and~\eqref{eq:termC}. Taking the maximum over $m$ and then expectations in~\eqref{eq:Jestimate}, and using these estimates, we obtain, using Lemma~\ref{lem:energyestimate},
\begin{align*}
	\E\max_{0\leq m\leq N-1}|(u^{m+1}_h-u_h^m,\phi) | & \leq 
	\sqrt{\Delta t} C\E \norm{\tu_h}_{L^\infty(0,T;L^2(\dom))}\norm{\tu_h}_{L^2(0,T;H^1_0(\dom))}\norm{w}_{W^{1,\infty}(\dom)} \\
	& +\sqrt{\Delta t}\nu C \E\norm{\Grad \tu_h}_{L^2([0,T]\times\dom)}\norm{w}_{H^2(\dom)}\\
	& +C  {h^{k+s}}  \left(\E\norm{u_h}_{L^\infty(0,T;L^2(\dom))}+ \E\norm{\tu_h}_{L^\infty(0,T;L^2(\dom))}\right)\norm{w}_{H^{k+1}}\\
	&+ \sqrt{\Delta t}C\norm{f}_{L^2(0,T;H^{-1}(\dom))}\norm{w}_{H^2}+ C\E\max_{0\leq m\leq N-1}\norm{G(u_h^m)\Delta W_m}_{L^2(\dom)}\norm{w}_{H^2} \\
	&+C\E \norm{u_h}_{L^\infty(0,T;L^2(\dom))}h^{\ell} \norm{g}_{H^{\ell+1}}+ C\E \norm{u_h}_{L^\infty(0,T;L^2(\dom))}h^{k+1} \norm{w}_{H^{k+1}}\\
	& \leq C (\sqrt{\Delta t} + h^{k+s} + h^\ell) + C\E\max_{0\leq m\leq N-1}\norm{G(u_h^m)\Delta W_m}_{L^2(\dom)}.
\end{align*}
So it remains to estimate $\E\max_{0\leq m\leq N-1}\norm{G(u_h^m)\Delta W_m}_{L^2(\dom)}$. For $1<r\leq 2^{r_0-1}$, we have
\begin{align*}
	\E\max_{0\leq m\leq N-1}\norm{G(u_h^m)\Delta W_m}_{L^2(\dom)}& \leq \left(\E\sum_{m=0}^{N-1} \norm{G(u_h^m)\Delta W_m}_{L^2(\dom)}^{2r}\right)^{\frac{1}{2r}}\\
	& \leq \left(\E\sum_{m=0}^{N-1} \norm{G(u_h^m)}^{2r}_{\mathcal{L}(\mathcal{K},L^2)}\norm{\Delta W_m}_{\mathcal{K}}^{2r}\right)^{\frac{1}{2r}}\\
	& \leq CK_1 \left(\E\sum_{m=0}^{N-1} \left(\norm{ u_h^m}^{2r}_{L^2(\dom)}+1\right)\norm{\Delta W_m}_{\mathcal{K}}^{2r}\right)^{\frac{1}{2r}}\\
	& \leq K_1C_r(\tr Q)^{\frac12} \Delta t^{\frac{r-1}{2r}}\left(\left(\E  \norm{ u_h }^{2r}_{L^{2r}(0,T;L^2(\dom))}\right)^{\frac{1}{2r}}+1\right),
\end{align*}
which goes to zero as $\Delta t\to 0$. Thus~\eqref{eq:easier} follows. Since we also have that $\bar{u}\in L^\infty(0,T;L^2(\dom))$ almost surely by Proposition~\ref{prop:skorokhod} (v), we can use~\cite[Lemma II.5.9]{Boyer2013} to conclude that $\bar{u}\in C([0,T];L^2(\dom)_w)$ almost surely.
\end{proof}
	Since $L^2_{\Div}(\dom)$ is a weakly closed subspace of	$L^2(\dom)$, the weak continuity of $\bar u$ and Lemma~\ref{lem:divfree} imply that $\bar u(t)\in L^2_{\Div}(\dom)$ for every $t\in[0,T]$, almost surely.
Next, we pass to the limit in all the terms in the scheme except for the noise term which will be handled later.
\begin{lemma}
	\label{lem:limits}
	Under Assumption~\ref{ass:stokes}, the condition $h^{2(k+s)} = o(\Delta t)$ where $k$ is the approximation order of the finite element space $\Uh$, and $r_0\geq 3$,
	we have along the subsequence of Proposition~\ref{prop:skorokhod} for any {$\varphi\in C^\infty_{c,\Div}(\dom)$} and for any $t\in [0,T]$, and any $1\leq \alpha <2^{r_0-1}$,
		\begin{align}
		\lim_{h,\Delta t\to 0} \bar{\E}\left|\int_0^t (\Grad\hu_h,\Grad \Ih \varphi) ds-\int_0^t (\Grad \bar{u},\Grad \varphi) ds \right|^\alpha&= 0,\label{eq:dissipationlimit}\\
		\lim_{h,\Delta t\to 0} \bar{\E}\left|\int_0^t (\bar{p}_h(s+\Delta t),\Div \Ih \varphi) ds \right|^{\alpha}& = 0,\label{eq:pressurelimit}\\
		\lim_{h,\Delta t\to 0} \bar{\E}\left|\int_0^t b(\hu_h(s-\Delta t),\hu_h,\Ih \varphi) ds +\int_0^t \int_{\dom}(\bar{u}\cdot\Grad )\varphi \cdot \bar{u} dx ds\right|^\alpha& =0,\label{eq:convectionlimit}\\
		\lim_{h,\Delta t\to 0} \bar{\E}\left|(\bar{u}_h(t),\Ih \varphi) - (\bar{u}(t),\varphi)\right|^\alpha& =0,\label{eq:fixedtimelimit}\\
			\lim_{h,\Delta t\to 0} (\bar{u}_h(0),\Ih \varphi) = (\bar{u}(0),\varphi)& =(u_0,\varphi),\quad \text{a.s.},\label{eq:initimelimit}\\
		\lim_{h,\Delta t\to 0} \bar{\E}\int_0^t \langle f_h,\Ih \varphi\rangle ds  =\lim_{h,\Delta t\to 0} \int_0^t \langle f_h,\Ih \varphi\rangle ds   & = \int_0^t \langle f,\varphi\rangle  ds,\label{eq:forcinglimit}
	\end{align}
	where $\bar{p}_h(t)\in \Ph$ is defined through
	\begin{equation*}
		\left(\frac{\bar{u}_h(t)-\hu_h(t-\Delta t)}{\Delta t}, v\right)= (\bar{p}_h,\Div v),\qquad (\Div \bar{u}_h(t),q)=0,
	\end{equation*}
	for all $(v,q)\in \Uh\times\Ph$ (cf.~\eqref{eq:projectionfully}). 
\end{lemma}
\begin{proof}
	We start by rewriting the dissipation term~\eqref{eq:dissipationlimit}:
	\begin{align*}
		\bar{\E}\left|\int_0^t \left((\Grad \hu_h,\Grad \Ih \varphi)-(\Grad \bar{u},\Grad \varphi) \right) ds\right|^\alpha&  \leq C_\alpha \underbrace{\bar{\E}\left|\int_0^t (\Grad \hu_h-\Grad \bar{u},\Grad\varphi) ds\right|^\alpha}_{\text{(a)}}\\
		&\quad + C_\alpha\underbrace{\bar{\E}\left|\int_0^t (\Grad\hu_h,\Grad (\Ih\varphi -\varphi))ds \right|^\alpha}_{\text{(b)}}.
	\end{align*}
	We show that both terms (a) and (b) converge to zero:
	For (a), we integrate by parts and then use that $\hu_h\to \bar{u}$ in $L^2([0,T]\times\dom)$ $\bar{\mathbb{P}}$-a.s. and the moment bounds from Lemma~\ref{lem:energyestimate} to conclude with Vitali's convergence theorem that (a)$\to 0$:
	\begin{equation*}
		\text{(a)} = \bar{\E}\left|\int_0^t (\hu_h-\bar{u},\Delta\varphi) ds\right|^\alpha\stackrel{h,\Delta t\to 0}{\longrightarrow} 0.
	\end{equation*}
	For term (b), we use Assumption~\ref{ass:approxorder} to obtain
	\begin{align*}
		\text{(b)}& \leq T^{\alpha/2}\bar{\E} \norm{\Grad \hu_h }_{L^2([0,T]\times\dom)}^\alpha \norm{\Grad \Ih\varphi-\Grad \varphi}^\alpha_{L^2}\\
		& \leq  C{\E} \norm{\Grad \tu_h}_{L^2([0,T]\times\dom)}^\alpha h^\alpha\norm{ \varphi}_{H^2}^\alpha\stackrel{h,\Delta t\to 0}{\longrightarrow} 0.
	\end{align*}
	Thus,~\eqref{eq:dissipationlimit} follows. We next consider the pressure term. Using that $\varphi$ is divergence free, we can rewrite using the Stokes projection $\varphi_h$ of $\varphi$,
	\begin{equation*}
		 \bar{\E}\left|\int_0^t (\bar{p}_h(s+\Delta t),\Div \Ih \varphi) ds\right|^\alpha  =  \bar{\E}\left|\int_0^t ( \bar{p}_h(s+\Delta t),\Div\Ih \varphi-\Div\varphi_h) ds\right|^\alpha.
	\end{equation*}
	Thus, we can estimate, using~\eqref{eq:rewrite}, Assumptions~\ref{ass:approxorder} and~\ref{ass:stokes}, Lemma~\ref{lem:stokesestimate}, and the energy estimate, Lemma~\ref{lem:energyestimate},
	\begin{align*}
		 \bar{\E}\left|\int_0^t (\bar{p}_h(s+\Delta t),\Div \Ih \varphi) ds\right|^\alpha  & \leq  C\bar{\E}\left|\int_0^t \frac{\norm{ \bar{u}_h(s+\Delta t)-\hu_h}_{L^2}}{\Delta t} ds \right|^\alpha \left(h^{(k+1)\alpha} \norm{\varphi}^\alpha_{H^{k+1}(\dom)}+h^{(k+s)\alpha} \norm{\varphi}_{H^{k+1}(\dom)}^\alpha\right)\\
&  = C {\E}\left|\int_0^t \frac{\norm{ {u}_h(s+\Delta t)-\tu_h}_{L^2}}{\Delta t} ds \right|^\alpha \left(h^{(k+1)\alpha} \norm{\varphi}_{H^{k+1}(\dom)}^\alpha+h^{(k+s)\alpha} \norm{\varphi}_{H^{k+1}(\dom)}^\alpha\right)\\
		 & \leq C\frac{h^{(k+s)\alpha}}{\Delta t^{\alpha/2}}\norm{\varphi}^\alpha_{H^{k+1}(\dom)}
	\end{align*}
	Under the condition that {$h^{2(k+s)} = o(\Delta t)$}, this term goes to zero and~\eqref{eq:pressurelimit} follows. Next, we consider the convection term, which we rewrite as follows using that $\bar{u}$ is divergence free almost surely (Lemma~\ref{lem:divfree}):
	\begin{align*}
		&\bar{\E} \left|\int_0^t \left(b(\hu_h(\cdot -\Delta t),\hu_h,\Ih \varphi)+\int_{\dom}(\bar{u}\cdot\Grad)\varphi \cdot \bar{u} dx\right) ds \right|^\alpha\\
		& \leq  C_\alpha\underbrace{\bar{\E}\left|\int_0^t\int_{\dom} (\bar{u}\cdot\Grad) (\hu_h - \bar{u} ) \cdot\varphi dx ds\right|^\alpha}_{\text{(a)}}\\
		& + {C_\alpha}  \underbrace{\bar{\E}\left|\int_0^t \int_{\dom} \left[((\hu_h(s-\Delta t)-\bar{u})\cdot\Grad) \hu_h\cdot \varphi -((\hu_h(s-\Delta t)-\bar{u})\cdot\Grad)\varphi \cdot\hu_h \right]dx ds\right|^\alpha}_{\text{(b)}}\\
		& +  C_\alpha \underbrace{\bar{\E}\left|\int_0^t \int_{\dom} \left[(\hu_h(s-\Delta t)\cdot\Grad) \hu_h\cdot (\Ih \varphi -\varphi) -(\hu_h(s-\Delta t)\cdot\Grad)(\Ih \varphi -\varphi)\cdot \hu_h \right]dx ds\right|^\alpha}_{\text{(c)}}
 	\end{align*}
	We can rewrite and estimate term (a) as follows, using that $\bar{u}$ is a.s. divergence free:
	\begin{align*}
		\text{(a)}& = \bar{\mathbb{E}}\left|\int_0^t\int_{\dom} (\bar{u}\cdot\Grad)\varphi\cdot (\hu_h-\bar{u}) dx ds\right|^\alpha\\
		& \leq \bar{\mathbb{E}}\norm{\bar{u}}_{L^2([0,T]\times\dom)}^\alpha\norm{\hu_h-\bar{u}}_{L^2([0,T]\times\dom)}^\alpha \norm{\Grad \varphi}_{L^\infty}^\alpha\\
		& \leq \left(\bar{\mathbb{E}}\norm{\bar{u}}_{L^2([0,T]\times\dom)}^{2\alpha}\right)^{1/2}\left(\bar{\mathbb{E}}\norm{\hu_h-\bar{u}}_{L^2([0,T]\times\dom)}^{2\alpha}\right)^{1/2} \norm{\Grad \varphi}_{L^\infty}^\alpha \stackrel{h,\Delta t\to 0}{\longrightarrow} 0,
	\end{align*}
	where we used the $\bar{\mathbb{P}}$-a.s. convergence of $\hu_h$ to $\bar{u}$ in $L^2([0,T]\times\dom)$ and the uniform integrability bounds from the discrete energy estimate, Lemma~\ref{lem:energyestimate}, to pass to the limit with Vitali's convergence theorem.
	For term (b), we note that thanks to the energy estimate, Lemma~\ref{lem:energyestimate}, and Proposition~\ref{prop:skorokhod}, we have
	\begin{align*}
		\bar{\E}\norm{\hu_h(\cdot-\Delta t)-\bar{u}}_{L^2}^2 &\leq  	2\bar{\E}\norm{\hu_h(\cdot-\Delta t)-\bar{u}_h}_{L^2}^2+2 \bar{\E}\norm{\bar{u}_h-\bar{u}}_{L^2}^2\\
		& = 2{\E}\norm{\tu_h(\cdot-\Delta t)-{u}_h}_{L^2}^2+2 \bar{\E}\norm{\bar{u}_h-\bar{u}}_{L^2}^2\\
		& \leq C \Delta t + 2\bar{\E}\norm{\bar{u}_h-\bar{u}}_{L^2}^2 \stackrel{h,\Delta t\to 0}{\longrightarrow} 0.
	\end{align*}
Therefore,
	\begin{align*}
		\text{(b)} & \leq C\left(\bar{\E}\norm{\hu_h(\cdot-\Delta t)-\bar{u}}_{L^2([0,T]\times\dom)}^{2\alpha}\right)^{1/2}\left(\bar{\E}\norm{\Grad\hu_h}_{L^2([0,T]\times\dom)}^{2\alpha}\right)^{1/2}\norm{\varphi}_{L^\infty}^\alpha\\ &+C\left(\bar{\E}\norm{\hu_h(\cdot-\Delta t)-\bar{u}}_{L^2([0,T]\times\dom)}^{2\alpha}\right)^{1/2}\left(\bar{\E}\norm{\hu_h}_{L^2([0,T]\times\dom)}^{2\alpha}\right)^{1/2}\norm{\Grad\varphi}_{L^\infty}^\alpha\\
		&= C\left(\bar{\E}\norm{\hu_h(\cdot-\Delta t)-\bar{u}}_{L^2([0,T]\times\dom)}^{2\alpha}\right)^{1/2}\left({\E}\norm{\Grad\tu_h}_{L^2([0,T]\times\dom)}^{2\alpha}\right)^{1/2}\norm{\varphi}_{L^\infty}^\alpha\\ &+C\left(\bar{\E}\norm{\hu_h(\cdot-\Delta t)-\bar{u}}_{L^2([0,T]\times\dom)}^{2\alpha}\right)^{1/2}\left({\E}\norm{\tu_h}_{L^2([0,T]\times\dom)}^{2\alpha}\right)^{1/2}\norm{\Grad\varphi}_{L^\infty}^\alpha\\		
		& \leq C \left(\bar{\E}\norm{\hu_h(\cdot-\Delta t)-\bar{u}}_{L^2([0,T]\times\dom)}^{2\alpha}\right)^{1/2}\stackrel{h,\Delta t\to 0}{\longrightarrow} 0.
	\end{align*}
	For term (c), we use Assumption~\ref{ass:approxorder},
	\begin{align*}
		\text{(c)}& \leq C\left(\bar{\E}\norm{\hu_h}_{L^2([0,T]\times\dom)}^{2\alpha}\right)^{1/2}\left(\bar{\E}\norm{\Grad \hu_h}_{L^2([0,T]\times\dom)}^{2\alpha}\right)^{1/2}\norm{\Ih\varphi -\varphi}_{L^\infty}^\alpha\\
		&\quad + C \bar{\E}\norm{\hu_h}_{L^2([0,T]\times\dom)}^{2\alpha}\norm{\Grad\Ih\varphi -\Grad \varphi}_{L^\infty}^\alpha\\
		& \leq C h^{k\alpha}\norm{\varphi}_{W^{k+1,\infty}(\dom)}^\alpha.
	\end{align*}
	Thus, also this term vanishes as $h,\Delta t\to 0$, and we obtain~\eqref{eq:convectionlimit}. We consider~\eqref{eq:fixedtimelimit}. We can bound it by
	\begin{equation}\label{eq:temp}
		\bar{\E}\left|(\bar{u}_h(t),\Ih \varphi) - (\bar{u}(t),\varphi)\right|^\alpha
		\leq C_\alpha \bar{\E}\left|(\bar{u}_h(t),\Ih \varphi-\varphi)\right|^\alpha + C_\alpha \bar{\E}\left|(\bar{u}_h(t)-\bar{u}(t), \varphi)\right|^\alpha
	\end{equation}
	The first term on the right-hand side converges to zero using the energy estimate, Lemma~\ref{lem:energyestimate}, and Assumption~\ref{ass:approxorder}. For the second term, we interpret $\bar{u}_h(t)$ as an element of $H^{-1}(\dom)$ and $\varphi$ as an element of $H^1_0(\dom)$ (since it is assumed to be smooth) and use the a.s. convergence of $\bar{u}_h$ in $\mathbb{D}([0,T];H^{-1}(\dom))$. Lemma 10.1 in Chapter 3 of~\cite{Ethier1986} then yields a.s. in $\omega\in \bar{\Omega}$,
	\begin{equation*}
		\limsup_{h,\Delta t\to 0}\sup_{0\leq t\leq T} \norm{\bar{u}_h(t)-\bar{u}(t)}_{H^{-1}(\dom)}\leq J(\bar{u},T),
	\end{equation*}
	where 
	\begin{equation*}
		J(\bar{u},T) = \sup_{0\leq t\leq T} \norm{\bar{u}(t)-\bar{u}(t-)}_{H^{-1}}.
	\end{equation*}
	Since $\bar{u}$ is continuous in time with values in $H^{-1}(\dom)$ almost surely by Lemma~\ref{lem:uhcontinuity}, we have $J(\bar{u},T)=0$ almost surely and thus $\limsup_{h,\Delta t\to 0}\sup_{0\leq t\leq T} \norm{\bar{u}_h(t)-\bar{u}(t)}_{H^{-1}(\dom)}=0$ which implies the second term in~\eqref{eq:temp} vanishes. Combined with the moment estimates, Lemma~\ref{lem:energyestimate} and Vitali's convergence theorem, this implies~\eqref{eq:fixedtimelimit}. The convergence at the initial time~\eqref{eq:initimelimit} follows similarly, by setting $t=0$. The second equality in~\eqref{eq:initimelimit} follows because $u_0$ is deterministic,  and we have
		\begin{equation*}
		|(u^0_h-\tu^0_h,\Ih\varphi)|=|(u^0_h-\tu^0_h,\Ih\varphi-\varphi_h)|\le C\|u^0_h-\tu^0_h\|_{L^2}\,h^{k+s}\|\varphi\|_{H^{k+1}}\le C\|u_0\|_{L^2}h^{k+s}\to0,
	\end{equation*}
	where $\varphi_h$ is the Stokes projection of $\varphi$ as in~\eqref{eq:stokesprojection}
	and 
	\begin{equation*}
		(\tu^0_h,\Ih\varphi)=(u_0,\Ih\varphi) \stackrel{h\to 0}{\longrightarrow} (u_0,\varphi).
	\end{equation*}
	Finally for~\eqref{eq:forcinglimit}, we have
	\begin{align*}
		\bar{\E}\int_0^t\langle f_h,\Ih\varphi\rangle ds & = \int_0^t\langle f_h,\Ih\varphi\rangle ds\\
		& =\int_0^t \langle f_h,\varphi\rangle ds + \int_0^t \langle f_h,\Ih\varphi-\varphi\rangle ds.
	\end{align*}
	The first term converges 
	\begin{equation*}
		\int_0^t \langle f_h,\varphi\rangle ds \stackrel{h,\Delta t\to 0}{\longrightarrow} \int_0^t \langle f,\varphi\rangle ds,
	\end{equation*}
	since time averages converge in $L^2([0,T])$. 
	We can bound the second term as follows:
	\begin{equation*}
		\left|\int_0^t \langle f_h,\Ih\varphi-\varphi\rangle ds\right|\leq C\norm{f}_{L^2(0,T;H^{-1}(\dom))}\norm{\Grad\varphi-\Grad\Ih\varphi}_{L^2}\leq C h\norm{\varphi}_{H^2}\stackrel{h,\Delta t\to 0}{\longrightarrow} 0.
	\end{equation*}
	This implies~\eqref{eq:forcinglimit} and concludes the proof of the lemma.
\end{proof}
Finally, we show that there exists a stochastic basis and a $Q$-Wiener process $\bar{W}$ such that $\bar{u}$ is a martingale solution of~\eqref{eq:SNSm}.
\begin{theorem}
	\label{lem:martingaleidentification}
	Under the assumptions of Proposition~\ref{prop:skorokhod} and Lemma~\ref{lem:limits} and $r_0\geq 3$,	there exists a stochastic basis $(\doublehat{\Omega},\doublehat{\mathcal{G}},(\doublehat{\mathcal{G}}_t)_{t\in [0,T]},\doublehat{\mathbb{P}})$ and a $Q$-Wiener process $\bar{W}$, such that $(\doublehat{\Omega},\doublehat{\mathcal{G}},(\doublehat{\mathcal{G}}_t)_{t\in [0,T]},\doublehat{\mathbb{P}},\bar{u},\bar{W})$ is a martingale solution of~\eqref{eq:SNSm}, with the limit $\bar{u}$ from Proposition~\ref{prop:skorokhod}.
\end{theorem}
\begin{proof}
	We have already identified a limit $\bar{u}$ via the Skorokhod representation theorem and passed to the limit in all the terms of the equation in Lemma~\ref{lem:limits} except for the noise term. To identify the noise term in the limit, we use the martingale representation theorem~\cite[Theorem 8.2]{DaPrato2014}. To do so, we first verify the martingale property at the discrete level and then pass to the limit using Lemma~\ref{lem:limits}.
So we define $\hM_h(t)\in \Uh$ via  
\begin{align}\label{eq:Mhdef}
	\begin{split}
	(\hM_h(t),v) &= (\bar{u}_h(t),v) - (\bar{u}_h(0),v)+\int_0^{\kappa_h(t)}b(\hu_h(s-\Delta t),\hu_h(s),v)ds +\nu \int_0^{\kappa_h(t)} (\Grad \hu_h(s),\Grad v) ds\\ 
	&\quad - \int_0^{\kappa_h(t)}(\bar{p}_h(s+\Delta t),\Div v) ds-\int_0^{\kappa_h(t)}\langle f_h,v\rangle ds
\end{split}
\end{align} 
for all $v\in \Uh$,
where 
\begin{equation*}
	\kappa_h(t) = t^m,\quad t\in [t^m,t^{m+1}),\quad m=0,\dots,N.
\end{equation*}
$\hM_h$ is a square integrable martingale with respect to the filtration generated by $\bar{\mathcal{G}}^h_{t^m}=\sigma( \bar{u}_h^0,\dots, \bar{u}_h^m,\allowbreak \bar{r}_h^0,\dots, \bar{r}_h^{m-1})$ (recall $\bar{r}_h^m =\hu_h^{m+1}-\bar{u}_h^m$) with predictable operator quadratic variation
\begin{equation*}
	\bar{R}_h(t) = \int_0^{\kappa_h(t)} \PUh G(\bar{u}_h(\tau))Q^{1/2}\left[\PUh G(\bar{u}_h(\tau))Q^{1/2}\right]^* d\tau,
\end{equation*}
since for any $\phi$ bounded and continuous on  $\mathbb{D}([0,s];H^{-1}(\dom))$, and $0\leq s\leq t\leq T$, we have, for $v\in \Uh$
\begin{align*}
	\bar{\E}\left[(\hM_h(t)-\hM_h(s),v) \phi(\bar{u}_h|_{[0,s]})\right] &= \E \left[\left(\int_{\kappa_h(s)}^{\kappa_h(t)}G(u_h)dW(\tau), v\right) \phi(u_h|_{[0,s]})\right] \\
	& =\E \left[\sum_{m=m_1}^{m_2-1}\left(G(u^m_h)\Delta W_m, v\right) \phi(u_h|_{[0,s]})\right] \\
	& =0,
\end{align*}
where $m_1$ and $m_2$ are such that $\kappa_h(t)=t^{m_2}$ and $\kappa_h(s)=t^{m_1}$ and we used that $u_h|_{[0,s]}$ is $\mathcal{G}_{t^{m_1}}$-measurable. 

The same calculation applies when $\phi(\bar u_h|_{[0,s]})$ is replaced by an arbitrary
bounded Borel function of the completed history $(\bar u_h^0,\ldots,\bar u_h^{m_1},\allowbreak 
\bar r_h^0,\ldots,\bar r_h^{m_1-1})$ where $\kappa_h(s)=t^{m_1}$. 
Indeed, as the corresponding $( u_h^0,\ldots, u_h^{m_1},\allowbreak 
r_h^0,\ldots,r_h^{m_1-1})$ are $\mathcal G_{t^{m_1}}$-measurable, the conditional mean-zero property of the subsequent Wiener increments still applies and the identity above results with arbitrary
bounded Borel function of the completed history $(\bar u_h^0,\ldots,\bar u_h^{m_1}, \allowbreak
\bar r_h^0,\ldots,\bar r_h^{m_1-1}).$ In the same way, the covariance calculation below can be extended to more general bounded Borel functions of $(\bar u_h^0,\ldots,\bar u_h^{m_1},\allowbreak 
\bar r_h^0,\ldots,\bar r_h^{m_1-1})$.
Note that we also have, using the BDG-inequality~\eqref{eq:BDG} and the growth assumption on $G$,~\eqref{eq:lingrowthG}, for $1\leq q\leq 2^{r_0}$
\begin{align}\label{eq:momentboundM}
	\begin{split}
		\bar{\E}\sup_{0\le t\le T}	\norm{\hM_h(t)}_{L^2}^{q}
	&=	\bar{\E}\max_{0\le m\le N}\|\hM_h(t^m)\|_{L^2}^{q}\\
	&\le\E\sup_{0\le t\le T}	\norm{
		\int_0^t\mathcal P_{\mathbb U_h}G(u_h(s))\,dW(s)}_{L^2}^{q}\\
	&\le	C_q\E	\left(	\tr (Q)\int_0^T	\norm{G(u_h(s))}_{\mathcal L(\mathcal K,L^2)}^2ds	\right)^{q/2}\\
	&\le	C_{q,T}	\E\max_{0\le m\le N}\norm{u_h^m}_{L^2}^{q}+C\\
	&\le C.
	\end{split}
\end{align}
Using that $(u_h,r_h)$ and $(\bar{u}_h,\bar{r}_h)$ have the same laws and that $W$ is a Wiener process, we also have for $v,w\in \Uh$, 
\begin{align*}
	&\bar{\E}\bigg[\bigg((\hM_h(t), v)(\hM_h(t), w)-(\hM_h(s), v)(\hM_h(s),w)\\
	&\quad -\int_{\kappa_h(s)}^{\kappa_h(t)}  \left( Q^{1/2} G(\bar{u}_h)^*\PUh  v,Q^{1/2}G(\bar{u}_h)^*\PUh  w\right) d\tau\bigg)\phi(\bar{u}_h|_{[0,s]})\bigg]\\
	& ={\E}\Bigg[\Bigg(\left(\int_{0}^{\kappa_h(t)}G(u_h)dW(\tau), v\right)\left(\int_{0}^{\kappa_h(t)}G(u_h)dW(\tau), w\right)\\
	&\quad-\left(\int_{0}^{\kappa_h(s)}G(u_h)dW(\tau), v\right)\left(\int_{0}^{\kappa_h(s)}G(u_h)dW(\tau), w\right)\\
	&\quad -\int_{\kappa_h(s)}^{\kappa_h(t)}  \left( Q^{1/2} G({u}_h)^* v,Q^{1/2}G({u}_h)^* w\right) d\tau\Bigg)\phi(u_h|_{[0,s]})\Bigg]\\
	& = {\E}\Bigg[\Bigg(\left(\sum_{m=0}^{m_2-1}G(u_h^m)\Delta W_m, v\right)\left(\sum_{m=0}^{m_2-1}G(u^m_h)\Delta W_m, w\right)\\
	&\quad-\left(\sum_{m=0}^{m_1-1}G(u^m_h)\Delta W_m, v\right)\left(\sum_{m=0}^{m_1-1}G(u^m_h)\Delta W_m, w\right)\\
	&\quad -\int_{\kappa_h(s)}^{\kappa_h(t)}  \left( Q^{1/2} G({u}_h)^*v,Q^{1/2}G({u}_h)^* w\right) d\tau\Bigg)\phi(u_h|_{[0,s]})\Bigg]\\
	& = {\E}\Bigg[\Bigg(\sum_{m,n=0}^{m_2-1}\left(\Delta W_m,G(u_h^m)^* v\right)\left(\Delta W_n,G(u^n_h)^* w\right)\\
	&\quad-\sum_{m,n=0}^{m_1-1}\left(\Delta W_m,G(u^m_h)^* v\right)\left(\Delta W_n,G(u^n_h)^* w\right)\\
	&\quad -\int_{\kappa_h(s)}^{\kappa_h(t)}  \left( Q^{1/2} G({u}_h)^*v,Q^{1/2}G({u}_h)^* w\right) d\tau\Bigg)\phi(u_h|_{[0,s]})\Bigg]\\
	& = {\E}\Bigg[\Bigg(\sum_{m,n=m_1}^{m_2-1}\left(\Delta W_m,G(u_h^m)^* v\right)\left(\Delta W_n,G(u^n_h)^* w\right)\\
	&\quad +\sum_{m=0}^{m_1-1}\sum_{n=m_1}^{m_2-1}\left(\Delta W_m,G(u^m_h)^* v\right)\left(\Delta W_n,G(u^n_h)^*w\right)\\
	&\quad 	+\sum_{m=m_1}^{m_2-1}\sum_{n=0}^{m_1-1}\left(\Delta W_m,G(u^m_h)^* v\right)\left(\Delta W_n,G(u^n_h)^* w\right)\\
	&\quad -\int_{\kappa_h(s)}^{\kappa_h(t)}  \left( Q^{1/2} G({u}_h)^* v,Q^{1/2}G({u}_h)^* w\right) d\tau\Bigg)\phi(u_h|_{[0,s]})\Bigg]\\
	& = {\E}\Bigg[\Bigg(\sum_{m=m_1}^{m_2-1}\left(\Delta W_m,G(u_h^m)^* v\right)\left(\Delta W_m,G(u^m_h)^* w\right)\\
	&\quad -\int_{\kappa_h(s)}^{\kappa_h(t)}  \left( Q^{1/2} G({u}_h)^* v,Q^{1/2}G({u}_h)^* w\right) d\tau\Bigg)\phi(u_h|_{[0,s]})\Bigg],
\end{align*}
using the adaptedness of $u_h^m$, the measurability of $\phi$ and the conditional mean-zero property of the Wiener process which makes the cross terms and all terms with $m\neq n$ vanish.
Next, we observe that since the covariance of the increment satisfies $\mathbb E[(\Delta W_m,a)(\Delta W_m,b)\mid\mathcal G_{t^m}]=\Delta t\,(Qa,b)$ for $m\geq m_1$, we have
\begin{multline*}
	{\E}\Bigg[\left(\Delta W_m,G(u_h^m)^* v\right)\left(\Delta W_m,G(u^m_h)^* w\right)\phi(u_h|_{[0,s]})\Bigg] \\
	= \Delta t \E\Bigg[\left(Q^{1/2}G(u_h^m)^* v,Q^{1/2}G(u^m_h)^* w\right)\phi(u_h|_{[0,s]})\Bigg], 
\end{multline*}
we obtain
\begin{multline*}
	{\E}\Bigg[\sum_{m=m_1}^{m_2-1}\left(\Delta W_m,G(u_h^m)^*v\right)\left(\Delta W_m,G(u^m_h)^* w\right)\phi(u_h|_{[0,s]})\Bigg] \\
	= \Delta t \E\Bigg[\sum_{m=m_1}^{m_2-1}\left(Q^{1/2}G(u_h^m)^* v,Q^{1/2}G(u^m_h)^* w\right)\phi(u_h|_{[0,s]})\Bigg], 
\end{multline*}
and thus for any $v,w\in \Uh$,
\begin{align}\label{eq:Rhquadraticvar}
	\begin{split}
		&\bar{\E}\bigg[\bigg((\hM_h(t), v)(\hM_h(t), w)-(\hM_h(s), v)(\hM_h(s), w)\\
		&\quad -\int_{\kappa_h(s)}^{\kappa_h(t)}  \left( Q^{1/2} G(\bar{u}_h)^* v,Q^{1/2}G(\bar{u}_h)^* w\right) d\tau\bigg)\phi(\bar{u}_h|_{[0,s]})\bigg]\\
		& = 0.
	\end{split}
\end{align}
The process $\bar R_h$ is predictable, since its jump at $t^{m+1}$ is $\bar{\mathcal G}_{t^m}^{h}$-measurable. Together with the preceding quadratic martingale identity, this identifies \(\bar R_h\) as the predictable operator quadratic variation of $\hM_h$.
Next, we pass to the limit in these two expressions. To do so, we use Lemma~\ref{lem:limits}. All the terms in the definition of $\hM_h$ in~\eqref{eq:Mhdef} are of the form of those considered in Lemma~\ref{lem:limits} with the integrals over $[\kappa_h(t),t]$ vanishingly small in $L^\alpha(\bar{\Omega})$ as $h,\Delta t\to 0$. Defining for $v\in C^\infty_{c,\Div}(\dom)$
\begin{equation}\label{eq:Mphidef}
	M^v(t): = (\bar{u}(t),v) - (u_0,v)+\int_0^{t}b(\bar{u}(s),\bar{u}(s),v)ds +\nu \int_0^{t} (\Grad \bar{u}(s),\Grad v) ds-\int_0^{t}\langle f,v\rangle ds,
\end{equation} 
it then follows from Lemma~\ref{lem:limits} that
\begin{equation*}
	\bar{\mathbb{E}}\left|(\hM_h(t),\Ih v)-M^v(t)\right|^\alpha \stackrel{h,\Delta t\to 0}{\longrightarrow } 0,
\end{equation*}
for $1\leq \alpha< 2^{r_0-1}$. Moreover, since $\bar{u}(t)$ has sample paths in $C([0,T];L^2(\dom)_w)$ a.s., and the integrands in this definition are in $L^1([0,T])$, $M^v(t)$ is continuous. 
Using this, we show that $M^v(t)$ is a scalar martingale.
To do so, it is enough to show that
\begin{equation}\label{eq:Mhatmartingale}
	\bar{\mathbb{E}}\left[(M^v(t)-M^v(s))\phi(\bar{u}|_{[0,s]})\right] = 0,
\end{equation}
for all $s\leq t\in [0,T]$ and all $\phi$ bounded and continuous functions on $\mathbb{D}([0,s];H^{-1}(\dom))$. 
\eqref{eq:Mhatmartingale} follows if we can show that
\begin{equation}\label{eq:hMlimit}
	\lim_{h,\Delta t\to 0} 	\bar{\mathbb{E}}\left[(M^v(t)-M^v(s))\phi(\bar{u}|_{[0,s]}) -(\hM_h(t)-\hM_h(s),\Ih v)\phi(\bar{u}_h|_{[0,s]})  \right]  = 0,
\end{equation}
since we have already shown that $\hM_h$ is a discrete martingale. So we compute
\begin{align*}
&	\left|	\bar{\mathbb{E}}\left[(M^v(t)-M^v(s))\phi(\bar{u}|_{[0,s]}) -(\hM_h(t)-\hM_h(s),\Ih v)\phi(\bar{u}_h|_{[0,s]})  \right]  \right|\\
 &\leq 
	\bar{\mathbb{E}}\left[\left|(M^v(t)-(\hM_h(t),\Ih v))\right| \left|\phi(\bar{u}_h|_{[0,s]})\right|  \right]  + \bar{\mathbb{E}}\left[\left|(M^v(s)-(\hM_h(s),\Ih v))\right| \left|\phi(\bar{u}_h|_{[0,s]})\right|  \right] \\
	&  + 	\bar{\mathbb{E}}\left[\left| (M^v(t)-M^v(s))\right| \left|\phi(\bar{u}_h|_{[0,s]})-\phi(\bar{u}|_{[0,s]})\right|\right]\\
	& \leq 	\bar{\mathbb{E}}\left[\left|(M^v(t)-(\hM_h(t),\Ih v))\right|  \right] \norm{\phi}_{\infty}  + \bar{\mathbb{E}}\left[\left|(M^v(s)-(\hM_h(s),\Ih v))\right|  \right]\norm{\phi}_{\infty}  \\
	&  + 	\bar{\mathbb{E}}\left[\left|M^v(t)-M^v(s)\right| \left|\phi(\bar{u}_h|_{[0,s]})-\phi(\bar{u}|_{[0,s]})\right|\right].
\end{align*}
The first two terms on the right-hand side converge to zero by Lemma~\ref{lem:limits} (choose any $1\leq \alpha< 2^{r_0-1}$). The third term converges to zero since $\phi$ is bounded and so the term inside the expectation is dominated by an integrable random variable independently of $h$ and $\Delta t$. Furthermore, since $\bar{u}\in C([0,T];H^{-1}(\dom))$ $\bar{\mathbb{P}}$-a.s. by Lemma~\ref{lem:uhcontinuity}, the Skorokhod convergence $\bar{u}_h\to \bar{u}$ is uniform a.s. (arguing as in the proof of~\eqref{eq:fixedtimelimit}), and therefore $\bar{u}_h|_{[0,s]}\to \bar{u}|_{[0,s]}$ in $\mathbb D([0,s];H^{-1}(D))$ $\bar{\mathbb{P}}$-a.s. Thus, $\phi(\bar{u}_h|_{[0,s]}) \to \phi(\bar{u}|_{[0,s]})$ $\bar{\mathbb{P}}$-a.s. since $\phi$ is continuous on  $\mathbb{D}([0,s];H^{-1}(\dom))$.
Therefore, we can apply the dominated convergence theorem to conclude that this term converges to zero.  
Thus~\eqref{eq:hMlimit} holds. 
Next, we aim to identify the quadratic variation of $M^v$. To do so, we first define $\tR$ by
\begin{equation*}
	\tR(t) = \int_0^t G(\bar{u}(\tau))Q^{1/2}(G(\bar{u}(\tau))Q^{1/2})^* d\tau.
\end{equation*}
We want to show that
\begin{equation}\label{eq:Rquadraticvar}
	\bar{\E}\bigg[\left(M^v(t)M^w(t)-M^v(s)M^w(s) - \left\langle(\tR(t)-\tR(s))v,w\right\rangle\right)	\phi(\bar{u}|_{[0,s]})\bigg] = 0,
\end{equation}
for all $s\leq t\in [0,T]$, all $v,w\in C^\infty_{c,\Div}(\dom)$ and all $\phi$ bounded and continuous functions on $\mathbb{D}([0,s];H^{-1}(\dom))$.
This follows if we can show that 
\begin{multline*}
	\lim_{h,\Delta t\to 0} \Bigg( 	\bar{\E}\bigg[\left(M^v(t)M^w(t)-M^v(s) M^w(s) - \left\langle(\tR(t)-\tR(s))v,w\right\rangle\right)	\phi(\bar{u}|_{[0,s]})\bigg] \\
	-  \bar{\E}\bigg[\bigg((\hM_h(t),\Ih v)(\hM_h(t),\Ih w)-(\hM_h(s),\Ih v)(\hM_h(s),\Ih w)\\
	-\left((\hR_h(t)-\hR_h(s))\Ih v,\Ih w\right)\bigg)\phi(\bar{u}_h|_{[0,s]})\bigg] \Bigg) =0.
\end{multline*}
To do so, we bound the terms as follows:
\begin{align*}
	&\Bigg|	\bar{\E}\bigg[\left(M^v(t)M^w(t)-M^v(s)M^w(s) - \left\langle(\tR(t)-\tR(s))v,w\right\rangle\right)	\phi(\bar{u}|_{[0,s]})\bigg] \\
&\quad 	-  \bar{\E}\bigg[\bigg((\hM_h(t),\Ih v)(\hM_h(t),\Ih w)-(\hM_h(s),\Ih v)(\hM_h(s),\Ih w)\\
&\qquad \qquad 	-\left((\hR_h(t)-\hR_h(s))\Ih v,\Ih w\right)\bigg)\phi(\bar{u}_h|_{[0,s]})\bigg] \Bigg|\\
& \leq  \underbrace{\bar{\E}\left[\left|(\hM_h(t),\Ih v)(\hM_h(t),\Ih w) -M^v(t)M^w(t) \right| \left|\phi(\bar{u}_h|_{[0,s]})\right|\right]}_{(i)}\\
&\quad + \underbrace{\bar{\E}\left[\left|(\hM_h(s),\Ih v)(\hM_h(s),\Ih w) - M^v(s)M^w(s) \right| \left|\phi(\bar{u}_h|_{[0,s]})\right|\right]}_{(ii)}\\
& \quad +  \underbrace{\bar{\E}\left[\left| \left( (\hR_h(t)-\hR_h(s))\Ih v,\Ih w\right) - \left\langle(\tR(t)-\tR(s))v,w\right\rangle\right|\left|\phi(\bar{u}_h|_{[0,s]})\right|\right]}_{(iii)}\\
& \quad +	\underbrace{\bar{\E}\bigg[\left|M^v(t)M^w(t)-M^v(s)M^w(s) - \left\langle(\tR(t)-\tR(s))v,w\right\rangle\right|	\left|\phi(\bar{u}|_{[0,s]}) -\phi(\bar{u}_h|_{[0,s]})  \right|\bigg]}_{(iv)}. 
\end{align*}
We estimate the terms on the right-hand side. For $(i)$, we have
\begin{align*}
	(i)&\leq \bar{\E}\left[\left|(\hM_h(t),\Ih v)(\hM_h(t),\Ih w) - M^v(t)M^w(t) \right| \right] \norm{\phi}_\infty\\
	& \leq \left(\bar{\E}\left|(\hM_h(t),\Ih v)-M^v(t)\right|^2\right)^{1/2}\left(\bar{\E}\left|(\hM_h(t),\Ih w)\right|^2\right)^{1/2}\norm{\phi}_\infty\\
	&\quad +\left(\bar{\E}\left|(\hM_h(t),\Ih w)-M^w(t)\right|^2\right)^{1/2}\left(\bar{\E}\left|M^v(t)\right|^2\right)^{1/2}\norm{\phi}_\infty\\
	& \leq C \left(\left(\bar{\E}\left|(\hM_h(t),\Ih v)-M^v(t)\right|^2\right)^{1/2} + \left(\bar{\E}\left|(\hM_h(t),\Ih w)-M^w(t)\right|^2\right)^{1/2}\right).
\end{align*}
The right-hand side goes to zero by Lemma~\ref{lem:limits} (choosing $\alpha=2$).
Term $(ii)$ is shown to converge to zero in the same way. For term $(iii)$, we have
\begin{align*}
	(iii)& = \bar{\E}\left[\left| \int_{\kappa_h(s)}^{\kappa_h(t)} \!\!
	 \left( Q^{1/2} G(\bar{u}_h)^*\PUh \Ih v,Q^{1/2}G(\bar{u}_h)^*\PUh \Ih  w\right) d\tau\!\! -\! \int_s^t\!\! \left( Q^{1/2} G(\bar{u})^*  v,Q^{1/2}G(\bar{u})^* w\right) d\tau \right|\left|\phi(\bar{u}_h|_{[0,s]})\right|\right]\\
	 & \leq  \bar{\E}\left[\! \int_s^t\!\! \left|\left( Q^{1/2} G(\bar{u}_h)^*\Ih v,Q^{1/2}G(\bar{u}_h)^*\Ih w\right) -\left( Q^{1/2} G(\bar{u})^*  v,Q^{1/2}G(\bar{u})^* w\right)\right| d\tau \right]\norm{\phi}_\infty\\
	 &\quad + \bar{\E}\left[ \int_{\kappa_h(s)}^{s} \left|
	 \left( Q^{1/2} G(\bar{u}_h)^*\Ih v,Q^{1/2}G(\bar{u}_h)^*\Ih w\right) \right| d\tau \right]\norm{\phi}_\infty \\
	 &\quad +  \bar{\E}\left[ \int_{\kappa_h(t)}^{t} \left|
	 \left( Q^{1/2} G(\bar{u}_h)^*\Ih v,Q^{1/2}G(\bar{u}_h)^*\Ih w\right) \right| d\tau \right] \norm{\phi}_\infty\\
	 & \leq C \bar{\E}\int_s^t\!\! \left|\left( Q^{1/2} \left(G(\bar{u}_h)^*-G(\bar{u})^*\right)\Ih v,Q^{1/2}G(\bar{u}_h)^*\Ih w\right)\right| d\tau \\
	 &\quad+ C \bar{\E} \int_s^t\!\! \left|\left( Q^{1/2} G(\bar{u})^*(\Ih v-v),Q^{1/2}G(\bar{u}_h)^*\Ih w\right) \right| d\tau\\
	 &\quad + C \bar{\E}\int_s^t \left|\left( Q^{1/2} G(\bar{u})^*v,Q^{1/2}\left(G(\bar{u}_h)^*-G(\bar{u})^*\right)\Ih w\right)  \right| d\tau\\
	 &\quad 
	  +C \bar{\E}\int_s^t \left|\left( Q^{1/2} G(\bar{u})^*v,Q^{1/2}G(\bar{u})^*\left(\Ih w-w\right)\right)  \right| d\tau \\
	 &\quad + C\bar{\E} \int_{\kappa_h(s)}^{s} \norm{G(\bar{u}_h)}^2_{\mathcal{L}(\mathcal{K},L^2)}\norm{\Ih v}_{L^2}\norm{\Ih w}_{L^2} d\tau\\
	 &\quad +  C\bar{\E} \int_{\kappa_h(t)}^{t} \norm{G(\bar{u}_h)}^2_{\mathcal{L}(\mathcal{K},L^2)}\norm{\Ih v}_{L^2}\norm{\Ih w}_{L^2} d\tau  \\
	  & \leq C \bar{\E}\int_s^t  \norm{G(\bar{u}_h)-G(\bar{u})}_{\mathcal{L}(\mathcal{K},L^2)}\norm{G(\bar{u}_h)}_{\mathcal{L}(\mathcal{K},L^2)}\norm{\Ih v}_{L^2}\norm{\Ih w}_{L^2}d\tau \\
	 &\quad+ C \bar{\E} \int_s^t \norm{G(\bar{u}_h)}_{\mathcal{L}(\mathcal{K},L^2)}\norm{G(\bar{u})}_{\mathcal{L}(\mathcal{K},L^2)}\norm{\Ih v -v}_{L^2}\norm{\Ih w}_{L^2} d\tau\\
	 &\quad + C \bar{\E}\int_s^t \norm{G(\bar{u})}_{\mathcal{L}(\mathcal{K},L^2)}\norm{v}_{L^2}\norm{G(\bar{u}_h)-G(\bar{u})}_{\mathcal{L}(\mathcal{K},L^2)}\norm{\Ih w}_{L^2} d\tau\\
	 &\quad 	 +C \bar{\E}\int_s^t \norm{G(\bar{u})}_{\mathcal{L}(\mathcal{K},L^2)}^2\norm{v}_{L^2}\norm{\Ih w-w}_{L^2} d\tau \\
	 &\quad + CK_1^2\bar{\E} \int_{\kappa_h(s)}^{s} \left(\norm{ \bar{u}_h}^2_{L^2}+1\right)d\tau \norm{ v}_{H^2}\norm{ w}_{H^2}+ CK_1^2\bar{\E} \int_{\kappa_h(t)}^{t} \left(\norm{ \bar{u}_h}^2_{L^2}+1\right)d\tau \norm{ v}_{H^2}\norm{ w}_{H^2}  \\
	  & \leq CK_2K_1\bar{\E}\int_s^t \norm{ \bar{u}_h-\bar{u}}_{L^2}\left(\norm{\bar{u}_h}_{L^2}+1\right)d\tau\norm{ v}_{H^2}\norm{w}_{H^2} \\
	  &\quad + C K_1^2 h \bar{\E} \int_s^t \left(\norm{\bar{u}_h}_{L^2}^2 +\norm{\bar{u}}_{L^2}^2+1\right)d\tau\norm{v}_{H^2}\norm{ w}_{H^2} \\
	 &\quad + CK_1 K_2 \bar{\E}\int_s^t \left(\norm{\bar{u}}_{L^2}+1\right)\norm{ \bar{u}_h-\bar{u}}_{L^2}d\tau \norm{w}_{H^2}  \norm{v}_{H^2}\\
	 &\quad 	 +C h K_1^2\bar{\E}\int_s^t \left(\norm{\bar{u}}_{L^2}^2+1\right)d\tau \norm{v}_{H^2}\norm{w}_{H^2}   + CK_1^2\Delta t\bar{\E} \left( \norm{ \bar{u}_h}^2_{L^\infty(0,T;L^2(\dom))}+1\right) \norm{ v}_{H^2}\norm{ w}_{H^2}\\
	   & \leq CK_2K_1\left(\bar{\E}\int_s^t  \norm{ \bar{u}_h-\bar{u}}_{L^2}^2 d\tau\right)^{1/2}\left(\bar{\E}\int_s^t  \left(\norm{\bar{u}_h}_{L^2}^2 + \norm{\bar{u}}_{L^2}^2+1\right)d\tau\right)^{1/2} \\
	   &\quad + C K_1^2 h \bar{\E} \int_s^t (\norm{\bar{u}_h}^2_{L^2}+\norm{\bar{u}}_{L^2}^2)d\tau  +CK_1^2\Delta t\bar{\E}  \norm{ \bar{u}_h}^2_{L^\infty(0,T;L^2(\dom))}+C(h+\Delta t),
\end{align*}
where we used the assumptions on $G$,~\eqref{eq:lingrowthG} and~\eqref{eq:Glipschitz} and the stability of the $L^2$-projection. Now using that $\bar{u}_h\to \bar{u}$ a.s. in $L^2([0,T]\times\dom)$ and the bounds from Lemma~\ref{lem:energyestimate}, we see that the right-hand side converges to zero. For the fourth term, $(iv)$, we use that $\phi(\bar{u}_h|_{[0,s]})\to \phi(\bar{u}|_{[0,s]})$ a.s. since $\bar{u}_h\to \bar{u}$ in $\mathbb{D}([0,T];H^{-1}(\dom))$ a.s. and $\phi$ is bounded.
Hence the term converges to zero by the dominated convergence theorem.
 So we obtain
 \begin{equation}\label{eq:Rlimit}
 	\bar{\E}\bigg[\left(M^v(t)M^w(t)-M^v(s)M^w(s) - \left\langle(\tR(t)-\tR(s))v,w\right\rangle\right)	\phi(\bar{u}|_{[0,s]})\bigg]=0.
 \end{equation} 
Identities~\eqref{eq:Mhatmartingale} and~\eqref{eq:Rlimit} are valid for any $\phi$ bounded and continuous on  $\mathbb{D}([0,s];H^{-1}(\dom))$. 
The square integrability of $M^v$ follows from the moment bound~\eqref{eq:momentboundM} for $q\geq 2$. 
Next, we let $V=H^1_{\Div}(\dom)$ and note that 
\begin{equation}\label{eq:apriori}
\sup_{0\leq t\leq T}	\left|M^v(t)\right| \leq C_T(\omega)\norm{v}_V,
\end{equation}
where $C_T(\omega)$ is given by
\begin{equation*}
	C_T(\omega) = C\left(\sup_{0\leq t\leq T}	\norm{\bar{u}(t)-u_0}_{H^{-1}(\dom)} +\nu\sqrt T \norm{\bar u}_{L^2(0,T;V)}+C\norm{\bar u}_{L^2(0,T;V)}^2+ \sqrt T\norm{f}_{L^2(0,T;V^*)}\right).
\end{equation*}
By the estimates from Proposition~\ref{prop:skorokhod}, $C_T(\omega)<\infty$ almost surely. By Lemma~\ref{lem:uhcontinuity}, $M^v(t)$ is almost surely continuous in time. On a set of full $\bar{\PP}$-measure, the linear map $v\mapsto M^v(t)$ has a unique extension to $V$ by the density of $C^\infty_{c,\Div}(\dom)$ in $V$. Then we can define $\hM(t)\in V^*$ via
\begin{equation*}
	\langle \hM(t),v\rangle_{V^*,V} : = \lim_{n\to \infty} M^{v_n}(t),\quad v_n \in C^\infty_{c,\Div}(\dom),\quad v_n\to v\,\quad \text{in}\,\, V,
\end{equation*}
for any $v\in V$.
The estimate~\eqref{eq:apriori} makes the limit independent of the chosen sequence and yields
\begin{equation*}
	\sup_{0\leq t\leq T} \norm{\hM(t)}_{V^*}\leq C_T(\omega),
\end{equation*}
almost surely and moreover, $\hM(t)$ satisfies the formula~\eqref{eq:Mphidef}. Since 
\begin{equation*}
	\norm{\hM(t)-\hM(s)}_{V^*}\leq \norm{\bar{u}(t)-\bar{u}(s)}_{H^{-1}(\dom)}+ \nu \int_s^t \norm{\bar{u}(r)}_V dr + C \int_s^t \norm{\bar{u}(r)}^2_V dr + \int_s^t \norm{f(r)}_{H^{-1}(\dom)} dr,
\end{equation*}
it follows that $\hM \in C([0,T];V^*)$ a.s. Thanks to the estimates from Proposition~\ref{prop:skorokhod} (v), we have
\begin{equation*}
	\bar{\E}\sup_{0\leq t\leq T} \norm{\hM(t)}_{V^*}^q\leq \bar{\E}C_T^q <\infty,\quad 1\leq q \leq 2^{r_0-1}.
\end{equation*}
This implies that we can also extend the martingale properties~\eqref{eq:Mhatmartingale} and~\eqref{eq:Rquadraticvar} to any $v,w\in H^1_{\Div}(\dom)$.
Therefore $\hM$ is a $V^*$-valued square integrable martingale with respect to the filtration $\bar{\mathcal{G}}^0_t = \sigma(\{\bar{u}(s),\, s\leq t\})$. We denote by $(\bar{\mathcal G}_t)_{t\in[0,T]}$ the usual augmentation of $(\bar{\mathcal G}_t^0)_{t\in[0,T]}$.
Thanks to continuity and the preceding integrability bounds, the linear and compensated quadratic martingale identities remain valid with respect to $(\bar{\mathcal G}_t)_{t\in[0,T]}$.
To identify the quadratic variation, we proceed as in~\cite[Lemma C2, Theorem C3]{Brzezniak2013}.  
We note first that~\eqref{eq:Rlimit} implies that $N^{v,w}(t)$ defined by
\begin{equation*}
	N^{v,w}(t):= \langle\hM(t),v\rangle \langle\hM(t),w\rangle - \left\langle\tR(t)v,w\right\rangle,
\end{equation*}
is a scalar martingale for any $v,w\in V$. We can view $\tR(t)$ as an operator from $V$ to $V^*$. We let $I:V\to V^*$ be the Riesz isomorphism defined by $I^{-1}(f)=u$ where for $f\in V^*$, $u$ is the solution of 
\begin{equation*}
	(u,v)_V = \langle f,v\rangle_{V^*,V} = \langle I(u),v\rangle,
\end{equation*}
which exists and is unique by the Lax--Milgram lemma. One can then define an inner product on $V^*$ by
\begin{equation*}
	(f,g)_{V^*}:=(I^{-1}(f),I^{-1}(g))_V,\quad \forall \, f,g\in V^*.
\end{equation*}
By~\cite[Lemma C2]{Brzezniak2013}, it holds 
\begin{equation*}
	\langle f,v\rangle_{V^*,V} = (f,I(v))_{V^*},\quad \forall\, f\in V^*,\quad \forall \, v\in V.
\end{equation*}
Now we define $\hR=\tR\circ I^{-1}:V^*\to V^*$. Then $\hR$ satisfies for $v,w\in V$,
\begin{equation*}
	\left\langle \tR(t)v,w\right\rangle = \left\langle \hR(t)I(v),w\right\rangle = (\hR(t)I(v),I(w))_{V^*}=(\hM(t),I(v))_{V^*}(\hM(t),I(w))_{V^*}-N^{v,w}(t).
\end{equation*}
Since $I:V\to V^*$ is an isomorphism and therefore
is onto, we have  for arbitrary $a,b\in V^*$, upon setting $a:=I(v)$  and $b=I(w)$ in this identity,
\begin{equation*}
	(\hR(t)a,b)_{V^*}=(\hM(t),a)_{V^*}(\hM(t),b)_{V^*} - N^{I^{-1}(a),I^{-1}(b)}(t).
\end{equation*}
Since $\hR$ 
is a predictable, increasing, positive trace-class-valued process with $\hR(0)=0$, which follows from its integral representation, the Hilbert-space characterization of quadratic variation implies
\begin{equation*}
	\hR(t) = \langle\!\langle \hM \rangle\!\rangle_t,
\end{equation*}
i.e., $\hR$ is the quadratic variation of $\hM$.
Equivalently, $\tR(t) = \langle\!\langle \hM\rangle\!\rangle_t\circ I$, cf. Theorem C3 in~\cite{Brzezniak2013}. Interpreting $G(\bar{u})k$ for $k\in \mathcal{K}$ as an element of $V^*$, we can now apply the martingale representation theorem~\cite[Theorem 8.2]{DaPrato2014} for the Hilbert space $V^*$ to infer the existence of a filtered  probability space  $(\widetilde{\Omega},\widetilde{\mathcal{G}}_t,\widetilde{\mathcal{G}},\widetilde{\mathbb{P}})$ and a $Q$-Wiener process $\bar{W}$ with values in $\mathcal{K}$ defined on $(\bar{\Omega}\times\widetilde{\Omega},\bar{\mathcal{G}}\otimes\widetilde{\mathcal{G}},\bar{\mathbb{P}}\otimes\widetilde{\mathbb{P}})$ adapted to $(\bar{\mathcal{G}}_t\otimes\widetilde{\mathcal{G}}_t)_{t\in [0,T]}$ such that
\begin{equation*}
	\hM(t,\omega,\widetilde{\omega}) = \int_0^tG(\bar{u}(s,\omega,\widetilde{\omega}))d\bar{W}(s,\omega,\widetilde{\omega}),\quad t\in [0,T],\quad  \text{a.s. } (\omega,\widetilde{\omega})\in \bar{\Omega}\times\widetilde{\Omega},
\end{equation*}
where 
\begin{equation*}
	\hM(t,\omega,\widetilde{\omega}) =\hM(t,\omega),\quad \bar{u}(t,\omega,\widetilde{\omega}) =\bar{u}(t,\omega) ,\quad t\in [0,T],\quad \text{a.s. } (\omega,\widetilde{\omega})\in \bar{\Omega}\times\widetilde{\Omega}.
\end{equation*}
$\bar{u}$ satisfies
\begin{equation*}
		(\bar{u}(t),v)-(u_0,v)+\nu\int_0^t(\Grad \bar{u},\Grad v) ds -\int_0^t\int_{\dom} (\bar{u}\cdot\Grad)v \cdot \bar{u} dx ds = \int_0^t \langle f,v\rangle  ds +\int_0^t (G(\bar{u}(s))d\bar{W}(s),v)
	\end{equation*}
	for any $v\in H^1_{\Div}(\dom)$, $\bar{\mathbb{P}}\otimes\widetilde{\mathbb{P}}$-a.s.	Let $(\bar{\bar{\Omega}},\bar{\bar{\mathcal G}},\bar{\bar{\mathbb P}})$ denote the completion of the product probability space
	$(\bar{\Omega}\times\tilde{\Omega},\bar{\mathcal G}\otimes\tilde{\mathcal G},\bar{\mathbb P}\otimes\tilde{\mathbb P}),$ and let $(\bar{\bar{\mathcal G}}_t)_{t\in[0,T]}$ denote the usual augmentation of the product filtration $(\bar{\mathcal G}_t\otimes\tilde{\mathcal G}_t)_{t\in[0,T]}$.
	Then $\bigl(\bar{\bar{\Omega}},\bar{\bar{\mathcal G}},(\bar{\bar{\mathcal G}}_t)_{t\in[0,T]},	\bar{\bar{\mathbb P}},\bar u,\bar W\bigr)$
 is a martingale solution of~\eqref{eq:SNSm}. 
\end{proof}
\section{Numerical experiments}\label{sec:numexp}
In this section, we test the performance of the numerical method~\eqref{eq:predictionfully} - \eqref{eq:projectionfully} in several numerical experiments. The code uses the FEniCSx software and can be found at \url{https://github.com/fraenschii/SNS_projection_method}. For the following experiments we will use the physical domain $\dom = [0,1]^2$ with homogeneous Dirichlet boundary conditions. We set the deterministic forcing $f=0$, and  let $\Delta t = T/N$ as before. We use lowest order Taylor--Hood (P2/P1) finite elements for all the following numerical examples. The mesh is a uniform $n\times n$ subdivision into $2n^2$ identical triangles, so $h=\sqrt{2}/n$.
\subsection{Convergence test}\label{sec:convtest}
Our analysis of the preceding sections only provides convergence along a subsequence to a martingale solution, however, in 2D and for sufficiently smooth initial data, convergence with a rate can be expected (cf.~\cite{Hausenblas2019} where error estimates for the time discretization were proved and~\cite{Feng2023,Carelli2012b} for error estimates for the projection method for the stochastic Stokes equations). We investigate the empirical spatial and temporal convergence behavior in the following example. We set the deterministic forcing $f=0$, choose $G(u)dW(t)=\sigma u d\beta(t)$, where $\sigma =0.3$ and $\beta$ is a one-dimensional real-valued Brownian motion. As initial data, we choose
\begin{equation}\label{eq:init1}
	u_0=\nabla^\perp\psi_0: = (\partial_{x_2} \psi_0, -\partial_{x_1} \psi_0),\quad  \psi_0 = \sin(\pi x_1)^2 \sin(\pi x_2)^2, \quad x=(x_1,x_2)\in \dom = [0,1]^2.
\end{equation}
 We set $T=0.5$ and $\nu=0.01$ for both the spatial and the temporal refinement test. We compute an approximation of the error at time $T$ using the error estimator
\begin{equation}
	\label{eq:numerror}
		\widehat E_{h,\Delta t}=	\left(\frac1J\sum_{j=1}^{J}	\norm{u_{h,\Delta t}^{(j)}(T)
	-u_{{\rm ref}}^{(j)}(T)}_{L^2}^{2}	\right)^{1/2}.
\end{equation}
Here $J$ is the number of independent samples and $u_{h,\Delta t}^{(j)}$ and $u_{{\rm ref}}^{(j)}$ are computed along the same Brownian path. $u_{{\rm ref}}^{(j)}$ is a reference solution computed with a finer resolution. In contrast to the previous section, we index the approximation by $h$ and $\Delta t$ instead of just $h$ to make explicit which parameter is refined.

\subsubsection{Temporal refinement}
We first check the temporal convergence rates. We compute a reference solution with $n=32$ and $N=4096$ time steps. Then we compute approximations with $n=32$ (i.e., $h$ stays fixed) and $N=16, 32,64,128,256$. We use $J=100$ samples. The errors $\widehat{E}_{h,\Delta t_i}$ are reported in Table~\ref{tab:temporal_convergence}. We compute estimates of the convergence rates using the formula
\begin{equation*}
r_i=	\frac{\log(\widehat{E}_{h,\Delta t_i}/\widehat{E}_{h,\Delta t_{i-1}})}{ \log(\Delta t_i / \Delta t_{i-1})}
\end{equation*}
\begin{table}
	\setlength{\tabcolsep}{6pt}
	\begin{tabular}[t]{cccc}
		\toprule
			\multicolumn{4}{c}{Temporal convergence} \\
		\midrule
	$h$ & 	$\Delta t$ & $\widehat{E}_{h,\Delta t}$ & Rate \\
		\midrule
		0.0442&$2^{-5}$  &	0.870 & --\\
		0.0442&$2^{-6}$	& 0.690 & 0.336\\
		0.0442&$2^{-7}$ 	& 0.500 & 0.463\\
		0.0442&$2^{-8}$ & 0.328 & 0.607\\
		0.0442&$2^{-9}$  &	0.195 & 0.754\\
		 \bottomrule
	\end{tabular}
	 \qquad
	\begin{tabular}[t]{cc c c}
		\toprule
		\multicolumn{4}{c}{Spatial convergence} \\
		\midrule
		$h$& $\Delta t$ & $\widehat{E}_{h,\Delta t}$ & Rate \\
		\midrule
		0.3536 & $2^{-9}$ &4.57E-01	 & --\\
		0.1768 & $2^{-9}$ & 8.83E-02&	2.36\\
		 0.0884 & $2^{-9}$ & 5.78E-03	&3.93  \\
		0.0442  & $2^{-9}$ & 6.96E-04 &	3.06 \\
		\bottomrule
	\end{tabular}
	\centering
	\caption{Temporal and spatial convergence rates at final time $T=0.5$ for the initial data~\eqref{eq:init1}.}
	\label{tab:temporal_convergence}
\end{table}
We see that the errors decrease as $\Delta t$ decreases at a rate of approximately 0.54 on average. This is faster than the rates proved in~\cite{Hausenblas2019,Feng2023} but we recall that those rates are a worst-case scenario and convergence can be faster under favorable conditions. Furthermore, we are testing convergence towards a reference solution computed on a finer grid, which can be faster than convergence to the (unknown) exact solution. 
\subsubsection{Spatial refinement}
For the spatial convergence test, we compute a reference solution with $n=128$ and $N=256$ time steps. Then we compute approximations with $N=256$ time steps and $n=4,8,16,32$. We estimate the error~\eqref{eq:numerror} with $J=100$ samples. The errors and estimated convergence rates
\begin{equation*}
	r_i=	\frac{\log(\widehat{E}_{h_i,\Delta t}/\widehat{E}_{h_{i-1},\Delta t})}{ \log(h_i /  h_{i-1})}
\end{equation*}
are reported in Table~\ref{tab:temporal_convergence}. We observe a convergence rate with respect to the $L^2$-spatial norm of the corrected velocity of about 3.12. This is consistent with the formal third-order ($L^2$)-approximation accuracy of quadratic finite elements for sufficiently smooth velocity fields. The reference solution uses the same time step as the approximations, so this experiment primarily measures the effect of the spatial refinement.
\subsection{Vortex dynamics}
\begin{figure}[t]
	\centering
	\begin{subfigure}[b]{1.0\textwidth}
		\includegraphics[width=\textwidth,trim={0 0 0 11mm},clip]{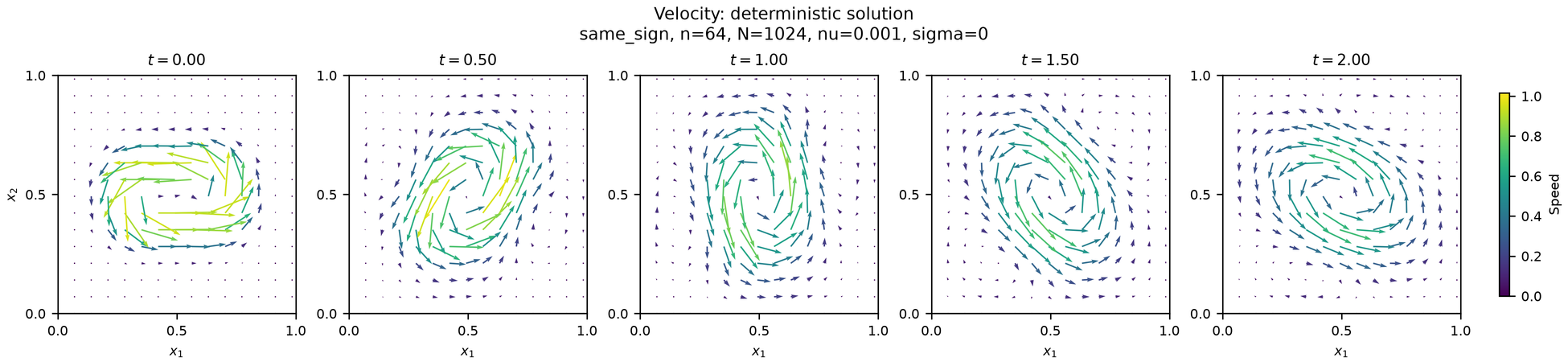}
		\caption{velocity $u_h$}
	\end{subfigure}
	
	\vspace{0.2cm}
	
	\begin{subfigure}[b]{1\textwidth}
		\includegraphics[width=\textwidth,trim={0 0 0 11mm},clip]{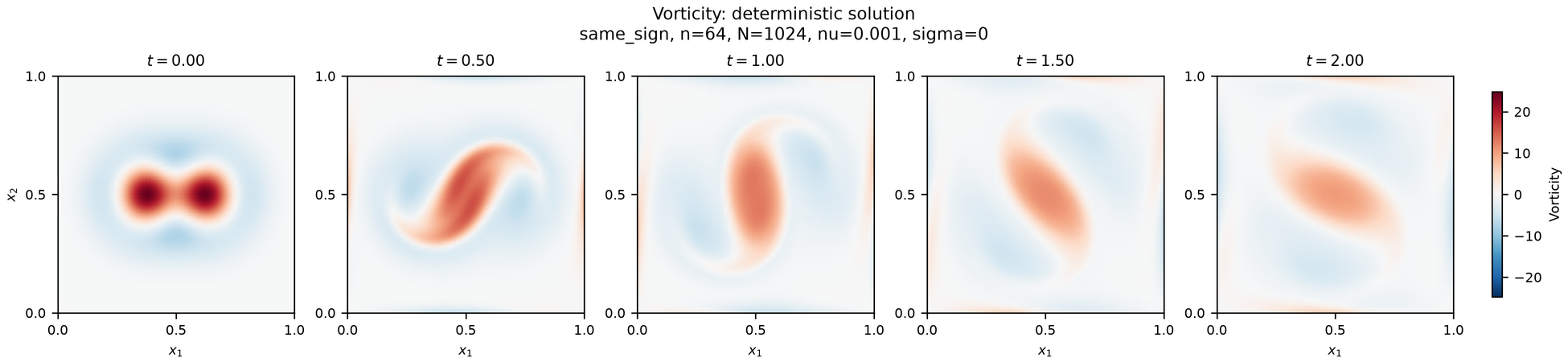}
		\caption{vorticity $\omega_h$}
	\end{subfigure}

	\caption{Approximation of the deterministic solution of~\eqref{eq:SNSm} for initial data~\eqref{eq:init2} (setting $G(u)=0$).}
	\label{fig:exp2_det}
\end{figure}

\begin{figure}[t]
	\centering
	\begin{subfigure}[b]{1.0\textwidth}
		\includegraphics[width=\textwidth,trim={0 0 0 11mm},clip]{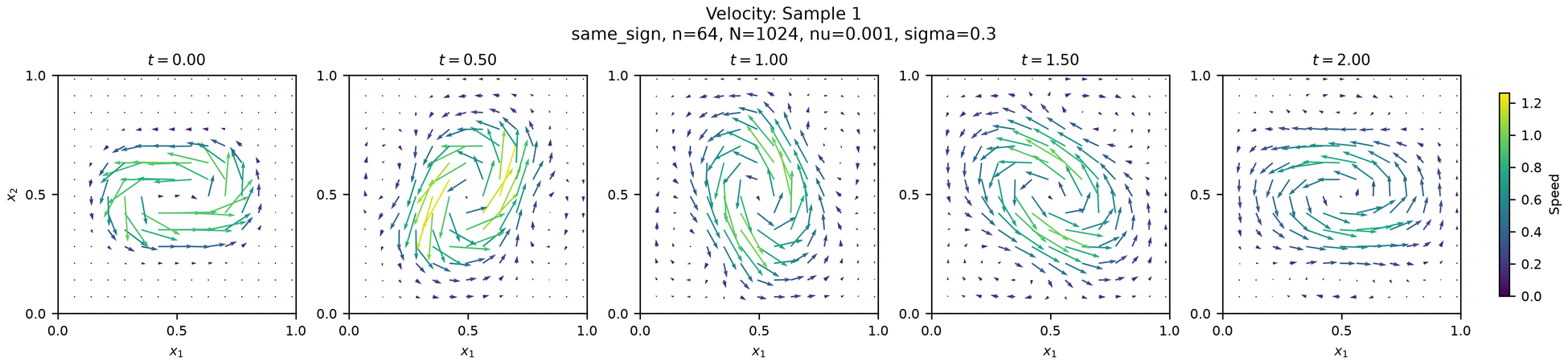}
		\caption{velocity $u_h$}
	\end{subfigure}
	
	\vspace{0.2cm}
	
	\begin{subfigure}[b]{1\textwidth}
		\includegraphics[width=\textwidth,trim={0 0 0 11mm},clip]{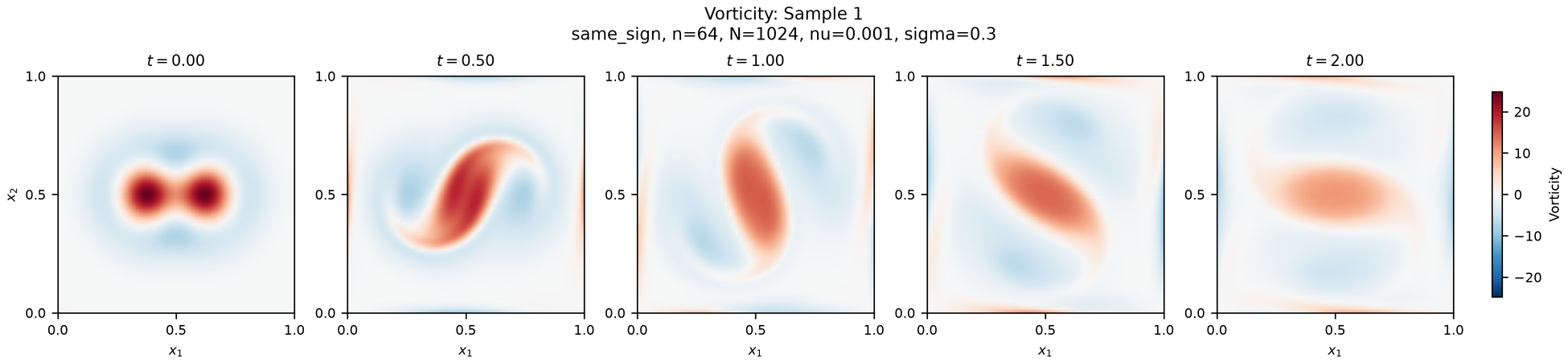}
		\caption{vorticity $\omega_h$}
	\end{subfigure}
	\caption{Sample of the approximation of~\eqref{eq:SNSm} for initial data~\eqref{eq:init2} with $G(u)=\sigma u$, $\sigma=0.3$.}
	\label{fig:exp2_sample}
\end{figure}

\begin{figure}[t]
	\centering
	\begin{subfigure}[b]{1.0\textwidth}
		\includegraphics[width=\textwidth,trim={0 0 0 11mm},clip]{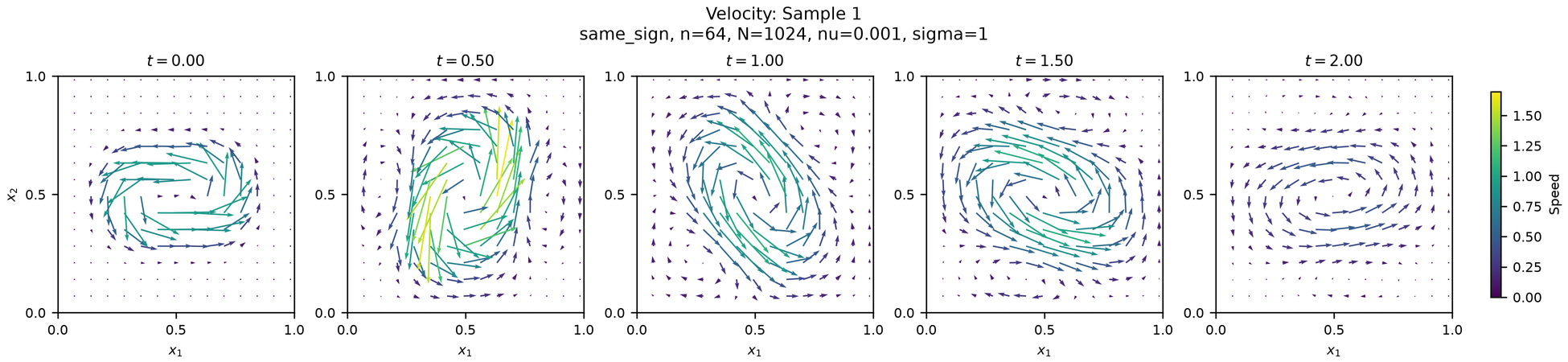}
		\caption{velocity $u_h$}
	\end{subfigure}	
	
	\vspace{0.2cm}
	
	\begin{subfigure}[b]{1\textwidth}
		\includegraphics[width=\textwidth,trim={0 0 0 11mm},clip]{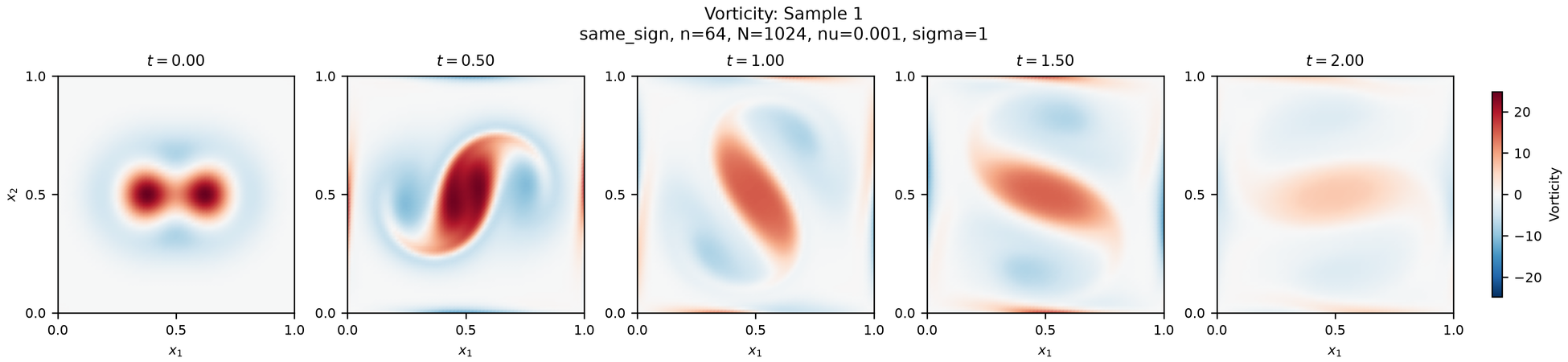}
		\caption{vorticity $\omega_h$}
	\end{subfigure}

	\caption{Sample of the approximation of~\eqref{eq:SNSm} for initial data~\eqref{eq:init2} with $G(u)=\sigma u$, $\sigma=1$.}
	\label{fig:exp2_samplesigma1}
\end{figure}

Next, we consider the influence of the noise on vortex dynamics. We take the initial data
\begin{equation}\label{eq:init2}
	u_0=\bigl(\partial_{x_2}\psi_0,-\partial_{x_1}\psi_0\bigr),\quad		\psi_0(x_1,x_2)=
		C\sin^2(\pi x_1)\sin^2(\pi x_2)	\left[	e^{-\frac{(x_1-0.35)^2+(x_2-0.5)^2}{2(0.12)^2}}	+e^{-\frac{(x_1-0.65)^2+(x_2-0.5)^2}{2(0.12)^2}}	\right]	,
\end{equation}
where $C=0.1843991698438689$ (to approximately normalize the maximum initial speed), and we chose final time $T=2$, viscosity $\nu = 0.001$, spatial subdivisions $n=64$, and $N=1024$ time steps. The approximation of the deterministic solution ($G(u)=0$) is shown in Figure~\ref{fig:exp2_det}. 
Then we add multiplicative noise of the same form as in the previous experiment, $G(u) dW = \sigma u d\beta$, first with $\sigma=0.3$ and then with $\sigma=1$. In Figures~\ref{fig:exp2_sample} and~\ref{fig:exp2_samplesigma1} the velocity field and the approximation $\omega_h$ of the vorticity $\omega = \partial_{x_1} u_2-\partial_{x_2} u_1$ for a single sample are shown at times $t=0, 0.5,1,1.5,2$ for $\sigma =0.3$ and $\sigma=1$ respectively. 

\begin{figure}[t]
	\centering
	\begin{subfigure}[b]{1.0\textwidth}
		\includegraphics[width=\textwidth,trim={0 0 0 11mm},clip]{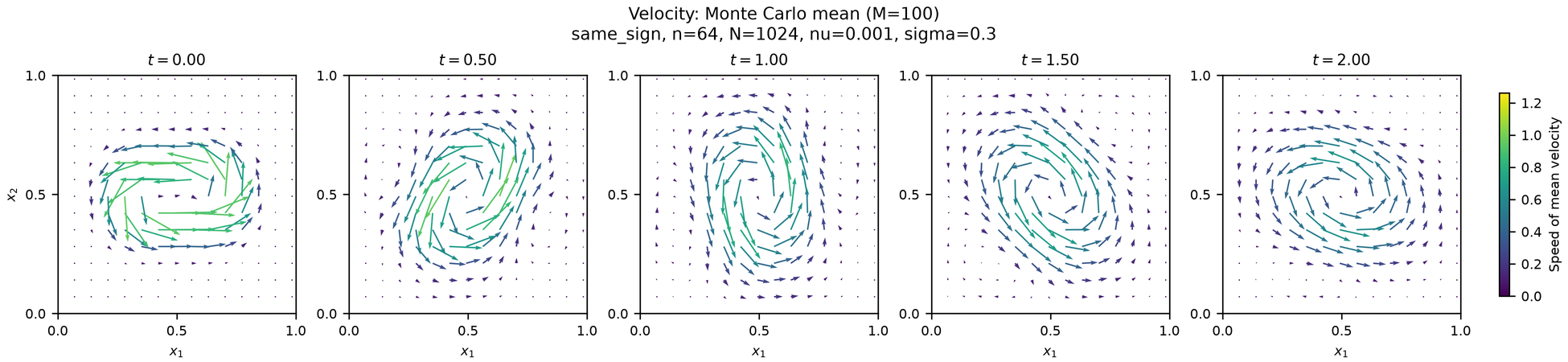}
		\caption{velocity mean $\widehat{E}_J u_h$}
	\end{subfigure}
	
		\vspace{0.2cm}
	
	\begin{subfigure}[b]{1.0\textwidth}
		\includegraphics[width=\textwidth,trim={0 0 0 11mm},clip]{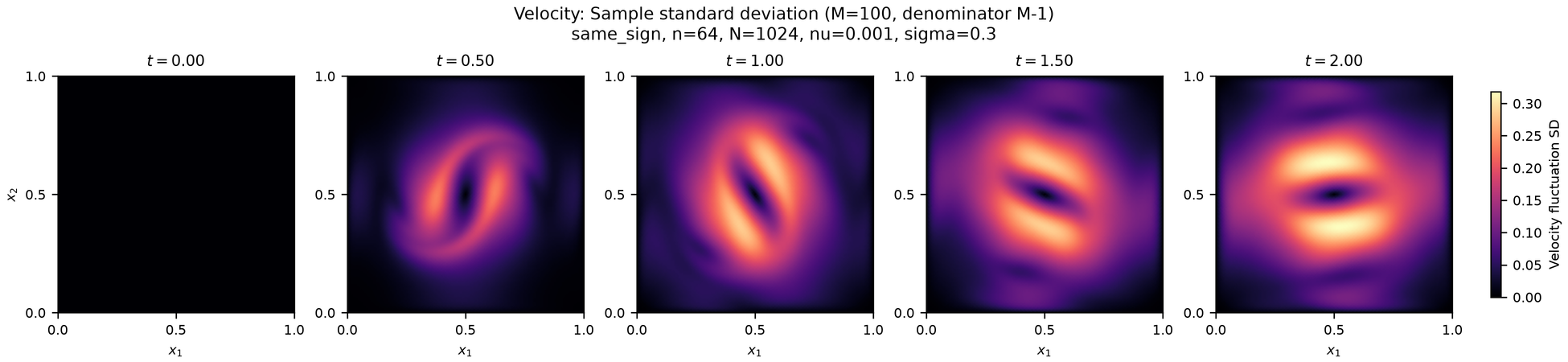}
		\caption{velocity fluctuations $\widehat{\mathrm{SD}}_J (u_h)$}
	\end{subfigure}
	
	\vspace{0.2cm}
	
	\begin{subfigure}[b]{1\textwidth}
		\includegraphics[width=\textwidth,trim={0 0 0 11mm},clip]{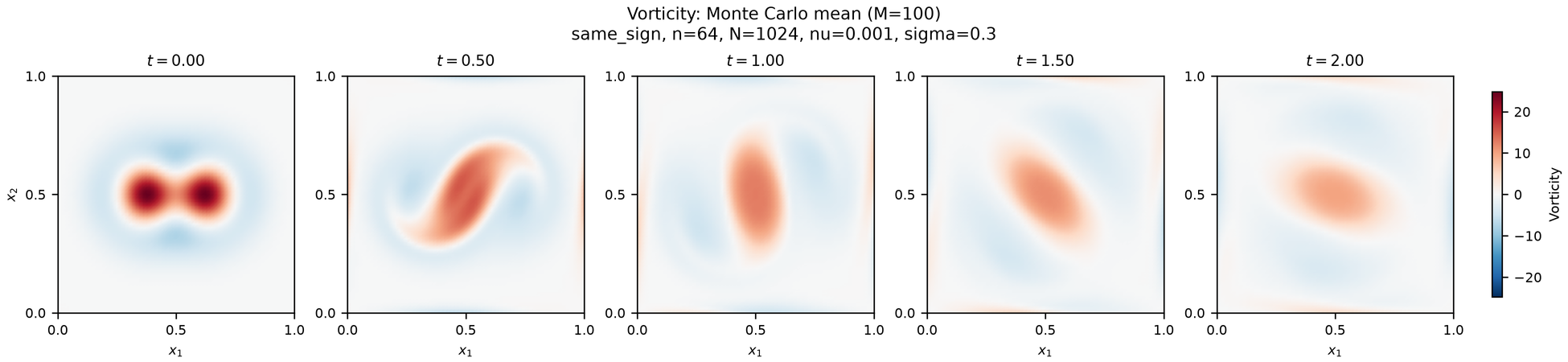}
		\caption{vorticity mean $\widehat{E}_J\omega_h$}
	\end{subfigure}
	
	\vspace{0.2cm}
	
	\begin{subfigure}[b]{1\textwidth}
		\includegraphics[width=\textwidth,trim={0 0 0 11mm},clip]{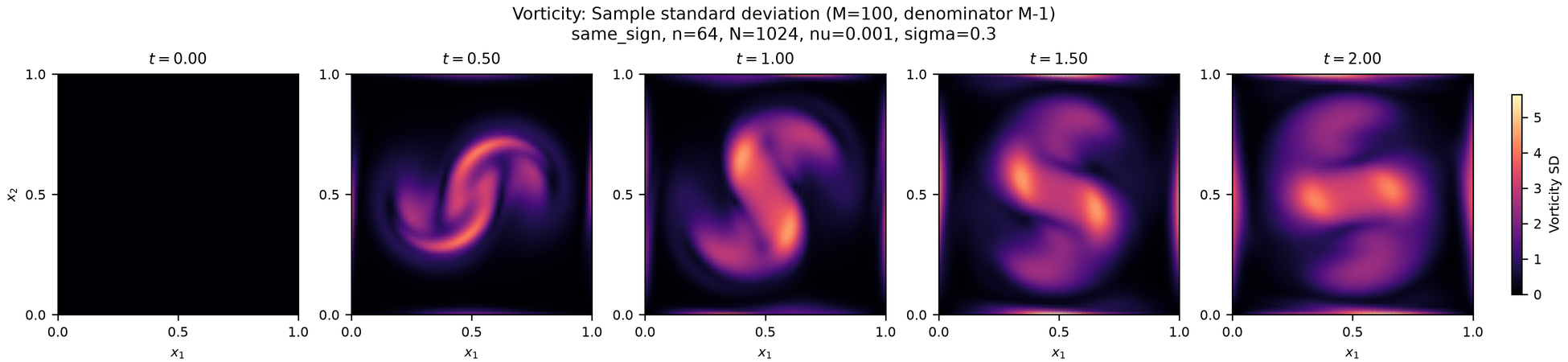}
		\caption{vorticity standard deviation $\widehat{\mathrm{SD}}_J(\omega_h)$}
	\end{subfigure}
	\caption{Approximated mean and fluctuations/standard deviation of the approximation of~\eqref{eq:SNSm} for initial data~\eqref{eq:init2} with $G(u)=\sigma u$, $\sigma=0.3$.}
	\label{fig:exp2_mean}
\end{figure}

\begin{figure}[t]
	\centering
	\begin{subfigure}[b]{1.0\textwidth}
		\includegraphics[width=\textwidth,trim={0 0 0 11mm},clip]{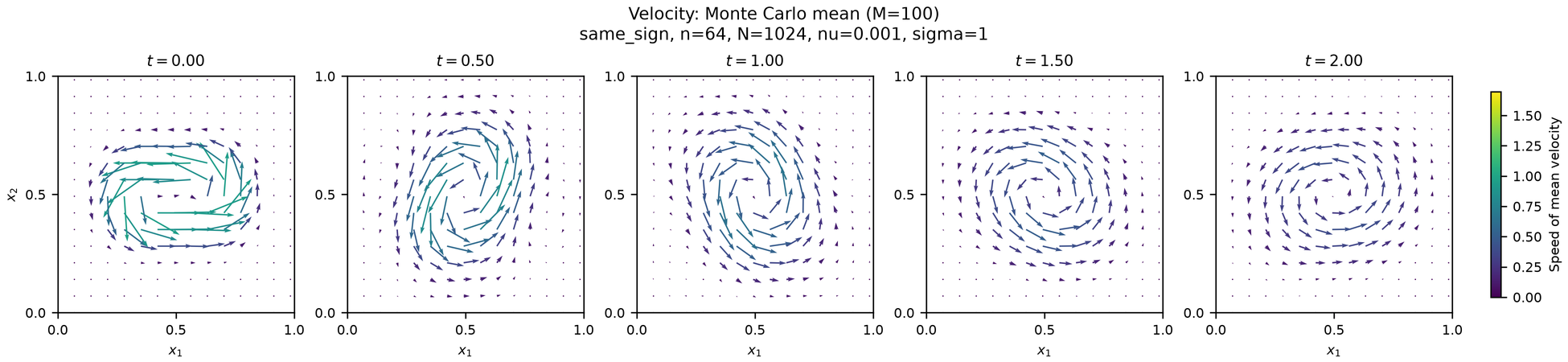}
		\caption{velocity mean $\widehat{E}_J u_h$}
	\end{subfigure}
	
	\vspace{0.2cm}
	
	\begin{subfigure}[b]{1.0\textwidth}
		\includegraphics[width=\textwidth,trim={0 0 0 11mm},clip]{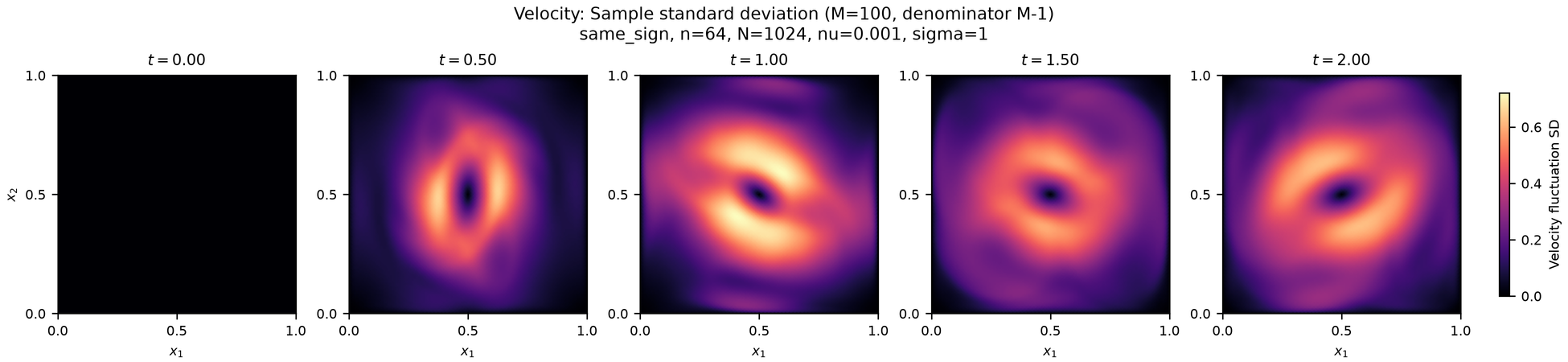}
		\caption{velocity fluctuations $\widehat{\mathrm{SD}}_J (u_h)$}
	\end{subfigure}
	
	\vspace{0.2cm}
	
	\begin{subfigure}[b]{1\textwidth}
		\includegraphics[width=\textwidth,trim={0 0 0 11mm},clip]{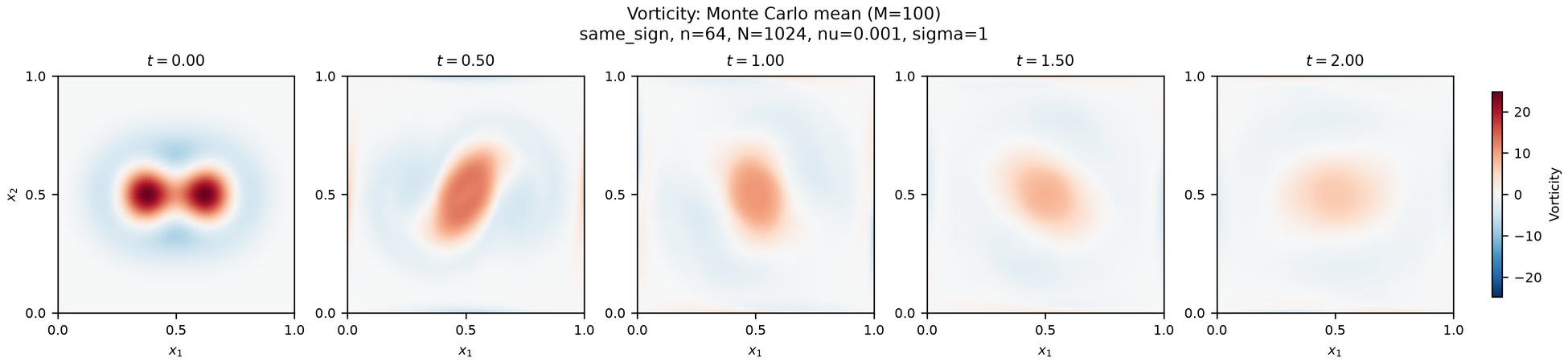}
		\caption{vorticity mean $\widehat{E}_J\omega_h$}
	\end{subfigure}
	
	\vspace{0.2cm}
	
	\begin{subfigure}[b]{1\textwidth}
		\includegraphics[width=\textwidth,trim={0 0 0 11mm},clip]{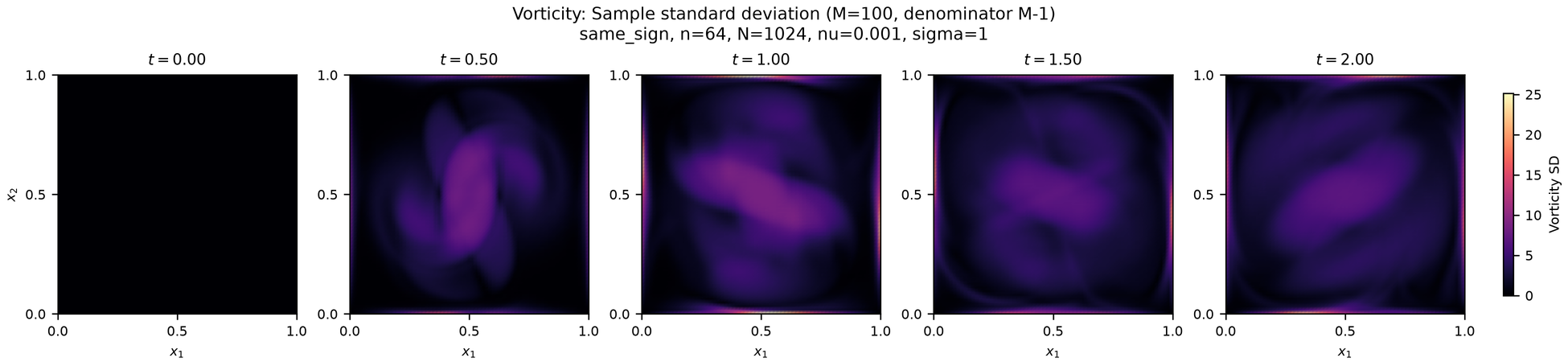}
		\caption{vorticity standard deviation $\widehat{\mathrm{SD}}_J(\omega_h)$}
	\end{subfigure}
	
	\caption{Approximated mean and fluctuations/standard deviation of the approximation of~\eqref{eq:SNSm} for initial data~\eqref{eq:init2} with $G(u)=\sigma u$, $\sigma=1$.}
	\label{fig:exp2_meansigma1}
\end{figure}

The estimated mean and standard deviation of the vorticity
\begin{equation}
	\label{eq:meanvorticity}
	 \widehat{E}_J \omega_h= \frac{1}{J}\sum_{j=1}^J \omega_h^{(j)},\quad \widehat{\mathrm{SD}}_J(\omega_h) = \sqrt{\frac{1}{J-1}\sum_{j=1}^J(\omega_h^{(j)}- \widehat{E}_J \omega_h )^2}
\end{equation}
where $\omega_h^{(j)}$, $j=1,\dots, J$ are independent samples
and the estimated mean and fluctuations of the velocity
\begin{equation}
	\label{eq:meanvelocity}
	\widehat{E}_J u_h = \frac{1}{J}\sum_{j=1}^J u_h^{(j)},\quad \widehat{\mathrm{SD}}_J(u_h) = \sqrt{\frac{1}{J-1}\sum_{j=1}^J |u_h^{(j)}- \widehat{E}_J u_h |^2} = \sqrt{\widehat{\mathrm{SD}}_J(u_{h,1})^2 + \widehat{\mathrm{SD}}_J(u_{h,2})^2 }
\end{equation}
 for $J=100$ are shown in Figures~\ref{fig:exp2_mean} and~\ref{fig:exp2_meansigma1} at times $t=0, 0.5,1,1.5,2$. 
We observe that the two circular regions of positive vorticity rotate around one another and over time form an elongated central structure. The structure appears to rotate and become broader and weaker over time. The empirical mean shows the same behavior but the finer features visible in individual realizations are less pronounced under the averaging. The velocity fluctuations appear concentrated around the areas of fastest rotation whereas the vorticity standard deviation is higher around the vortex cores. 

As we increase the noise amplitude to $\sigma=1$, we observe more pronounced fluctuations, especially near the boundary. The rotating structure remains visible in the empirical mean. Hence, it appears that the stronger multiplicative forcing affects both the magnitude and the spatial distribution of the fluctuations.
\subsection{Rough vortex sheet}
\begin{figure}[t]
	\centering
	\begin{subfigure}[b]{1.0\textwidth}
		\includegraphics[width=\textwidth,trim={0 0 0 11mm},clip]{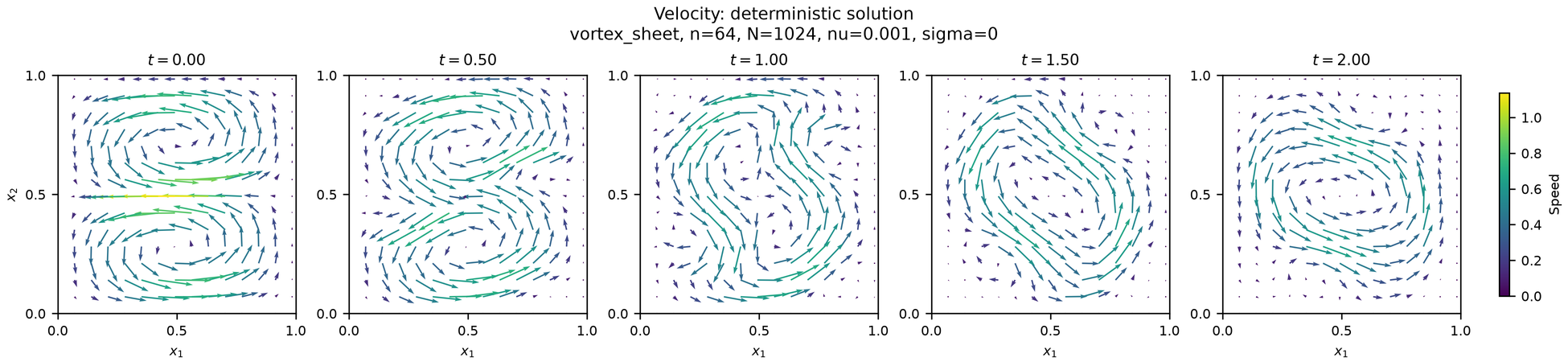}
		\caption{velocity $u_h$}
	\end{subfigure}
	\vspace{0.2cm}
	
	\begin{subfigure}[b]{1\textwidth}
		\includegraphics[width=\textwidth,trim={0 0 0 11mm},clip]{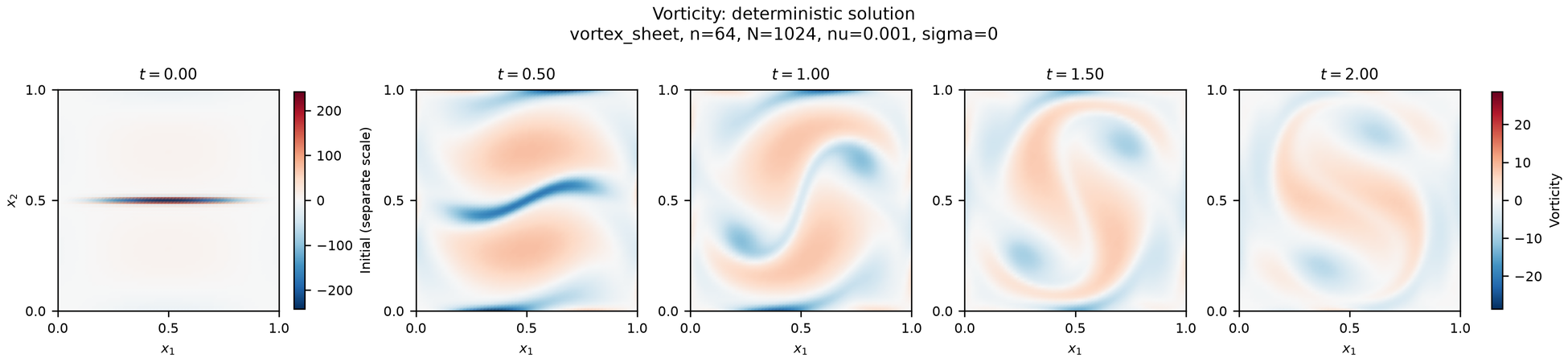}
		\caption{vorticity $\omega_h$}
	\end{subfigure}
	\caption{Approximation of the deterministic solution of~\eqref{eq:SNSm} for initial data~\eqref{eq:init3} ($G(u)=0$).}
	\label{fig:exp3_deterministic}
\end{figure}

\begin{figure}[t]
	\centering
	\begin{subfigure}[b]{1.0\textwidth}
		\includegraphics[width=\textwidth,trim={0 0 0 11mm},clip]{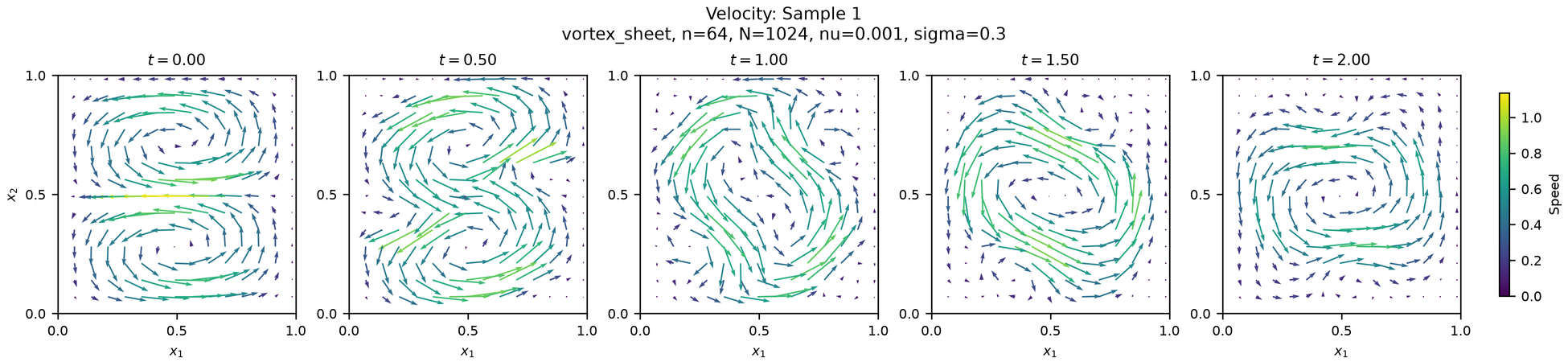}
		\caption{velocity $u_h$}
	\end{subfigure}
	\vspace{0.2cm}
	
	\begin{subfigure}[b]{1\textwidth}
		\includegraphics[width=\textwidth,trim={0 0 0 11mm},clip]{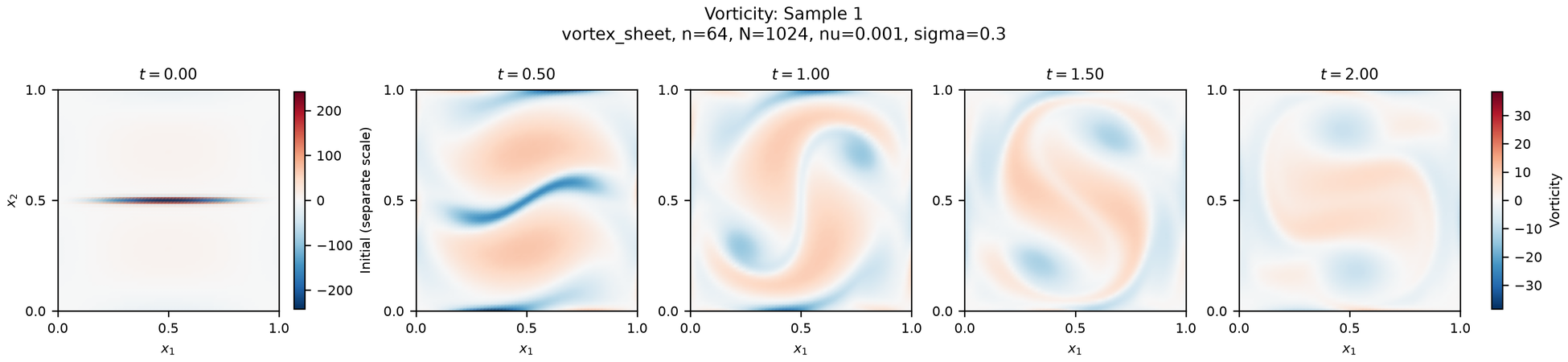}
		\caption{vorticity $\omega_h$}
	\end{subfigure}
	\caption{Sample of the approximation of~\eqref{eq:SNSm} for initial data~\eqref{eq:init3} with $G(u)=\sigma u$, $\sigma=0.3$.}
	\label{fig:exp3_sample}
\end{figure}

\begin{figure}[t]
	\centering
	\begin{subfigure}[b]{1.0\textwidth}
		\includegraphics[width=\textwidth,trim={0 0 0 11mm},clip]{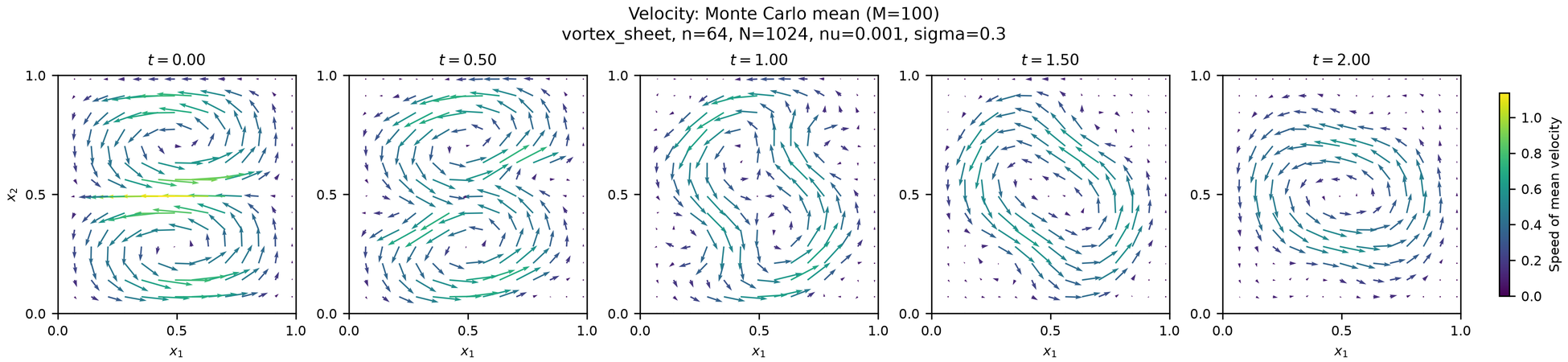}
		\caption{velocity mean $\widehat{E}_J u_h$}
	\end{subfigure}
	
	\vspace{0.2cm}
	
	\begin{subfigure}[b]{1.0\textwidth}
		\includegraphics[width=\textwidth,trim={0 0 0 11mm},clip]{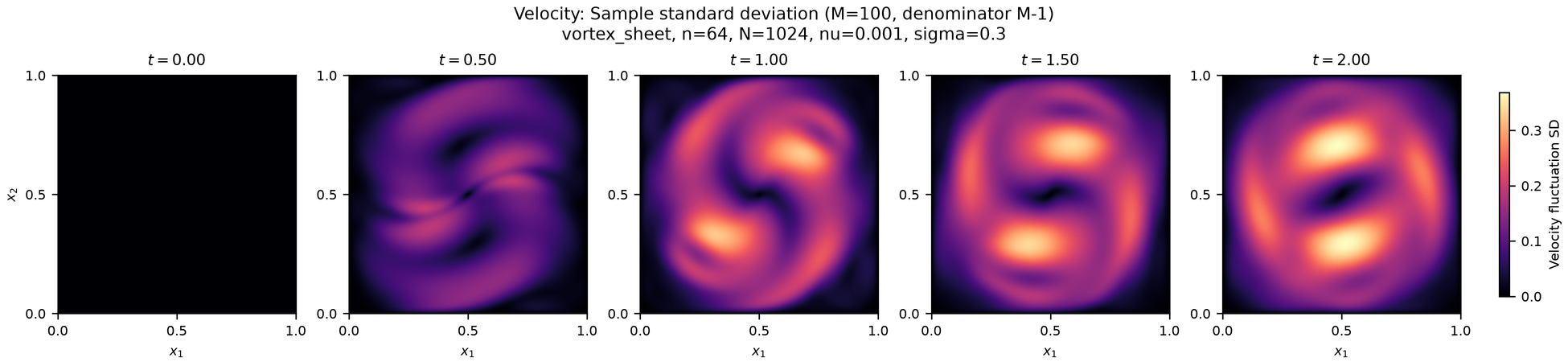}
		\caption{velocity fluctuations $\widehat{\mathrm{SD}}_J (u_h)$}
	\end{subfigure}

	\vspace{0.2cm}
	
	\begin{subfigure}[b]{1\textwidth}
		\includegraphics[width=\textwidth,trim={0 0 0 11mm},clip]{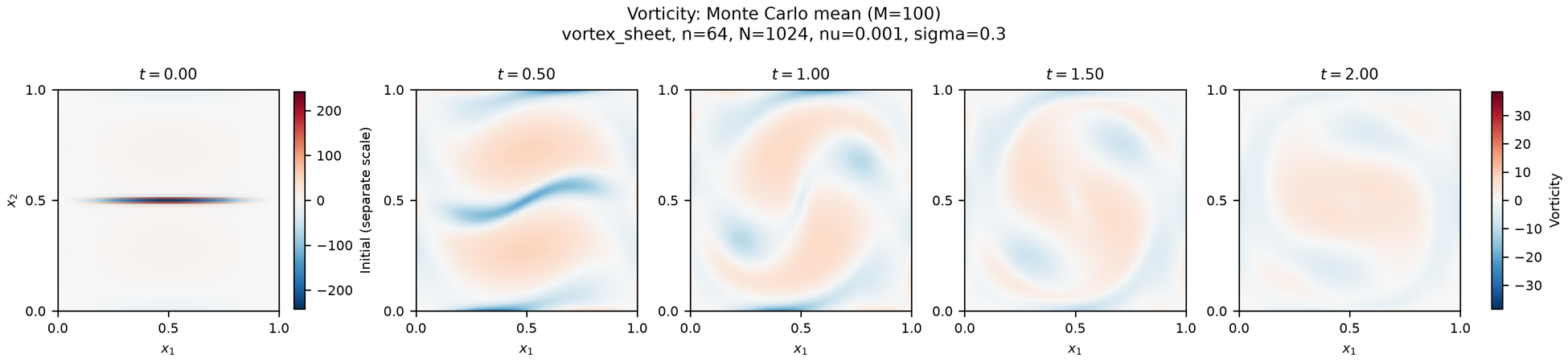}
		\caption{vorticity mean $\widehat{E}_J\omega_h$}
	\end{subfigure}
	
	\vspace{0.2cm}
	
	\begin{subfigure}[b]{1\textwidth}
		\includegraphics[width=\textwidth,trim={0 0 0 11mm},clip]{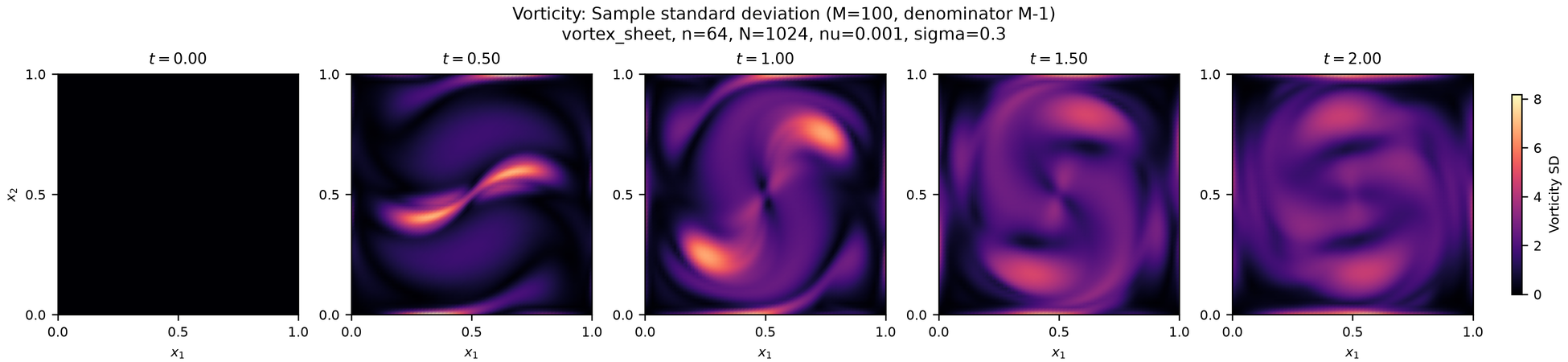}
		\caption{vorticity standard deviation $\widehat{\mathrm{SD}}_J(\omega_h)$}
	\end{subfigure}
	\caption{Approximated mean and standard deviation/fluctuations of the approximation of~\eqref{eq:SNSm} for initial data~\eqref{eq:init3} with $G(u)=\sigma u$, $\sigma=0.3$.}
	\label{fig:exp3_mean}
\end{figure}

Finally, we test the scheme on a problem with less regular initial data, a vortex sheet adapted to the homogeneous Dirichlet boundary conditions. Specifically, we take
\begin{equation}\label{eq:init3}
u_0=\bigl(\partial_{x_2}\psi_0,-\partial_{x_1}\psi_0\bigr),\quad		\psi_0(x_1,x_2)=	256 x_1^2(1-x_1)^2x_2^2(1-x_2)^2	\left|x_2-\tfrac12\right|
\end{equation}
with $n=64$, $N=1024$, $T=2$, and $\nu = 0.001$. The deterministic solution ($G(u)=0$) is shown in Figure~\ref{fig:exp3_deterministic}. Then we add multiplicative noise of the same form as in the previous experiment $G(u) dW = \sigma u d\beta$, with $\sigma=0.3$. A single sample is shown in Figure~\ref{fig:exp3_sample} and the estimated mean and standard deviation and velocity fluctuations for $J=100$ samples are shown in Figure~\ref{fig:exp3_mean}.


In contrast to the previous examples, the initial velocity is only in $L^2_{\Div}(\dom)$ instead of $H^1(\dom)$. This experiment therefore demonstrates the performance of the scheme for low regularity initial data for which our convergence analysis still holds. We observe that the initial narrow sheet deforms into a rotating curved structure that eventually spreads over a larger region. The empirical mean shows similar behavior, whereas the vorticity standard deviation is initially concentrated along the curved sheet and eventually spreads out. The initial vorticity contains a singular contribution supported on $x_2=1/2$. Its approximation therefore exhibits mesh-dependent peaks near the line of discontinuity of the initial velocity field that get steeper as the mesh is refined. For this reason we use a separate color scale for the initial snapshot.

\section*{AI usage statement}
Claude (Anthropic) and Codex (OpenAI) assisted with proofreading, literature search and identifying potential errors and inconsistencies in the notation and mathematical arguments of this manuscript. Furthermore, the code for the numerical experiments was primarily written by Codex, with a few improvements by Claude Code. The author independently reviewed and verified the suggestions and code and takes full responsibility for the manuscript, code and the reported results. The key mathematical ideas and proofs were developed by the author.

\appendix
\section{Auxiliary technical results}
We recall from~\cite{Benavides2026}:
\begin{restatable}{lemma}{rellich}
	\label{lem:Rellich2}
	For any $\eta>0$ there exists $L\in \N$ and $\phi_1,\dots, \phi_L\in L^2(\dom)$ such that for any $v\in L^2(\dom)$, we have
	\begin{equation}
		\label{eq:rellich}
		\norm{v}_{H^{-1}(\dom)}\leq \sqrt{\sum_{k=1}^L (v,\phi_k)^2} + \eta \norm{v}_{L^2(\dom)}.
	\end{equation}	
\end{restatable}
We use this to prove the generalization in Lemma~\ref{lem:Rellich3}, restated here for convenience:
\rellichsmooth*

\begin{proof}
	Given $\zeta>0$, we apply Lemma~\ref{lem:Rellich2} with $\eta=\zeta/2$, to obtain $L\in \N$ and $\widetilde{\phi}_1,\dots, \widetilde{\phi}_L\in L^2(\dom)$ such that for any $v\in L^2(\dom)$, we have
	\begin{equation}\label{eq:prelim}
		\norm{v}_{H^{-1}(\dom)}\leq \sqrt{\sum_{k=1}^L (v,\widetilde{\phi}_k)^2} + \frac{\zeta}{2} \norm{v}_{L^2(\dom)}.
	\end{equation}	
		$C^\infty_{c,\Div}(\dom)$ is dense in $L^2_{\Div}(\dom)$ and $\Grad C^\infty(\overline{\dom})$ is dense in $\Grad H^1(\dom) = \{u\in L^2(\dom),\, u=\Grad p,\, p\in H^1(\dom)\}$ which together span $L^2(\dom)$, i.e., $\overline{\Grad H^1(\dom)}^{L^2}  = L^2_{\Div}(\dom)^\perp$ (cf. Theorem 1.4 in Chapter 1 of~\cite{Temam1977}). Thus for each given $\widetilde{\phi}_k$, $k=1,\dots, L$ we can find $\phi_k\in C^\infty_{c,\Div}(\dom)\oplus \Grad C^\infty(\overline{\dom})$ such that 
		\begin{equation*}
			\norm{\widetilde{\phi}_k-\phi_k}_{L^2(\dom)}\leq \frac{\zeta}{2\sqrt{L}}.
		\end{equation*}
		So we have
		\begin{equation*}
			|(v,\widetilde{\phi}_k)| \leq |(v,\phi_k)| + \norm{\phi_k - \widetilde{\phi}_k}_{L^2}\norm{v}_{L^2} \leq |(v,\phi_k)| + \frac{\zeta}{2\sqrt{L}}\norm{v}_{L^2} 
		\end{equation*}
		Using this in~\eqref{eq:prelim}, we get
		\begin{align*}
			\norm{v}_{H^{-1}(\dom)}&\leq \sqrt{\sum_{k=1}^L (v,{\phi}_k)^2} +\sqrt{\sum_{k=1}^L \norm{v}_{L^2}^2\frac{\zeta^2}{4 L}}+ \frac{\zeta}{2} \norm{v}_{L^2(\dom)}\\
				& =  \sqrt{\sum_{k=1}^L (v,{\phi}_k)^2} +\frac{\zeta}{2}\norm{v}_{L^2}+ \frac{\zeta}{2} \norm{v}_{L^2(\dom)},
		\end{align*}
		which proves the result.
\end{proof}
\bibliographystyle{abbrv}
\bibliography{SNSprojectionbib}

\end{document}